\magnification\magstep1


\input epsf


\input amssym.def
\input amssym.tex


\newcount\secno
\newcount\subsecno
\newcount\subno
\newcount\subsubno
\newcount\refno
\newbox\numbox
\newif\ifcounting
\newif\ifproclaiming
\newif\ifdemoing
\newif\ifshowingtags
\newif\ifdeepsection
\newdimen\partindent
\newdimen\preheadskip
\newdimen\postheadskip
\newread\oldtags
\newwrite\newtags

\let\nonumbers=\countingfalse
\countingtrue

\showingtagsfalse
\deepsectionfalse

\global\partindent = 44pt
\global\preheadskip = 18pt
\global\postheadskip = 0pt
\global\secno=0
\global\subsecno=0
\global\subno=0
\global\subsubno=0
\global\refno=0
\xdef\secstring{\the\secno}
\xdef\currentno{\secstring}


\font\commentfont=cmbx6


\let\Afont=\bf
\let\Bfont=\it
\let\titlefont=\bigbold
\let\authorfont=\bf
\let\addressfont=\small
\let\titleskip=\relax 
\let\sectionsymbol=\relax 
\let\headstyle=\leftline 
\parindent=0pt
\parskip=4pt plus 1pt
\baselineskip 14pt


\def\subnumber{
     \global\advance\subno by 1
     \global\subsubno=0
     \xdef\currentno{\secstring.\the\subno}\currentno}

\def\subsubnumber{
     \global\advance\subsubno by 1
     \xdef\currentno{\secstring.\the\subno.\the\subsubno}\currentno}

\outer\def\introduction{
     \vskip\preheadskip
     \message{Introduction}
     \headstyle{\bf Introduction}
     \vskip\postheadskip\nobreak\noindent}

\outer\def\section#1{
     \deepsectionfalse
     \global\advance\secno by 1
     \message{Section \the\secno}
     \global\subsecno=0
     \global\subno=0
     \global\subsubno=0
     \goodbreak
     \vskip\preheadskip\par
     \vskip0pt plus .1\vsize\penalty-50
     \vskip0pt plus-.1\vsize
     \headstyle{\Afont{\sectionsymbol\the\secno.~~#1}}
     \xdef\secstring{\the\secno}\xdef\currentno{\secstring}
     \vskip\postheadskip\par}

\outer\def\deepsection#1{
     \deepsectiontrue
     \global\advance\secno by 1
     \message{Section \the\secno}
     \global\subsecno=0
     \global\subno=0
     \global\subsubno=0
     \goodbreak
     \vskip\preheadskip\par
     \vskip0pt plus .1\vsize\penalty-50
     \vskip0pt plus-.1\vsize
     \headstyle{\Afont{\sectionsymbol\the\secno.~~#1}}
     \xdef\secstring{\the\secno.\the\subsecno}\xdef\currentno{\the\secno}
     \vskip\postheadskip\par}

\outer\def\subsection#1{
     \ifdeepsection
       \global\advance\subsecno by 1
       \global\subno=0
       \global\subsubno=0
       \message{Subsection \the\secno.\the\subsecno}
       \goodbreak
       \vskip\preheadskip\par
       \vskip0pt plus .1\vsize\penalty-50
       \vskip0pt plus-.1\vsize
       \headstyle{\Afont{\sectionsymbol\the\secno.\the\subsecno.~~#1}}
       \xdef\secstring{\the\secno.\the\subsecno}\xdef\currentno{\secstring}
       \vskip\postheadskip\par
     \else
       \errmessage{Only deep sections can have subsections.}
     \fi}

\outer\def\paragraph{
     \vskip\preheadskip\par
     \vskip0pt plus .05\vsize\penalty-50
     \vskip0pt plus-.05\vsize
     {\Afont\subnumber.}\spacefactor=500\enskip\ignorespaces}

\outer\def\appendix#1#2{
     \message{Appendix #1}
     \global\subno=0
     \goodbreak
     \vskip\preheadskip
     \vskip0pt plus .1\vsize\penalty-50
     \vskip0pt plus-.1\vsize
     \headstyle{\Afont{Appendix #1.~~#2}}
     \xdef\secstring{\hbox{#1}}\xdef\currentno{\secstring}
     \vskip\postheadskip\par}

\outer\def\references{\vskip20pt
     \vskip0pt plus .1\vsize\penalty-50
     \vskip0pt plus-.1\vsize
     \bigskip\vskip\parskip
     \message{References}
     \centerline{\bf Bibliography}
     \bigskip
     \frenchspacing
     \tolerance=500
     }


\edef\t@mp{\catcode`\noexpand\#=\the\catcode`\#}%
    \catcode`\#=12
    \def\h@sh{#}%
    \t@mp

\def\tagdef#1#2{\expandafter\xdef\csname#1\endcsname{#2}}

\def\usetagfile{
  \openin\oldtags=\jobname.tag
  \ifeof\oldtags
    \message{I did not find file \jobname.tag. Cross-references will not appear.}
    \closein\oldtags
  \else
    \message{Loading cross-reference file}
    \closein\oldtags
    \input \jobname.tag
  \fi
  \immediate\openout\newtags=\jobname.tag
}

\def\tag#1{\ifmmode\expandafter\xdef\csname#1\endcsname{(\currentno)}
     \else\expandafter\xdef\csname#1\endcsname{\currentno}\fi%
     \unskip%
     \unskip%
     \immediate\write\newtags{\string\tagdef{#1}{\currentno}}%
     \ifshowingtags{\rm $\langle$#1$\rangle$ }\fi\ignorespaces}

\def\xref#1{\expandafter\ifx\csname#1\endcsname\relax
 ??\else\csname#1\endcsname\unskip\fi}


\def\ref#1{
      \global\advance\refno by 1
      \xdef#1{\the\refno}}

\def\cite#1{\csname#1\endcsname}

\def\key#1{\part{[\csname#1\endcsname]}}


\def\caption#1{
     \medbreak
     \vskip0pt plus .1\vsize\penalty-50
     \vskip0pt plus-.1\vsize
     \leftline{\Bfont #1}
     \medskip}

\def\proclaim#1{\checkends\proclaimingtrue
     \vskip0pt plus .05\vsize\penalty-50
     \vskip0pt plus-.05\vsize
     \ifvmode\medbreak\noindent\fi
     \begingroup\sl
     {\Afont#1\ifcounting\Afont\subnumber\fi.}\enspace\spacefactor=500\ignorespaces}

\def\endproclaim{\par\endgroup
     \ifdim\lastskip<\medskipamount \removelastskip\penalty55\medskip\fi
     \proclaimingfalse\countingtrue}

\def\part#1{\par\noindent\hangindent\partindent
     \hbox to \partindent{\hskip .5\partindent minus .5\partindent
     #1\enspace\hfill}\ignorespaces}

\def\demo#1{\checkends\demoingtrue
     \ifvmode\noindent\fi{\Bfont#1.\enspace}}

\def\enddemo{\demoingfalse\smallskip}

\def\checkends{\ifproclaiming\errmessage{Missing endproclaim}\fi
     \ifdemoing\errmessage{Missing enddemo}\fi}

\def\text#1 {{\qquad\hbox{\rm #1}\qquad}}


\def\title#1\par{\xdef\thetitle{#1}}

\def\author#1\par{\xdef\theauthor{#1}}

\def\support#1\par{\ \footnote{}{#1}}

\def\address{\medskip\begingroup\addressfont
     \parindent=0pt
     \parskip=0pt
     \baselineskip=3pt
     \leftskip=4em
     \obeylines}

\def\endaddress{\endgroup\par\bigskip}

\def\today{\ifcase\month\or January\or February\or March\or April\or May\or
	June\or July\or August\or September\or October\or November\or
	December\fi
	\space\number\day, \number\year}
%


\def\showtitle{
\vglue1truein
\halign to\hsize{\titleskip\titlefont##\hfil\cr\thetitle\cr\hskip\hsize\cr}
\smallskip
\halign to\hsize{\titleskip\authorfont##\hfil\cr\theauthor\cr\hskip\hsize\cr}
\medskip
}

\def\comment{\par\bgroup\commentfont}
\def\endcomment{\egroup\par}


\def\Definition{\proclaim{Definition}\ \rm}
\def\EndDefinition{\endproclaim}
\def\Theorem{\proclaim{Theorem}}
\def\EndTheorem{\endproclaim}

\def\Lemma{\proclaim{Lemma}}
\def\EndLemma{\endproclaim}
\def\Proposition{\proclaim{Proposition}}
\def\EndProposition{\endproclaim}
\def\Corollary{\proclaim{Corollary}}
\def\EndCorollary{\endproclaim}
\def\EndProof{{\hfill$\square$\medskip}}
\def\Proof{{{\it Proof.\enskip}\/}}
\def\NoProof{{\hfill$\square$}}

\input veggie.tag

\def\genus{\mathop{\rm genus}}
\def\int{\mathop{\rm int}}
\def\frontier{\mathop{\rm Fr}}
\def\im{\mathop{\rm im}}
\def\VT{\mathop{\rm VT}}
\def\lcm{\mathop{\rm lcm}}
\def\split(#1,#2){#1_{#2}}
\def\rTo^#1{{\;\buildrel #1 \over {\hbox to 30pt{\rightarrowfill}}\; }}
\def\bdry{\partial}
\def\from{\colon}
\def\homeo{\cong}

\def\lint{\wedge_{\cal L}}
\def\pint{\wedge_{\cal P}}
\def\dotlint{\dot\wedge_{\cal L}}
\def\R{{\Bbb R}}
\def\Z{{\Bbb Z}}
\def\Q{{\Bbb Q}}
\def\C{{\Bbb C}}

\def\v{{\varphi}}
\def\a{{\alpha}}
\def\b{{\beta}}

\def\D{{\Delta}}

\def\t{{\tau}}
\def\u{{\upsilon}}

\title
Characteristic subsurfaces and Dehn filling

\author
Steve Boyer$^1$, Marc Culler$^2$, Peter B. Shalen$^2$ and Xingru Zhang$^3$

\support
\vbox{1. Partially supported by NSERC grant OGP0009446 and FCAR grant ER-68657.}
\vbox{2. Partially supported by NSF grant DMS 0204142.}
\vbox{3. Partially supported by NSF grant DMS 0204428.}

\showtitle

\introduction

Many results in the theory of Dehn surgery (see [\cite{GordonSurvey}])
assert that if $M$ is a compact, orientable, atoroidal, irreducible
$3$-manifold whose boundary is an incompressible torus, and if two
Dehn fillings $M(\alpha)$ and $M(\beta)$ have specified properties,
then the distance $\Delta(\alpha,\beta)$ of the slopes $\alpha$ and
$\beta$ is bounded by a suitable constant.  (The reader is referred to
the body of this paper for the definition of ``slope'' and
``distance'', as well as for precise versions of many definitions,
statements and proofs that are hinted at in this introduction.)

In this paper, we define a closed $3$-manifold to be {\it very small}
if its fundamental group has no non-abelian free subgroup.  The
motivating result of the paper, Corollary
\xref{VerySmallReducible}, asserts that if $M(\alpha)$ is very small
and $M(\beta)$ is a reducible manifold other than $S^2\times S^1$ or
$P^3\# P^3$, then $\Delta(\alpha,\beta)\le5$.  This
follows from the following stronger result which deals with
essential planar surfaces in $M$ which are not semi-fibers
(see \xref{Splittings}).

\nonumbers
\proclaim{Corollary \xref{PlanarDisk}}
Let $M$ be a simple knot manifold and $F\subset M$ an essential planar
surface with boundary slope $\beta$ which is not a semi-fiber.  Let
$\alpha$ be a slope in $\partial M$. If $M(\alpha)$ is very small, or
more generally if $F\subset M\subset M(\alpha)$ is not
$\pi_1$-injective in $M(\alpha)$, then $\D(\a,\b)\le 5.$
\endproclaim

Corollary \xref{HighGenusDisk} gives a qualitatively similar
conclusion when $M(\alpha)$ is very small and $\beta$ is the boundary
slope of an essential bounded surface $F$ of arbitrary genus which is
not a semi-fiber. Here the upper bound for $\Delta(\alpha,\beta)$ is
$15+(20g-15)/m$, where $g$ is the genus of $F$ and $m$ is the number
of components of $\partial F$. Note that, while this result does
provide an upper bound in the case where $F$ is planar, the bound
given by Corollary \xref{PlanarDisk} is much stronger. In a follow-up paper
we will examine the case where $g=1$ and show that in this situation,
the bounds obtained in Corollary \xref{HighGenusDisk}, and those described
in the results mentioned below, can be significantly improved.

The proofs of these results begin with the observation that a
bounded essential surface $F$ with boundary slope $\beta$ may be
regarded as a non-properly embedded surface of negative Euler
characteristic in $M(\alpha)$, and that if $M(\alpha)$ is very
small then the inclusion homomorphism from $\pi_1(F)$ to
$\pi_1(M(\alpha))$ cannot be injective. From this one can deduce
that there is a map of a disk into $M(\alpha)$ which maps the
boundary of the disk into $M-F$ but cannot be homotoped rel
boundary into $M-F$. After normalizing such a map and restricting
it to the inverse image of $M$, one obtains a ``singular surface''
in $M$, ``well-positioned'' with respect to $F$; according to the
precise definitions given in Section \xref{SingularSurfaces}, such
a singular surface is defined by a map $h$ of a surface $S$ into
$M$ having certain properties. In the case we are discussing here,
$S$ is planar, and each component of $\partial S$ is either mapped
into $M-F$ by $h$, or mapped homeomorphically onto a curve in
$\partial M$ of slope $\alpha$.
This last property is expressed by saying
that the singular surface has boundary slope $\alpha$.

The main results of the paper, Theorems \xref{HighGenusTheorem} and
\xref{PlanarTheorem}, give bounds on the distance between two slopes
$\alpha$ and $\beta$ in terms of the data involving an essential
surface $F$ in $M$ which is not a semi-fiber and has boundary slope
$\beta$, and a singular surface which is well-positioned with respect
to $F$ and has boundary slope $\alpha$. Applying this in the case of a
planar singular surface we obtain such results as Corollaries
\xref{HighGenusDisk} and \xref{PlanarDisk}.

By applying Theorems \xref{HighGenusTheorem} and
\xref{PlanarTheorem} to other kinds of singular surfaces, we
obtain different kinds of information about boundary slopes and
Dehn filling.  For instance we prove:

\nonumbers
\proclaim{Corollary \xref{PlanarSeifert}}
Let $M$ be a simple knot manifold and $F\subset M$ an essential planar
surface with boundary slope $\beta$ which is not a semi-fiber. Let
$\alpha$ be a slope in $\partial M$. If $M(\a)$ is a Seifert fibered
space or if there exists a $\pi_1$-injective map from $S^1\times S^1$
to $M$ then
$\D(\a,\b)\le 6.$
\endproclaim

This implies
Corollary \xref{ReducibleSeifert}, which asserts that if
$M(\alpha)$ is a Seifert fibered space and $M(\beta)$ is a
reducible manifold other than $S^2\times S^1$ or $P^3\# P^3$, then
$\Delta(\alpha,\beta)\le6$. Corollary \xref{PlanarSeifert},
like Corollary \xref{PlanarDisk}, has a high-genus analogue:
Corollary \xref{HighGenusSeifert} asserts that if $M(\alpha)$ is a
Seifert fibered space and $\beta$ is the boundary slope of an
essential surface $F$ in $M$ which is not a semi-fiber, then
$\Delta(\alpha,\beta)\le18+(24g-18)/m$, where $g$ is the genus of
$F$ and $m$ is the number of its boundary components. These
results are proved by observing that if $M(\alpha)$ is Seifert
fibered then either it is very small, in which case the
conclusions follow from Corollaries \xref{HighGenusDisk} and
\xref{PlanarDisk}, or it contains a $\pi_1$-injective singular
torus. Such a torus can be used to construct a genus-$1$ singular
surface in $M$ having boundary slope $\alpha$, to which Theorems
\xref{HighGenusTheorem} and \xref{PlanarTheorem} can be applied.

Still another type of application of Theorems
\xref{HighGenusTheorem} and \xref{PlanarTheorem} can be obtained
by observing that an essential surface in $M$ with boundary slope
$\alpha$ is special case of a singular surface with boundary slope
$\alpha$. This leads to upper bounds for the distance between
boundary slopes of two essential surfaces (not both semi-fibers)
in terms of the genera and numbers of boundary components of the
surfaces. Such bounds are given in Corollary
\xref{CameronCorollary} in the general case, and 
Corollary \xref{GordonLitherland}, which recovers a result of Gordon
and Litherland [\cite{GLi}, Proposition 6.1] under the additional hypothesis
that one of the surfaces is planar and is not a semi-fiber.  Corollary
\xref{CameronCorollary}, which is qualitatively similar to an unpublished
result due to Cameron Gordon, strengthens a result due to Torisu
[\cite{To}], but in turn has been strengthened slightly by Agol [\cite{Agol}],
using results of Cao and Meyerhoff [\cite{CaoMeyerhoff}].

The constructions of singular surfaces that are needed to pass from
Theorems \xref{HighGenusTheorem} and \xref{PlanarTheorem} to their various
corollaries are given in detail in Section \xref{SingularSurfaces}.

The statements of Theorems \xref{HighGenusTheorem} and
\xref{PlanarTheorem} involve an essential surface $F$ which is not a semi-fiber.
We shall sketch the proofs under the simplifying assumption that $M-F$
has two components, whose closures we shall denote by $M_F^+$ and
$M_F^-$.  The first step in the proofs, which is carried out in
Section \xref{ReducedHomotopies}, involves graph-theoretical
arguments. Suppose that we are given an essential surface $F\subset M$
with boundary slope $\beta$ and a singular surface which is
well-positioned with respect to $F$ and has boundary slope
$\alpha$. Such a singular surface is defined by a certain kind of map
$h$ of a compact $2$-manifold $S$ into $M$. We obtain a surface $\hat
S$ from $S$ by identifying certain components of $\partial S$ to
points, and the images of the arc components of $h^{-1}(S)$ are the
edges of a graph $G\subset \hat S$. Each vertex of $G$ has valence
$m\Delta(\alpha,\beta)$, where $m$ is the number of boundary
components of $F$. By using certain non-degeneracy properties of $G$
we find a family of parallel edges in $G$ whose size is bounded below
in terms of topological data about $F$ and the valence
$m\Delta(\alpha,\beta)$.

A parallel family of edges in $G$ gives $3$-dimensional information
about how the essential surface $F$ sits in $M$.  While the edges of
$G$ do not map to properly embedded arcs in $F$, each edge of $G$ does
give rise to an essential path in $(F,\partial F)$, which can be
extended to a map of a ``pair of glasses'' (see Figure
\xref{GlassesFig}) into $F$ that maps the rims homeomorphically to
components of $\partial F$. A parallel family of $k+1$ edges in $G$
defines a sequence of $k+1$ such ``singular pairs of glasses'' and $k$
essential homotopies in $M_F^+$ and $M_F^-$ between the successive
pairs of glasses in the sequence. (Under the homotopies, the images of
the rims of the glasses stay in $\partial M$.)  Furthermore, these
homotopies alternate strictly between homotopies in $M_F^+$ and
homotopies in $M_F^-$.  According to the precise definition given in
Section \xref{ReducedHomotopies} such a sequence of homotopies determines a
{\it reduced homotopy of length} $k$. The graph-theoretical arguments
that we have sketched here are used in Section
\xref{ReducedHomotopies} to show that upper bounds for the length of a
reduced homotopy of singular pairs of glasses in $M$ imply theorems of
the type of \xref{HighGenusTheorem} and \xref{PlanarTheorem}. The rest
of the paper is devoted to obtaining such bounds for the lengths of
reduced homotopies, in the more general context of a map of a
polyhedron into $M$ which is ``large'' in the sense that the induced
homomorphism of fundamental groups has a non-abelian image.

A reduced homotopy of length $1$ is by definition an essential
homotopy in $M_F^+$ or $M_F^-$ whose time-$0$ and time-$1$ maps are
maps of the domain into $F$. We are interested in reduced homotopies
whose time-$0$ maps are large.  Such homotopies can be understood in
terms of the characteristic submanifold theory ([\cite{JacoShalen}],
[\cite{Johannson}]). This theory provides a (possibly disconnected)
$2$-manifold $\Phi_\pm\subset M_F^\pm$, which is ``large'' in the
sense that the fundamental group of each component of $\Phi_\pm$ is
non-abelian and maps injectively into $\pi_1(F)$.  Any large map of a
polyhedron into $F$ which is the time-$0$ map of an essential homotopy
in $M_F^\pm$ is homotopic in $F$ to a map into
$\Phi_\pm$. Furthermore, the identity map of $\Phi_\pm$ is itself the
time-$0$ map of an essential homotopy in $M_F^\pm$.

In Section \xref{GeneralBounds} we generalize this to reduced
homotopies of length $k$: we define large $2$-dimensional submanifolds
$\Phi_k^\pm$ of $F$ for $k=1,2,\ldots$. If a large map $f$ of a
polyhedron into $F$ is the time-$0$ map of a length-$k$ reduced
homotopy in $M$ which ``begins'' in $M_F^+$ (or $M_F^-$), then $f$ is
homotopic in $F$ to a map into $\Phi_k^+$ (respectively
$\Phi_k^-$). Furthermore, the identity map of $\Phi_k^\pm$ is itself
the time-$0$ map of a reduced homotopy in $M$ which begins in
$M_F^\pm$. The $\Phi_k^\pm$ for $k>1$ are defined inductively, using
the notion of an essential intersection of subsurfaces of $F$, which
is presented in [\cite{Jaco}]. In Section \xref{EssentialIntersections} we
give a self-contained account of a version of the theory of essential
intersections that is adapted to the study of large subsurfaces.

The $\Phi_k^\pm$ give a natural tool for bounding the length of a
reduced homotopy: if $n$ is a positive integer such that
$\Phi_+^n=\Phi_-^n=\emptyset$, it is clear that any reduced homotopy
with large time-$0$ map has length $<n$.  (Of course there can be no
such bound in the case that $F$ is a semi-fiber.)

The $\Phi_k^\pm$ may be taken to be nested:
$\Phi_1^+\supset\Phi_2^+\supset\dots$, and similarly for the
$\Phi_k^-$. A crucial step in the argument is provided by Proposition
\xref{StrictContainment}, which asserts that when $F$ is not a
semi-fiber, and if $\Phi_k^+$ (say) is non-empty for a given $k$, then
$\Phi_{k+2}^+$ is not isotopic to $\Phi_k^+$. Hence in the sequence
$\Phi_1^+\supset\Phi_3^+\supset\Phi_5^+\supset\dots$, the successive
subsurfaces are always non-isotopic until one of them becomes
empty. This means that to bound the length of a reduced homotopy whose
time-$0$ map is large, it suffices to bound the length of a nested
sequence of subsurfaces of $F$ in which successive subsurfaces are
non-isotopic. This is a matter of elementary surface topology, and the
bound can be improved by a factor of $2$ using Corollary
\xref{EvenChar}, which asserts that for odd $k$ the $\Phi_k^\pm$ all
have even Euler characteristic.  This leads to Theorem
\xref{LengthBound}, which gives a bound of $8g+3m-8$ for the length of
a reduced homotopy with a large time-$0$ map, where $g$ is the genus
of $F$ and $m$, as above, denotes the number of its boundary
components.

While Theorem \xref{LengthBound} is significant for general large
maps, it is far from optimal for the case of singular pairs of glasses
arising from essential paths. In Section \xref{DottedPhis} we
introduce a variant of $\Phi_k^\pm$ which we denote by
$\dot\Phi_k^\pm$; it is simply the union of the ``outer components of
$\Phi_k^\pm$, i.e. those components which have at least one boundary
component which is homotopic to a component of $\partial F$. If a
singular pair of glasses $f$ arising from an essential path in $F$ is
the time-$0$ map of a length-$k$ reduced homotopy in $M$ which begins
in $M_F^+$ (or $M_F^-$), then $f$ is homotopic in $F$ to a map into
$\dot\Phi_k^+$ (respectively $\dot\Phi_k^-$). In Section
\xref{DottedPhis}, under the assumption that there exists a reduced
homotopy of length $m$ whose time-$0$ map is an essential path, we
establish analogues for the $\dot\Phi_k^\pm$ of all the properties of
the $\Phi_k^\pm$ that are established in Section 4, including
Corollary \xref{DottedEvenChar} and Proposition
\xref{DottedStrictContainment} which are the analogues of Corollary
\xref{EvenChar} and Proposition \xref{StrictContainment}. By
definition the $\dot\Phi_k^\pm$ have the additional property that they
are outer subsurfaces, in the sense that all their components are
outer components. Hence to bound the length of a reduced homotopy
whose time-$0$ map is an essential path, it suffices to bound the
length of a nested sequence of outer subsurfaces of $F$ in which
successive subsurfaces are non-isotopic. The restriction to outer
subsurfaces turns out to improve the bound, almost by another factor
of $2$. The upshot is Theorem \xref{DottedLengthBound}, which gives a bound of
$4g+3m-4$ for the length of a reduced homotopy whose time-$0$ map is
an essential path, where $g$ and $m$ are defined as above. Theorem
\xref{HighGenusTheorem} is proved by combining Theorem \xref{DottedLengthBound}
with the results of Section 2.

In order to prove Theorem \xref{PlanarTheorem} we need to improve
the conclusion of Theorem \xref{DottedLengthBound} in the special
case where $F$ is planar, i.e. $g=0$.
 
\nonumbers
\proclaim{Theorem \xref{TightLengthBound}}
Let $F$ be an essential planar surface in a simple knot manifold $M$.
Suppose that $F$ is not a semi-fiber.  Set $m =
|\bdry F|$ and let $H$ be any reduced homotopy in the pair $(M,F)$ such
that $H_0$ is an essential path in $F$ and $H_t(\partial I)\subset
\partial M$ for each $t\in I$. Then the length of $H$ is at
most $m-1$.
\endproclaim

The proof formally proceeds by contradiction, beginning with the
assumption that there is a length-$m$ reduced homotopy whose time-$0$
map is an essential path, and applying some machinery that is set up
in Section \xref{TightSurfaces}. We shall sketch this machinery in the
case where $F$ is planar, although many of the results of Section
\xref{TightSurfaces} are stated more generally. We introduce yet
another variant of $\Phi_k^\pm$ which we denote by
$\breve\Phi_k^\pm$. It differs from $\dot\Phi_k^\pm$ in that it
contains $\partial F$, but may have annular components.
Under the assumption that there is a length-$m$ reduced homotopy
whose time-$0$ map is an essential path, we again obtain analogues
for the $\breve\Phi_k^\pm$ of the properties that are established
in the preceding sections for the $\Phi_k^\pm$ and the
$\dot\Phi_k^\pm$, even though the surfaces $\breve\Phi_k^\pm$ need
not be large. In Proposition \xref{SmileyStrictContainment}, which
is the analogue of Propositions \xref{StrictContainment} and
\xref{DottedStrictContainment}, the analogue of the conditions
$\Phi_k^\pm=\emptyset$ or $\dot\Phi_k^\pm=\emptyset$ is that
$\breve\Phi_k^\pm$ is a regular neighborhood of $\partial F$,
which is in fact equivalent to the condition
$\dot\Phi_k^\pm=\emptyset$.

The planarity of $F$ implies that if $F$ is not a semi-fiber then some
component of $\breve\Phi_1^\pm$ is {\it tight} in the sense that its
frontier in $F$ is a single simple closed curve. We define the {\it
size} of a tight component of $\breve\Phi_k^\pm$ to be the number of
components of $\partial F$ that it contains.  We call a component of
$\breve\Phi_1^+$ or $\breve\Phi_1^-$ {\it very tight} if its size is
at most the minimum size of any tight component of $\breve\Phi_1^+$ or
$\breve\Phi_1^-$. We may assume by symmetry that $\breve\Phi_1^+$ has
a very tight component.  Lemma \xref{NumVTEven}, the proof of which is
based on the same ideas as that of Corollary \xref{EvenChar}, implies
that the number $|VT(\breve\Phi_k^+)|$ of very tight components of
$\breve\Phi_k^+$ is always even. A key step in proving Theorem
\xref{TightLengthBound} is Lemma \xref{NumVTIncreasesStrictly},
which implies that if $F$ is not a semi-fiber then increasing $k$
by $2$ always strictly increases $|VT(\breve\Phi_k^+)|$, unless
$\breve\Phi_k^+$ is already a regular neighborhood of $\partial
F$. From this it is not hard to deduce (cf. Proposition
\xref{BigSmile}) that $\breve\Phi_{m-1}^+$ is a regular
neighborhood of $\partial F$; this in turn easily implies that there
is no length-$m$ reduced homotopy whose time-$0$ map is an essential
path, a contradiction which completes the proof of Theorem
\xref{TightLengthBound}.

In Section \xref{SeifertSurgeries} we investigate further the
situation where $M$ is a simple knot manifold and $\alpha$ and $\beta$
are slopes such that $M(\alpha)$ is Seifert fibered while $M(\beta)$
is reducible.  In Proposition \xref{ThreeToSix} we establish
restrictions on which Seifert fibered spaces can arise in
this situation when $\Delta(\alpha,\beta) > 3$.  The proof begins by
applying Corollary \xref{ReducibleSeifert} to deduce that
$\Delta(\alpha,\beta)$ is equal to 4, 5 or 6.  We then use the
$PSL_2(\C)$ character variety together with some observations from
algebraic number theory to describe the Seifert fibered structure on
$M(\alpha)$. For instance we prove:

\nonumbers
\proclaim{Corollary \xref{ThreeToSixCorollary}}
Let $M$ be a simple knot manifold and fix slopes $\alpha$ and $\beta$
on $\partial M$.  If $M(\beta)$ is reducible, though not $S^1 \times
S^2$ or $P^3 \# P^3$, and $M(\alpha)$ is a Seifert fibered space, then
$\Delta(\alpha, \beta) \leq 5$ unless perhaps $M(\beta) \cong P^3 \#
L(p,q)$ and $M(\alpha)$ is a small Seifert manifold with base orbifold
$S^2(a,b,c)$ where $(a,b,c)$ is a hyperbolic triple and $6$ divides
$\lcm(a,b,c)$.
\endproclaim

Similar methods lead to restrictions, in Proposition
\xref{SixToTen}, of the possible Seifert fibrations of
$M(\alpha)$ when $M(\beta)$ is a Seifert fibered space containing an
(embedded) incompressible torus and $5 < \Delta(\alpha,\beta)\le 10$.
The method of proof is similar to that used in Proposition
\xref{ThreeToSix}, except that the role of Corollary
\xref{ReducibleSeifert} is played by a theorem due to Agol [\cite{Agol}] and
Lackenby [\cite{Lackenby}] which implies that $\Delta(\alpha,\beta)\le 10$
in this situation.

This paper had its origins in unpublished work done by M. Culler,
C. MacA. Gordon and P. B. Shalen in 1984. This work established a
preliminary version of the Cyclic Surgery Theorem [\cite{CGLS}], giving an
upper bound of $5$ for the distance between two cyclic filling slopes
for a simple knot manifold, which was afterwards superseded by the
bound of $1$ established in [\cite{CGLS}]. The proof used character variety
techniques to reduce the result to a complicated topological
statement, which was in turn proved by combining the graph-theoretical
construction of Section \xref{ReducedHomotopies} and the use of the
subsurfaces $\Phi_k^\pm$ in the same way as is done in the proof of
Theorem \xref{HighGenusTheorem}. The techniques that were then
available could have produced a result qualitatively similar to
Theorem \xref{HighGenusTheorem}, but much weaker.

Independently and more recently, techniques similar to ours have been
used in a related context by Cooper and Long [\cite{CooperLong}] and Li
[\cite{Li}].  While our main results bound $\Delta(\alpha,\beta)$ under the
assumption that $\beta$ is a boundary slope and $\pi_1(M(\alpha))$
does not contain a non-abelian free group, the results in
[\cite{CooperLong}] and [\cite{Li}] imply a weaker bound for
$\Delta(\alpha,\beta)$ under the weaker assumption that $\beta$ is a
boundary slope and $\pi_1(M(\alpha))$ does not contain the fundamental
group of a closed surface of genus $> 1$.  Li proves a result which is
qualitatively similar to our Theorem \xref{DottedLengthBound}, but
with a bound of $6g+4m-6$ where Theorem \xref{DottedLengthBound}
provides the stronger bound of $4g+3m-4$.

We are indebted to Cameron Gordon for his role in the development 
of the ideas in this paper.

\section{Terminology and notation.}
We describe here various notational conventions that will be used
throughout the paper.

\paragraph
\tag{terminology}
If $X$ is a topological space, $|X|$ will denote the number of
components of $X$.

A (continuous) map $f\colon X\to Y$ of topological spaces will be
called {\it $\pi_1$-injective} if for each $x_0 \in X$, the
homomorphism $f_\sharp\colon\pi_1(X; x_0) \to \pi_1(Y; f(x_0))$ is
injective. A subset $A$ of a space $Y$ will be called {\it
$\pi_1$-injective} if the inclusion map $A\to Y$ is $\pi_1$-injective.

A {\it homotopy with domain $X$ and target $Y$} is a map
$H\colon X\times I \to Y$.  For each $t\in [0,1]$ we define $H_t\from
X \to Y$ by $H_t(x) = H(x,t)$. We shall sometimes refer to $H_t$ as
the {\it time-$t$ map} of $H$.

Let $H^1, \ldots, H^n$ be homotopies with domain $X$ and
target $Y$.  A homotopy $H$ with domain $X$ and target $Y$ will be
said to be a {\it composition} of $H^1, \ldots, H^n$ if there exist
numbers $0 = x_0 < x_1 \cdots < x_n = 1$ and monotone increasing linear
homeomorphisms $\alpha_i \colon [x_{i-1},x_i] \to [0,1]$ such that
$H(x,t) = H^i(x,\alpha_i(t))$ whenever $t\in [x_{i-1},x_i]$.

We shall say that a map $f \colon X \to Y$ between spaces is {\it
homotopic into a subset $B$ of $Y$} if $f$ is homotopic to a map $g
\colon X \to Y$ for which $g(X) \subset B$. A subset $A$ of a space
$Y$ is said to be {\it homotopic into a subset $B$ of $Y$} if the
inclusion map from $A$ to $Y$ is homotopic into $B$.

Let $f\from (X, Y) \to (Z, W)$ be a map of topological pairs, where
$Z$ is a connected $n$-manifold and $W\subset\bdry Z$ is an
$(n-1)$-manifold. In the case where $X$ is pathwise connected, we
shall say that $f$ is {\it essential} if it is $\pi_1$-injective as a
map from $X$ to $Z$ and is not homotopic, as a map of pairs, to a map
$f'\from (X, Y) \to (Z, W)$ where $f'(X) \subset W$. In general we
shall say that $f$ is {\it essential} if $X\ne\emptyset$ and $f$
restricts to an essential map from $(C,C\cap Y)$ to $(Z,W)$ for every
component $C$ of $X$. A map $f$ from a space $X$ to a manifold $Z$
will be termed {\it essential} if $f:(X,\emptyset)\to(Z,\partial Z)$
is essential.

We will denote the unit interval $[0,1]$ by $I$.  By an {\it essential
path} in a surface $F$ with non-empty boundary we shall mean
an essential map $f\from (I,\bdry I) \to (F,\bdry F)$.

A manifold $M$ is said to be {\it orientable} if all of its components
are orientable.  An {\it orientation} of $M$ is defined by a
choosing an orientation for each component; a manifold $M$ is said
to be {\it oriented} if an orientation of $M$ has been fixed.  Every
codimension $0$ submanifold of an oriented manifold $M$ inherits an
orientation.

Suppose that $A$ is a codimension $0$ submanifold of a manifold $M$
and that $h:A \to M$ is an embedding.  We say that $h$ {\it preserves
orientation} if it carries the orientation of $A$ inherited from $M$
to the orientation of $h(A)$ inherited from $M$; we say that $h$ {\it
reverses orientation} if it carries the orientation of $A$ inherited
from $M$ to the orientation of $h(A)$ which is the opposite of that
inherited from $M$.  In general, if $A$ is disconnected, there may
exist embeddings which neither preserve nor reverse orientation.

If $S$ is a compact surface then $\chi(S)$ denotes the Euler
characteristic of $S$ and $\genus(S)$ denotes the total genus of $S$,
i.e. the sum of the genera of the components of $S$.

\paragraph
A compact orientable $3$-manifold $M$ is said to be {\it irreducible}
if it is connected and every $2$-sphere in $M$ bounds a ball.  A
compact orientable $3$-manifold $M$ is said to be {\it boundary
irreducible} if it is connected and for any properly embedded disk $D$
in $M$ the curve $\partial D$ bounds a disk in $\partial M$.

Let $M$ be a compact, orientable, irreducible $3$-manifold, and let
$Q\subset\bdry M$ be a compact $\pi_1$-injective surface.  We define
an {\it essential surface} in $(M,Q)$ to be a compact surface $F$ in
$M$ such that (i) $\bdry F=F\cap\bdry M\subset Q$, and (ii) the
inclusion map $(F,\partial F)\to(Z,Q)$ is essential.  Note that
condition (ii) is equivalent to the condition that $F$ is
incompressible and not parallel to a subsurface of $Q$.  By an {\it
essential surface} in a $3$-manifold $M$ we mean an essential surface
in $(M,\bdry M)$.

Let $M$ be a compact, orientable, irreducible $3$-manifold and let $Q$
be a $\pi_1$-injective subsurface of $\partial M$. We will say that the
pair $(M,Q)$ is {\it acylindrical} if there does not exist any
essential map from $(S^1\times I, S^1 \times
\partial I)$ to $(M,\partial M)$ which sends $S^1 \times \partial I$
into $Q$.  By the Annulus Theorem [\cite{JacoShalen} IV.3.1], $(M,Q)$ is
acylindrical if and only if there is no essential annulus in $(M,Q)$.

We will say that a compact $3$-manifold $M$ is {\it atoroidal} if
there exists no essential map from $(S^1\times S^1,\emptyset)$ to
$(M,\partial M)$.

A closed orientable $3$-manifold will be said to be {\it very
small} if its fundamental group contains no non-abelian free subgroup.

\paragraph
\tag{Splittings}
Given a compact orientable $3$-manifold $M$ and a connected, properly
embedded, transversely oriented surface $F$ in $M$, we will implicitly
fix a regular neighborhood $N(F)$ of $F$ and a homeomorphism $f_F$
from $F\times [-1,1]$ to $N(F)$ which sends the standard orientation
of $[-1,1]$ to the transverse orientation of $F$ and maps
$F\times\{0\}$ to $F$.  The sets $f_F(F\times(-1,0])$ and
$f_F(F\times[0,1))$ will be denoted $N_-(F)$ and $N_+(F)$
respectively.

We will use the notation $M_F$ to denote the compact manifold $M -
f_F(F\times (-1/2,1/2))$.  The boundary of $M_F$ contains the two
subsurfaces $F\times\{-1/2\}$ and $F\times\{1/2\}$ which will be
denoted $F_-$ and $F_+$ respectively.  By restricting the projection
map $N(F) \to F$ we obtain two standard homeomorphisms $i_-\colon
F_- \to F$ and $i_+\colon F_+ \to F$.  Whenever it is convenient to do
so we will identify the manifold $M$ with the quotient of $M_F$
obtained by gluing $F_-$ to $F_+$ via the homeomorphism $i_+^{-1}\circ
i_-$.

The surface $F$ will be called a {\it semi-fiber} in $M$ if the
pair $(M_F,F_-\cup F_+)$ is an $I$-pair in the sense of [\cite{Jaco}];
that is, if there is an $I$-bundle $E$ over a surface and a
homeomorphism $h\from M_F\to E$ such that $h(F_-\cup F_+)$ is the
$\partial I$-bundle associated to $E$. (Note that $E$ may be
either a trivial $I$-bundle over a connected orientable surface,
or a twisted $I$-bundle over a  nonorientable surface of two
components.)

\paragraph
{\it Simple knot manifolds and slopes.}
We will say that a compact connected orientable $3$-manifold $M$ is
{\it simple} provided that (1) $M$ is irreducible and boundary
irreducible, (2) $M$ contains no essential surface of Euler characteristic
$0$, and (3) $M$ is not Seifert-fibered. If in addition the boundary of
$M$ is a torus we will say that $M$ is a {\it simple knot manifold}.

If $M$ is a simple knot manifold then an unoriented isotopy class of
homotopically non-trivial simple closed curves on $\partial M$ will be
called a {\it slope}. We will write $\Delta(\alpha,\beta)$ for the
geometric intersection number of two slopes $\alpha$ and $\beta$. We
will denote by $M(\alpha)$ the Dehn filling of $M$ determined by
$\alpha$.  If $F$ is a bounded essential surface in $M$ then the
boundary curves of $F$ all have the same slope, which will be called
the {\it boundary slope} of $F$.  A slope will be called {\it a
boundary slope} if it is the boundary slope of some essential surface.
A slope will be called a {\it strict boundary slope} if it is the
boundary slope of some essential surface which is not a semi-fiber in
$M$.

\section{Singular surfaces}
\tag{SingularSurfaces}
Let $M$ be a simple knot manifold. By a {\it singular surface} in $M$
we will mean a triple $(S,X,h)$, where $S$ is a compact, connected,
orientable surface, $X$ is a non-empty union of components of
$\partial S$, and $h:(S,X)\to(M,\partial M)$ is a map of pairs such
that (i) $h(S-X)\subset\int M$ and (ii) $h$ maps the components of $X$
homeomorphically onto disjoint, homotopically non-trivial simple
closed curves in $\partial M$.

If $(S,X,h)$ is a singular surface in
$M$, the components of $h(X)$ are all simple closed curves with the
same slope, which we shall call the {\it boundary slope} of $(S,X,h)$.

We shall say that two simple closed curves $\gamma$ and $\gamma'$ on a
torus $T$ are in {\it standard position} if there is a covering map
$p:\R^2\to T$ for which $p^{-1}(\gamma)$ and $p^{-1}(\gamma')$ are
Euclidean lines.

Let $F$ be a bounded essential surface in a simple knot manifold
$M$. A singular surface $(S,X,h)$ will be said to be {\it
well-positioned with respect to $F$} if (i) $h$ is transverse to $F$,
(ii) $h(\partial S-X)\cap F=\emptyset$, (iii) $h(X)$ is in standard
position with respect to $\partial F$, and (iv) each component of
$h^{-1}(F)$ is mapped by $h$ to an essential path or to a
homotopically non-trivial (possibly singular) closed curve in $F$. It
follows from (ii) that the arc components of $h^{-1}(F)$ have their
endpoints in $X$. Note that (iv) implies in particular that no simple
closed curve component of $h^{-1}(F)$ bounds a disk in $S$, and that
no arc component of $h^{-1}(F)$ is parallel in $S$ to an arc in $X$.

In this section we shall give several ways of constructing singular
surfaces that are well-positioned with respect to a given bounded
essential surface in a simple knot manifold.

\Proposition
\tag{Cameron}
Suppose that $F$ and $S$ are bounded essential surfaces in a simple knot
manifold $M$. Then the inclusion map from $S$ to $M$ is isotopic to an
embedding $h:S\to M$ such that the singular surface $(S,\partial S,h)$
is well-positioned with respect to $F$.
\EndProposition

\Proof
After an isotopy we may assume that $S$ and $F$ meet transversely and
that $\partial S$ is in standard position with respect to $F$. We
claim that every arc component of $S\cap F$ is essential in $F$. If
this is not the case, then there is a disk $D\subset F$ whose frontier
in $F$ is an arc $A$ such that $A=D\cap S$. Since $\partial M$ is a
torus, the essential surface $S$ is boundary-incompressible; hence $A$
is the frontier in $S$ of a disk $E\subset S$. Now $D\cup E$ is a
properly embedded disk in the simple knot manifold $M$ and is
therefore parallel in $M$ to a disk $J\subset\partial M$. As $\partial
J$ is made up of an arc in $\partial F$ and an arc in $\partial S$, we
have a contradiction to standard position, and the claim is proved.

To prove the proposition it now suffices to show that if some
component of $S\cap F$ is a homotopically trivial simple closed curve
$C$ in $F$, then $S$ is isotopic rel boundary to a surface $S'$ such
that $|S'\cap F|<|S\cap F|$. We may suppose $C$ to be chosen so that
there is a disk $D\subset F$ such that $\partial D=C=D\cap S$. Since
$S$ is essential, $C$ also bounds a disk $E\subset S$; since $M$ is
irreducible the disks $D$ and $E$ are isotopic by an isotopy that
fixes $C$. We may thus obtain the required surface $S'$ by moving
$(S-\int E)\cup D$ into general position with respect to $F$.
\EndProof

\Proposition
\tag{Irving}
Let $M$ be a simple knot manifold.  Let $F$ be a bounded essential
surface in $M$, and let $\alpha$ be a slope on $\bdry M$.  Suppose
that there exist a compact orientable surface $T$ and a
$\pi_1$-injective map $f\colon T \to M(\alpha)$ such that

\part{(1)}  $f(\partial T)\subset M-F\subset M\subset M(\alpha)$;

\part{(2)} there exists no map from $ T $ to $M$ which
agrees with $f$ on $\partial T$ and induces an injection from
$\pi_1(T)$ to $\pi_1(M(\alpha))$.

Then there exists a singular surface $(S,X,h)$ in $M$ which has
boundary slope $\alpha$ and is well-positioned with respect to $F$.
Moreover, we have $\genus S=\genus T$, and $|\partial S-X|=|\partial
T|$.

\EndProposition

\Proof
Let us write $M(\alpha)=M\cup V$, where $V$ is a solid torus and a
meridian curve of $V$ is identified with a curve of slope $\alpha$ in
$\partial M$.

By transversality and uniqueness of regular neighborhoods, there
exists a map $h\from T\to M(\alpha)$ such that

\part {(i)} $h$ is $\pi_1$-injective and agrees with $f$ on ${\bdry T}$;

\part {(ii)} The components of $h^{-1}(V)$ are disks in the interior
of $T$ which are mapped homeomorphically by $h$ to disjoint meridian disks of
$V$, whose boundaries are all in standard position with respect to
all of the components of $\partial F$;

\part {(iii)} each component of $h^{-1}(F)$ is a properly embedded
$1$-manifold in the bounded surface $h^{-1}(M)$.

Define the {\it complexity} of a map $h$ satisfying conditions
(i)---(iii) to be the ordered pair $(v(h),b(h))$ where $v(h)$
is the number of components of $h^{-1}(V)$ and $b(h)$ is the number of
components of $h^{-1}(F)$. Note that hypothesis (2) and condition (i)
imply that $v(h)$ must be strictly positive. 

Among all maps satisfying (i)---(iii), we suppose $h$ to be chosen so
that the complexity of $h^{-1}(F)$ is minimal with respect to
lexicographical order.

We set $S=T-h^{-1} (\int V)=h^{-1}(M)$ and $X=\partial S-\partial
T=h^{-1}(\partial M)$. By (ii), the components of $h(X)$ are simple
closed curves of slope $\alpha$.

We shall show that the map $h$ satisfies the following additional
condition:

\part {(iv)} Each component of $h^{-1}(F)$ is mapped by $h$
to an essential path or to a homotopically non-trivial (possibly
singular) closed curve in $F$. 

It follows from the definitions that if $h$ satisfies (i)--(iv) then
the triple $(S,X,h|S)$ is a singular surface, well-positioned with
respect to $F$ and having boundary slope $\alpha$. Furthermore, we
have $\genus S=\genus T$ and $|\partial S-X|=|\partial
T|$, since $S$ was obtained from $T$ by removing a collection of
disjoint disks bounded by $X$. Hence the proof of the proposition will
be complete when we have shown that $h$ satisfies (iv).

To prove (iv), first suppose that some curve component $C$ of
$h^{-1}(F)$ is mapped by $h$ to a homotopically trivial closed curve
in $F$. Then, since $h$ is $\pi_1$-injective, the curve $C$ bounds a
disk $D$ in $T$.  Thus there is a map $h'\colon T \to M(\alpha)$ which
agrees with $h$ on the complement of $D$, and maps $D$ into $F$.
Moving $h'$ into general position would produce a map of lower
complexity than $h$, giving a contradiction.

To complete the proof of (iv), we consider an arc component $A$ of
$h^{-1}(F)$ , and let $p,q\in X$ be the endpoints of $A$. We must show
that the path $h(A)$ is essential in $F$.

\def\u{{\bf u}}
\def\v{{\bf v}}
\def\w{{\bf w}}
It will be convenient to work in the smooth category for this
argument. Fix an orientation on the manifold $M$ and give $\partial M$
the induced orientation.  Given an ordered pair $(\gamma_1, \gamma_2)$
of oriented 1-manifolds on $\partial M$, intersecting transversely at
a point, we define the sign of their intersection at that point to be
the sign of the frame $(\u_1,\u_2)$ with respect to the orientation of
$\partial M$, where $\u_i$ is a tangent vector to $\gamma_i$ which is
positive with respect to its orientation.

Next we fix orientations of $F$ and $T$. The orientation of $T$
restricts to an orientation of $S$. Give the
components of $\partial F$ and $X$ the orientations induced from
those of $F$ and $S$, and then push the
orientation of $X$ forward under $h$ to obtain  orientations of the
meridian curves that make up $h(X)$.  The points $h(p)$
and $h(q)$ are transverse intersection points of the oriented
1-manifolds $\partial F$ and $h(X)$ on $\partial M$.  A
key step in the proof that $h(A)$ is an essential path is to show that
the signs of these intersections are opposite.

The orientations of $M$ and $F$ define a transverse orientation of
$F$. Because the map $h$ is transverse to $F$, we can pull back the
transverse orientation of $F$ to a transverse orientation of $A$ in
$S$.  Let $\u_p$ and $\u_q$ be tangent vectors to $X$
at $p$ and $q$ which are positive with respect to the orientation of
$X$. Let $\v_p$ and $\v_q$ be tangent vectors to
$\partial F$ at the points $h(p)$ and $h(q)$ which are positive with
respect to the orientation of $\partial F$.  Let $\w_p$ and $\w_q$ be
tangent vectors to $\partial M$ which are transverse to $\v_p$ and
$\v_q$ and positive with respect to the transverse orientation of $F$.
Observe first that, since $A$ is a properly embedded arc in
$S$, the transverse orientation of $A$ is consistent
with the orientation of $X$ at $p$ if and only if it is
inconsistent with the orientation of $X$ at $q$.
Second, observe that the signs of $(\v_p,\w_p)$ and $(\v_q,\w_q)$
agree.  It follows that the signs of $(Dh({\bf u}_p),\v_p)$ and
$(Dh(\u_q),\v_q)$ are opposite.  In other words, the signs of the
intersections of $\partial F$ and $h(X)$ are opposite
at $h(p)$ and $h(q)$.

Now assume that the path $h(A)$ is inessential in $F$.  In particular
the points $h(p)$ and $h(q)$ then lie on the same component of
$\partial F$.  But the signs of the intersections of $\partial F$ and
$h(X)$ are opposite at $h(p)$ and $h(q)$, and by (ii) the
components of $h(X)$ are all in standard position with
respect to all the components of $\partial F$. It follows that $p$ and
$q$ lie on distinct components $C_1$ and $C_2$ of $X$, and that
the orientations inherited from $h(X)$ by the meridian curves $h(C_1)$ and
$h(C_2)$ are opposite on the torus $\partial M$. Let $D_1,D_2\subset T$
be the disk components of $h^{-1}(V)$ bounded by $C_1$ and $C_2$, and define
$D$ to be a regular neighborhood of $D_1\cup A\cup D_2$ in $T$. Then
$h(\partial D)$ is homotopic in $M$ to a (possibly singular) closed
curve on the torus $\partial M$ which is the composition of conjugates
of the oppositely oriented meridian curves $h(C_1)$ and
$h(C_2)$.  In short, $h(\partial D)$ is homotopically trivial
in $M$.  Hence there is a map $h'\colon T \to M(\alpha)$ which agrees
with $h$ on the complement of $D$, and maps $D$ into $M$.  Then $h'$
has lower complexity than $h$.  This contradiction completes the proof
that $A$ is essential and hence that $h$ satisfies (iv).
\EndProof

The following results are corollaries to Propositions \xref{Cameron}
and \xref{Irving}. 

\Corollary
\tag{Disk}
Let $M$ be a simple knot manifold and $F\subset M$ an essential
bounded surface. Let $\alpha$ be a slope in $\partial M$. If
$M(\alpha)$ is very small, or more generally if $F\subset
M\subset M(\alpha)$ is not $\pi_1$-injective in $M(\alpha)$, then
there exists a singular surface $(S,X,h)$, well-positioned with
respect to $F$ and having boundary slope $\alpha$, such that $S$ is planar and
$\partial S-X$ is non-empty and connected.
\EndCorollary

\Proof
It follows from the definition of a simple manifold that if $F$ is a
bounded essential surface in a simple knot manifold $M$ then
$\chi(F)<0$, so that $\pi_1(F)$ is a non-abelian free group. Hence if
$M(\alpha)$ is very small then
$F\subset M(\alpha)$ is not $\pi_1$-injective in $M(\alpha)$. We must
show that if $F\subset M(\alpha)$ is not $\pi_1$-injective in
$M(\alpha)$ then there exists a singular surface $(S,X,h)$ with the
stated properties.
We will apply Proposition \xref{Irving}, taking $T$ to be a disk and
constructing the map $f\colon T \to
M(\alpha)$ as follows. Since $F$ is not $\pi_1$-injective in
$M(\alpha)$, there exists a homotopically non-trivial closed curve
$\gamma$ on the surface $F$ which is null-homotopic in $M(\alpha)$.
Deform the curve $\gamma$ by a homotopy to a curve $\gamma'$ which is
disjoint from $F$.  Let $f$ be a null-homotopy of $\gamma'$ in
$M(\alpha)$. Then $f$ is $\pi_1$-injective since $\pi_1(T)$ is
trivial, and condition (1) of Proposition \xref{Irving} holds by
construction. Since $F$ is essential, there does not exist a
null-homotopy of $\gamma'$ in $M$, so condition (2) of Proposition
\xref{Irving} as well. 
\EndProof

\Corollary
\tag{Sphere}
Let $M$ be a simple knot manifold and $F\subset M$ an essential
bounded surface. Let $\alpha$ be a slope in $\partial M$. If
$M(\alpha)$ is reducible then
there exists a singular surface $(S,X,h)$, well-positioned with
respect to $F$ and having boundary slope $\alpha$, such that $S$ is planar and
$X=\partial S$.
\EndCorollary

\Proof
If $M(\alpha)$ is reducible then there is an essential planar surface
in $M$ with boundary slope $\alpha$. The assertion therefore follows
from Proposition \xref{Cameron}.
\EndProof

\Corollary
\tag{Genusg}
Let $M$ be a simple knot manifold and $F\subset M$ an essential
bounded surface. Let $\alpha$ be a slope in $\partial M$ and let $g$
be a positive integer. Suppose that $\pi_1(M(\alpha))$ contains a
subgroup which is isomorphic to the fundamental group of a closed
connected surface of genus $g$, and that $\pi_1(M)$ does not.  Then
there exists a singular surface $(S,X,h)$, well-positioned with
respect to $F$ and having boundary slope $\alpha$, such that $\genus
S=g$ and $X=\partial S$.
\EndCorollary

\Proof If $\pi_1(M(\alpha))$ contains a genus-$g$ surface group,
then there is a $\pi_1$-injective map from a surface of genus $g$
to $M(\alpha)$. Note that condition (1) of Proposition
\xref{Irving} holds vacuously, and condition (2) holds because of
our hypothesis on $\pi_1(M(\alpha))$.  The assertion therefore follows
from Proposition \xref{Irving}. \EndProof

\Corollary
\tag{Seifert}
Let $M$ be a simple knot manifold and $F\subset M$ an essential
bounded surface. Let $\alpha$ be a slope in $\partial M$. If $M(\a)$
is a Seifert fibered space then there exists a singular surface
$(S,X,h)$, well-positioned with respect to $F$ and having boundary
slope $\alpha$, such that either (i) $\genus S=1$ and $X=\partial S$,
or (ii) $S$ is planar and $\partial S-X$ is non-empty and connected.
\EndCorollary

\Proof
If $\pi_1(M(\alpha))$ is finite, or if $M(\alpha)$ is reducible
and hence homeomorphic to $S^1\times S^2$ or $P^3\#P^3$,
then $M(\alpha)$ is very small and conclusion (ii) holds by
Corollary \xref{Disk}.

If $\pi_1(M(\alpha))$ is infinite and $M(\alpha)$ is irreducible,
then there exists a $\pi_1$-injective map $f: T\rightarrow
M(\alpha)$, where $T$ is a torus. Since $M$ is atoroidal, any
$\pi_1$-injective map from $T$ to $M$ is homotopic to a map into
$\partial M$.  But such a map cannot be $\pi_1$-injective as a map
from $T$ to $M(\alpha)\supset M$, since $\partial M$ bounds a
solid torus in $M(\alpha)$.  It follows $f$ is not homotopic to a
map whose image is contained in $M$.  Therefore $f$ satisfies the
hypotheses of Proposition \xref{Irving}, which asserts the
existence of a singular surface satisfying condition (i) of the
statement. \EndProof

\section{Reduced homotopies}
\tag{ReducedHomotopies}
\paragraph
\tag{HomotopyDef}
{\it Basic homotopies.\/}
Let $M$ be a simple knot manifold and let $F$ be a transversely
oriented essential surface in $M$.  A {\it homotopy} in $(M, F)$ with
domain $K$ is a homotopy $H$ with domain $K$ and target $M$ such that
$H(K\times \partial I)
\subset F$. A homotopy $H$ in $(M,F)$ is a {\it basic homotopy} if
$H^{-1}(F) = K \times \partial I$.

We shall say that a basic homotopy $H\from (K \times I, K \times
\partial I) \to (M, F)$ {\it starts on the $+$ side (or the $-$ side)}
if $H(K \times [0,\delta])$ is contained in $N_+(F)$ (or in $N_-(F)$)
for sufficiently small $\delta > 0$.  Similarly, $H$ will be said to
{\it end on the $+$ side (or the $-$ side)} if $H(K \times
[1-\delta,1])$ is contained in $N_+$ (or in $N_-$) for sufficiently
small $\delta > 0$. Of course in the case where $F$ separates $M$, a
basic homotopy starts on the $+$ side (or, respectively, the $-$ side)
if and only if it ends on the $+$ side (or the $-$ side).  In this
case we shall simply refer to it as a homotopy on the $+$ side
(or $-$ side).

As a notational convention we will treat the symbols $+$ and $-$ as
abbreviations for $1$ and $-1$ respectively.  Thus if
$\epsilon\in\{\pm1\}$ we may say that a homotopy $H$ starts on the
$\epsilon$-side, meaning that it starts on the $+$ side if
$\epsilon=1$ or that it starts on the $-$ side if $\epsilon = -1$.

There is a natural one-to-one correspondence between basic homotopies
$H$ in $(M,F)$ with domain $K$ which start on the $\epsilon$ side, and
homotopies $H'$ in $(M_F, F_-\cup F_+)$ with domain $K$ such that
$H'^{-1}(F_-\cup F_+) = K\times\partial I$ and $H'_0(K)\subset F_\epsilon$.
A basic homotopy $H$ is said to be {\it essential} if the corresponding
homotopy $H'\colon (K\times I, K\times \partial I) \to (M,F_-\cup F_+)$
is an essential map of pairs.

\paragraph
\tag{LengthDef}
{\it Reduced homotopies.\/}
Let $M$ be a simple knot manifold and $F$ an essential surface in $M$.
Let $K$ be a finite polyhedron. Suppose that $n$ is a positive
integer.  A homotopy $H\from (K \times I, K \times \partial I)\to (M,
F)$ is said to be a {\it reduced homotopy of length $n$ in $(M,F)$} if
for some $\epsilon\in\{\pm1\}$ we may write $H$ as a composition of
$n$ essential basic homotopies $H^1, \ldots, H^n$ in such a way that,
for $i=1,\ldots,n-1$, if $H^i$ ends on the $\epsilon$ side then
$H^{i+1}$ starts on the $-\epsilon$ side. In this case, if $H_1$
starts on the $\epsilon$ side then we shall also say that $H$ {\it
starts on the $\epsilon$ side.}

We define a {\it reduced homotopy of length $0$} in $(M,F)$ to be a
map $H$ from $K$ to $F$. The {\it time-$0$ and time-$1$ maps} of
$H$ are defined to be $H$ itself. If $H$ is a reduced homotopy of
length $0$ and $H'$ is a reduced homotopy of length $\ge0$ whose
time-$1$ (or time-$0$) map is equal to $H$, we define the {\it
composition} of $H$ with $H'$ (or of $H'$ with $H$) to be $H'$.

\paragraph
\tag{Glasses}
Let $\Gamma$ denote the 1-complex shown below, consisting of two vertices
$v_1$ and $v_2$, and three oriented edges $l_1$, $l_2$ and $b$.
\bigskip
\centerline{\epsfysize=1.25truein\epsfbox{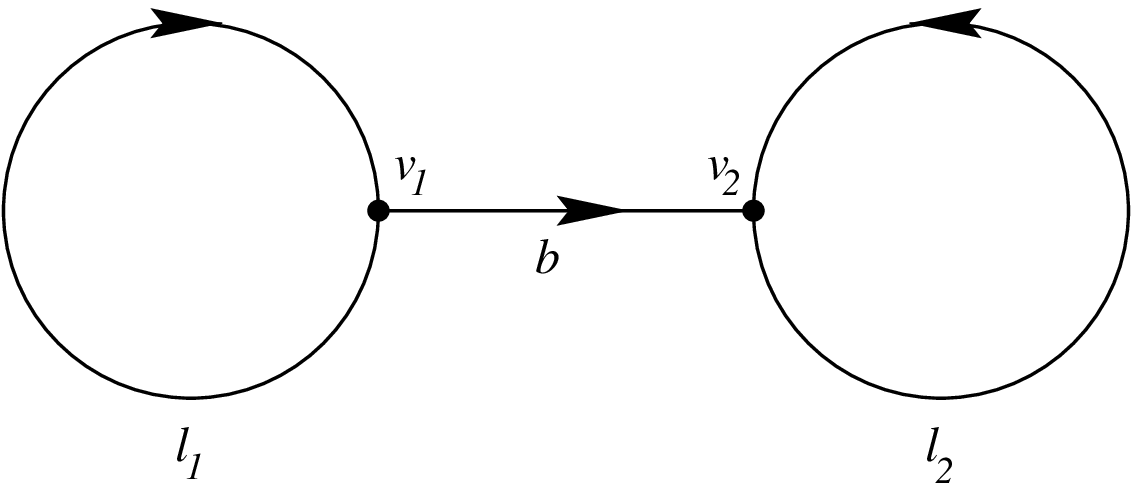}}
\centerline{ Figure \subsubnumber}
\tag{GlassesFig}
\medskip
By an {\it admissible pair of glasses} in a surface $F$ we will
mean a map $\gamma\colon \Gamma \to F$ such that the restriction
of $\gamma$ to $\bar b$ is an essential path, and, for $i = 1, 2$,
the restriction of $\gamma$ to $\bar l_i$ is a homeomorphism onto
a boundary component of $F$.  It is easy to show that if $\gamma$
is an admissible pair of glasses in a surface $F$ of negative
Euler characteristic then $\gamma\from(\Gamma,\bar l_1\cup \bar
l_2)\to(F,\partial F)$ is an essential map of pairs; and
furthermore that if $\alpha$ is any map from $S^1$ to $\Gamma$
which is not homotopic into $l_1$ or $l_2$,
then $g=\gamma\circ\alpha:S^1\to F$ is essential. In particular
this remark applies when $F$ is an essential surface in a simple
knot manifold, since the definition of simplicity then implies
that $\chi(F)<0$.

\paragraph
\tag{ArcToGlasses}
Let $F$ be an essential surface in a simple knot manifold $M$.
Suppose that $f \from (I,\bdry I) \to (F,\bdry F)$ is an essential
path and that
$$H\from (I\times I,I \times \partial I) \to (M,F)$$ is a reduced
homotopy of length $n$ such that $H_0 = f$ and $H_t(\partial I)
\subset \partial M$ for all $t\in I$.  Let us identify the set $\bar
b\subset\Gamma$ with $I$, respecting the orientations, so that $f$
becomes a map from $\bar b$ to $F$. Let us fix a product structure on
$\bdry M\homeo S^1\times S^1$ such that each simple closed curve in
$\bdry F$ has the form $\{x\}\times S^1$ for some $x \in S^1$, and for
$i=1,2$ let us fix orientation-respecting identifications of $\bar
l_i$ with $S^1$.  In this situation we shall describe a canonical way
to extend $f$ to an admissible pair of glasses $\hat f:\Gamma\to F$,
and to extend $H$ to a reduced homotopy $\hat H$ of length $n$ in the
pair $(M,F)$ such that $\hat H_0 = \hat f$.  The reduced homotopy
$\hat H$ will have the additional property that for all $t\in [0,1]$
the two closed curves $H_t|\bar l_i$, $i=1,2$, are homotopic to
boundary components of $F$.

For $i=1,2$ the point $f(v_i)$ lies in a component $\{x_i\}\times S^1$
of $\partial F$. As we have identified $l_i$ with $S^1$, there is a
unique rotation $\rho_i$ of $S^1$ such that $f(v_i) = (x_i,
\rho_i(v_i))$.  We extend $f$ to a map $\hat f \from \Gamma \to F$ by
sending $\theta \in \bar l_i$ to $(x_i, \rho_i(\theta))$ for $i=1,2$.
Since the map $f$ is an essential path, it follows that $\hat f$ is an
admissible pair of glasses.  Similarly, we extend $H$ to a homotopy
$\hat H$ as follows.  For $i=1,2$ we map each $(\theta,t)
\in \bar l_i\times I$ to $(x_{i,t},\rho_{i,t}(\theta))$ where
$x_{i,t}$ is chosen so that $H(v_i,t)\in\{x_{i,t}\}\times S^1$,
and $\rho_{i,t}$ is the unique rotation of $S^1$ such that
$H(v_i,t) = (x_{i,t},\rho_{i,t}(v_i))$. Since $H$ is a reduced
homotopy of length $n$, it is clear that $\hat H$ is also a
reduced homotopy of length $n$.

\paragraph
\tag{GlassesToCurve}
Suppose that $F$ is an essential surface in a simple knot manifold $M$
and that $\gamma \colon \Gamma\to F$ is an admissible pair of glasses.
Any length-$n$ reduced homotopy $\hat H$ in $(M,F)$ with domain
$\Gamma$ and time-$0$ map $\gamma$ can be used to produce a length-$n$
reduced homotopy in $(M,F)$ whose domain is $S^1$, and whose time-$0$
map is an essential map of $S^1$ into $F$.  We map $S^1$ to the
circuit in $\Gamma$ corresponding to the edge path $l_1 b l_2 b^{-1}$.
Composing this map with $\gamma$ we obtain a map $g\from S^1 \to F$.
The map $g$ is essential by the remark in \xref{Glasses}.  Composing
$\hat H$ with $g \times \hbox{\rm id}$ we obtain a reduced homotopy of
length $n$ whose time-$0$ map is $g$.

\Definition Let $G$ be a (possibly disconnected) graph in the
interior of a compact surface $S$, and let $E$ and $E'$ be distinct
edges of $G$. We shall say that $E$ and $E'$ are 
{\it adjacent parallel edges} if there exists a topological disk $D$,
whose boundary is a union of two arcs $A_1$ and $A_2$ with $A_1\cap
A_2=\partial A_1=\partial A_2$, and a map $i:D\to S$, such that
$i|(D-\partial A_1)$ is one-to-one, $i(A_j)=E_j$ for $j=1,2$, and
$i(\int D)\cap G=\emptyset$.  We shall say that $E$ and $E'$ are
{\it parallel edges} if there exist edges $E = E_1, \ldots, E_n = E'$
such that $E_i$ and $E_{i+1}$ are adjacent parallel edges for $i = i,
\ldots, n-1$.  Note that parallelism is an equivalence relation on the
set of edges of a graph.
\endproclaim

\paragraph
\tag{Nsnvf}
For the statement of the next lemma we introduce the following notation,
which will be used throughout the rest of the paper. If $s\ge0$,
$n\ge0$ and $v>0$ are integers, we define
$$N(s,n,v)=\max(1,6 + \left[{12s+6n-12\over v}\right]).$$

\Lemma
\tag{ValenceBound}
Let $S$ be a compact, connected orientable surface
having genus $s\ge0$ and $n\ge0$ boundary components. Let
$G\subset\int S$ be a non-empty (but possibly disconnected) graph with
$v>0$ vertices. Assume that no two (distinct) edges of $G$ are
parallel, and that no loop in $G$ bounds a disk whose interior is
disjoint from $G$. Then $G$ has a vertex of valence at most
$N(s,n,v)$.
\EndLemma

\Proof We may assume without loss of generality that $G$ has no vertex
of valence $0$ or $1$.

Let $\phi$ be any component of $S-G$. Since $G$ has no
valence-$0$ vertices, there exist a compact surface (with boundary)
$\hat \phi$ and a map $i_\phi:\hat \phi\to S$, such that $i_\phi$ maps
$\int\hat\phi$ homeomorphically onto $\phi$, and $\partial\hat\phi$
has a cell decomposition in which every edge or vertex is mapped
homeomorphically by $i_\phi$ onto an edge or vertex of $G$. We denote
by $o(\phi)$ the number of edges in the cell decomposition of
$\partial\hat\phi$.

We claim that 
$$o(\phi)\ge3\chi(\hat\phi)\leqno{(*)}$$ for every component
$\phi$ of $S-G$.  This is clear if $\bar\phi$ is not a disk,
as in that case $\chi(\hat\phi)\le0$. Now suppose that $\bar\phi$ is a
disk; we need to show that $o(\phi)\ge3$.  If $o(\phi)=1$ then there
is a loop in $G$ bounding a disk whose interior is disjoint from $G$;
this contradicts the hypothesis. If $o(\phi)=2$, and $i_\phi$ maps the
edges of $\hat\phi$ to distinct edges $E_1$ and $E_2$ of $G$, then
$E_1$ and $E_2$ are parallel, and again the hypothesis is
contradicted. If  $o(\phi)=2$, and
$i_\phi$ maps the edges of $\hat\phi$ to the same edge  of $G$, then
$S\cong S^2$ and the graph $G$ has one edge and two vertices. This
contradicts our assumption that $G$ has no vertices of valence
$1$. Thus ($*$) is proved in all cases.

Summing ($*$) over the components of $S- G$, we find that
$$2e=\sum_\phi o(\phi)\ge 3t,$$ where $e$ denotes the number of edges
of $G$ and $$t=\sum_\phi \chi(\hat\phi).$$ Now $$2-2s-n=\chi(S)=v-e+t\le
v-{e\over3}.$$ If $k$ denotes the minimum valence of any vertex of $G$
then we have $2e\ge kv$, and so $$2-2s-n\le v-{kv\over6}.$$
Since $k$ is an integer it follows that $k\le N(s,n,v)$.
\EndProof

\Proposition
\tag{GraphProposition}
Let $M$ be a simple knot manifold.  Let $F$ be a bounded connected
essential surface in $M$ with boundary slope $\beta$ and $m$ boundary
components.  Let $(S,X,h)$ be a singular surface in $M$ which has
boundary slope $\alpha\ne\beta$ and is well-positioned with respect to
$F$. Set $s=\genus S$, $n=|\partial S-X|$, $v=|X|$.

Then there exists an essential homotopy $H:I\times I\to M$ having
length at least $${m\D(\a,\b)\over N(s,n,v)}-1.$$ such that
$H_0$ is an essential path in $F$ and $H_t(\partial I)\subset
\partial M$ for all $t\in I$.
\EndProposition

\Proof
The definition of a singular surface gives $X\neq\emptyset$, so that
$v>0$.  Let $\hat S$ denote the surface obtained from $S$ by
identifying each component of $X$ to a point.  The surface $\hat S$
contains a non-empty graph whose vertices are the $v$ points in the
image of $X$, and whose edges are the images of the arc components of
$h^{-1}(F)$.  Each component of $h(X)$ is
a curve on $\bdry M$ of slope $\a$, and each component of $\bdry F$ is
a curve of slope $\b$.  Thus it follows from standard position that
each component of $X$ meets $h^{-1}(F)$ in $m\Delta(\alpha,
\beta)$ points.  In other words, each vertex of the graph $G$ has
valence $m\Delta(\alpha, \beta)$, which is strictly positive since
$\alpha\ne\beta$.  Since $(S,X,h)$ is well-positioned with respect to
$F$, there is no loop in $G$ bounding a disk whose interior is
disjoint from $G$.

We may apply Lemma \xref{ValenceBound} to a subgraph of $G$ containing
exactly one edge from each class of parallel edges.  We conclude that
the edges emanating from some vertex of $G$ fall into at most
$N(s,n,v)$ parallel classes.  Thus there must exist some class of
parallel edges containing at least $m\Delta(\alpha,\beta)/N(s,n,v)$
edges.  Label the edges of this class $E_1, \ldots, E_k$ where $k \ge
m\Delta(\alpha,\beta)/N(s,n,v)$.  For each $i=1, \ldots k$ let $A_i$
be the arc component of $h^{-1}(F)$ which maps to $E_i$ under the
quotient map from $S$ to $\hat S$.  For each $i = 1, \ldots, k-1$
there is a disk $Q_i$ in $S$ bounded by two subarcs of $X$ together
with $A_i$ and $A_{i+1}$.  Because $(S,X,h)$ is well-positioned with
respect to $F$, the interior of $Q_i$ is disjoint from $f^{-1}(F)$.
The restriction of $h$ to this quadrilateral disk defines a basic
homotopy $H^i$ whose time-$0$ map is the path $h(A_i)$ and whose
time-$1$ map is the path $h(A_{i+1})$.  The paths $h(A_i)$ are
essential because $(S,X,h)$ is well-positioned with respect to $F$.
It is an immediate consequence of standard position that each of the
basic homotopies $H^i$ is essential; since $h$ is transverse to $F$,
the composition of $H^1
\ldots H^{k-1}$ is a reduced homotopy of length $k-1$.
The time-$0$ map of this reduced homotopy is the essential path
$h(A_1)$.  Since $k \ge m\Delta(\alpha,\beta)/N(s,n,v)$ this completes
the proof of the proposition.
\EndProof

\Corollary
\tag{GeneralAcylindricalCorollary}
Let $M$ be a simple knot manifold that contains an essential surface
$F$ with boundary slope $\beta$.  Suppose that the pair $(M_F,F_-\cup
F_+)$ is acylindrical. Let $(S,X,h)$ be a singular surface with
boundary slope $\alpha$ which is well-positioned with respect
to $F$. Set $s=\genus S$, $n=|\partial S-X|$, $v=|X|$, and
$m=|\partial F|$. Then $$\Delta(\alpha,\beta)\le{N(s,n,v)\over m}.$$
\EndCorollary

\Proof
We may assume $\alpha\ne\beta$.
According to Proposition \xref{GraphProposition}, there exist an
essential path $\gamma:I\to F$ and a reduced homotopy in $(M,F)$ which
has time-$0$ map $\gamma$ and length at least $$l\ge{m\D(\a,\b)\over
N(s,n,v)}-1.$$
It therefore follows from \xref{ArcToGlasses} and
\xref{GlassesToCurve} that there exists a reduced homotopy of length
$l$ in $(M,F)$ whose time-$0$ map is a (possibly singular) essential
curve in $M$.  But since $(M_F,F_-\cup F_+)$ contains no properly
embedded essential annuli, the Annulus Theorem implies
that there can exist no essential basic homotopy in $M$ whose time-$0$
map is an essential curve in $F$. Hence we must have
$${m\D(\a,\b)\over N(s,n,v)}-1\le0,$$ which is equivalent to the
conclusion of the corollary.
\EndProof

\Corollary
\tag{AcylindricalCorollary}
Let $M$ be a simple knot manifold that contains an essential surface
$F$ with boundary slope $\beta$.  Suppose that the pair $(M_F,F_-\cup
F_+)$ is acylindrical. Let $\alpha$ be a slope on $\partial M$.

\part{(1)} If $M(\alpha)$ is a very small manifold, then
$\Delta(\alpha,\beta)\le 5$.

\part{(2)} If $M(\alpha)$ is a Seifert fibered space, then
$\D(\alpha,\beta)\le 6$.

\EndCorollary

\Proof
To prove (1) we invoke Corollary \xref{Disk} to obtain a singular
surface $(S,X,h)$, well-positioned with respect to $F$, such that
$\genus S=0$ and $|\partial S - X|=1$. The conclusion now follows from Corollary
\xref{GeneralAcylindricalCorollary} because for any $v\ge1$ we have
$N(0,1,v)\le 5$.  To prove (2) we invoke Corollary \xref{Seifert} to
obtain a singular surface $(S,X,h)$, well-positioned with respect to
$F$, such that either $\genus S=0$ and $|\partial S - X|=1$,
or $\genus S=1$ and $|\partial S - X|=0$. The conclusion now
follows from Corollary \xref{GeneralAcylindricalCorollary} because for
any $v\ge1$ we have $N(0,1,v)\le 5$ and $N(1,0,v)=6$.
\EndProof

\section{Essential intersections}
\tag{EssentialIntersections}
This section introduces a version of the notion of essential
intersection for subsurfaces of a $2$-manifold, a notion which has
appeared implicitly in much of the literature on the characteristic
submanifold of a Haken manifold. The version we present here is
adapted to the case of ``large'' subsurfaces, which we now define.

Let $S$ be a compact orientable surface. We say that a subsurface $A$
of $S$ is {\it large} if each component of $A$ is $\pi_1$-injective
and has negative euler characteristic. Note that the empty
set is considered to be a large subsurface according to this
definition. If the components of $S$ have Euler characteristic $\ge0$
then the empty set is the only large subsurface.

The {\it large part} of a $\pi_1$-injective subsurface $A$ of $S$,
denoted by ${\cal L}(A)$, is the union of all the large components of
$A$.

The next lemma is preliminary to the proof of Proposition
\xref{EssInt}, which will provide the definition of the ``large
intersection'' of two large subsurfaces.

\Lemma
\tag{SurfaceTwo}
Let $A$ and $B$ be large subsurfaces of a compact orientable surface
$S$, and suppose that $A$ is homotopic into $B$.

\part {$(1)$} $A$ is  isotopic in $S$ to a subsurface of $B$.

\part {$(2)$} If $B$ is
homeomorphic to a large subsurface of $A$
then $A$ and $B$ are isotopic subsurfaces of $S$.

\part {$(3)$} If  $B$ is homotopic
into $A$ then $A$ and $B$ are isotopic subsurfaces of $S$.

\EndLemma

\Proof
Without loss of generality we may assume that $S$ is connected.

We begin with the proof of (1).  We may assume that $A$ has been
chosen within its isotopy class so that $\partial A$ meets $\partial
B$ transversely and in the minimal number of points.  Moreover, we may
assume that $A$ has been chosen, among all surfaces in its isotopy
class which minimize $|\partial A \cap \partial B|$, to  minimize
the number components of $\partial A$ which are not contained in
$B$. Under these assumptions we will show that $A\subset B$, proving (1).

Let $A_0$ be any component of $A$. By hypothesis, $A_0$ is homotopic
into a component $B_0$ of $B$. We will prove that $A_0\subset B_0$. As
$A_0$ is an arbitrary component of $A$, the assertion will follow.

Fix a base point in $B_0$, and consider the covering $p \colon \tilde
S \to S$ determined by $\im(\pi_1(B_0)\to\pi_1(S))$.  There is a
subsurface $\tilde B_0$ of $\tilde S$ which is mapped homeomorphically
onto $B_0$ by $p$. Since $B_0$ is large, $\partial B_0$ is
$\pi_1$-injective in $S$, and hence each component of $X={\tilde
S-\int\tilde B_0}$ is a half-open annulus. Since $A_0$ is homotopic
into $B_0$, the inclusion $i:A_0\to S$ is homotopic to a map which
admits a lift to $\tilde S$. Hence $i$ itself admits a lift to $\tilde
S$; that is, we have a subsurface $\tilde A_0$ of $\tilde S$ which is
mapped homeomorphically onto $A_0$ by $p$. In order to show that
$A_0\subset B_0$ it suffices to show that $\tilde A_0\subset\tilde
B_0$.

As a preliminary, we will prove that $\partial\tilde A_0$ is disjoint
from $\partial\tilde B_0$.  Assume this is false. Then $\partial\tilde
A_0$ contains a properly embedded arc $\tilde\alpha\subset X$. Since
the component $X_0$ of $X$ containing $\tilde\alpha$ is a half-open
annulus, $\partial\tilde\alpha$ is the boundary of an arc
$\tilde\beta\subset \partial X_0\subset\tilde B_0$, and
$\tilde\alpha\cup\tilde\beta$ bounds a disk $\tilde D\subset\tilde
S$. Among all disks in $\tilde S$ whose boundaries are made up of an
arc in $p^{-1}(\partial A)$ and an arc in $p^{-1}(\partial B)$, choose
one, say $\tilde D_0$, which is minimal with respect to
inclusion. Write $\partial \tilde
D_0=\tilde\alpha_0\cup\tilde\beta_0$, where $\tilde\alpha_0\subset
p^{-1}(\partial A)$ and $\tilde\beta_0\subset p^{-1}(\partial B)$ are
arcs with $\partial\tilde\alpha_0=\partial\tilde\beta_0$. We claim
that $p|\partial\tilde D_0$ is one-to-one. It is clear from the minimality
of $\tilde D_0$ that $p(\int\tilde\alpha_0)\cap
p(\int\tilde\beta_0)=\emptyset$. If $p|\int\tilde\alpha_0$
is not one-to-one, then $\alpha=p(\int\tilde\alpha_0)$ is an entire
component of $\partial A$, which contains $p(\partial\tilde\alpha_0)$
and therefore meets $B$; hence $\int\tilde\alpha_0$ meets $p^{-1}(B)$,
and the minimality of $\tilde D_0$ is contradicted. Thus
$p|\int\tilde\alpha_0$ is one-to-one, and similarly
$p|\int\tilde\beta_0$ is one-to-one.   Hence if $p|\partial\tilde D_0$ is
not one-to-one then $\alpha_0=p(\tilde\alpha_0)$ and
$\beta_0=p(\tilde\beta_0)$ are simple closed curves meeting in a
single point; as these curves must be components, the intersection is
transverse. But this is impossible, as the hypothesis that $A$ is
homotopic into $B$ implies that $\alpha_0$ and $\beta_0$ are homotopic
to disjoint curves. Thus $p|\partial D_0$ must be one-to-one.

According to [\cite{Epstein}, Lemma 1.6], it follows that $p|\tilde D_0$ is
one-to-one; hence $D_0 = p(\tilde D_0)$ is a disk in $S$ whose
boundary consists of two arcs, one in $B$ and one in $A$, and $\int
D_0$ is disjoint from $\partial A$ and from $\partial B$. It follows
that $\partial B$ is isotopic to a curve system that meets $\partial
A$ transversely in fewer points than $\partial B$, a contradiction to
our choice of $B$. This proves that $\partial\tilde
A_0\cap\partial\tilde B_0=\emptyset$.

We are now ready to prove that $\tilde A_0\subset \tilde B_0$. If this
is false, then either $\tilde A_0\subset X$ or $\tilde A_0\cap\partial
X\not=\emptyset$. If $\tilde A_0\subset X$, then since the components of
$X$ are half-open annuli, $\im(\pi_1(A_0\to\pi_1(S))$ is cyclic, a
contradiction since $A$ is a large subsurface of $S$. There remains
the possibility that $\tilde A_0$ meets the boundary of some component
$X_1$ of $X$. Since $\partial
\tilde A_0\cap\partial \tilde B_0=\emptyset$, 
the compact subsurface $\tilde Z=\tilde A_0\cap X_1$ must have
$\tilde\gamma_1=\partial X_1$ as one boundary component. If
$\tilde\gamma$ is any other component of $\partial \tilde Z$, then
$\tilde\gamma$ lies in the interior of the half-open annulus $X_1$;
furthermore, $\tilde\gamma$ is a component of $\partial \tilde A_0$,
and is therefore homotopically non-trivial (since the large subsurface
$A_0$ of $S$ is in particular $\pi_1$-injective). It follows that
$\tilde Z$ is an annulus bounded by $\tilde\gamma_1$ and a single
component $\tilde\gamma_2$ of $\partial\tilde A_0$.  Since $p$ maps
$\tilde A_0$ homeomorphically onto $A_0$, it in particular maps
$\tilde Z$ homemorphically onto an annulus $Z\subset A_0$. There is an
isotopy supported on a small neighborhood $W$ of $Z$ in $S$ which carries
$\gamma_1=p(\gamma_1)$ into $\int B_0$. Since $W$ may be taken to meet
$\partial A$ only in $\gamma_1$, it is clear that the image of $A$
under the isotopy has the same intersection with $B$ as $A$ has, and
has fewer components $\partial A$ which are not contained in
$B$ than $A$ has. This contradiction to the minimality of $A$
completes the proof of (1).

We now turn to the proof of (2).  By part (1) we may assume that $A$
is a subsurface of $\int B$.  Let $B'$ be a $\pi_1$-injective
subsurface of $\int A$ which is homeomorphic to $B$.  Any disk
component of $\overline{B - A}$ would also be a component of
$\overline{S - A}$.  Since $A$ is large it follows that no
component of $\overline{B - A}$ is a disk.  Hence $\chi(B)
\le \chi(A)$.  Similarly, $\chi(A) \le \chi(B') = \chi(B)$.  Thus each
component of $\overline{B - A}$ is an annulus, and each
component of $\overline{A - B'}$ is an annulus.  Since $A$ and
$B$ are large it follows that each component of $B$ contains at least
one component of $A$ and that each component of $A$ contains at least
one component of $B'$.  Since $B$ is homeomorphic to $B'$, a given
component of $B$ can contain only one component of $B'$. Hence each
component of $B$ contains exactly one component of $A$.

We may therefore index the components of $A$ as $A_1, \ldots, A_n$,
and index the components of $B$ as $B_1, \ldots, B_n$, in such a way
that $A_i$ is a subsurface of $B_i$ for $i=1,\ldots, n$, and each
component of $\overline{B_i- A_i}$ is an annulus.

To complete the proof of (2) it suffices to show that there is no
annulus component of $\overline{B_i- A_i}$ whose boundary is
contained in $A_i$.  If such an annulus component did exist the genus
of $B_i$ would be strictly greater than the genus of $A_i$.  In
particular we would have $\genus B=\sum\genus B_i>\sum\genus
A_i=\genus A$. (Recall that $\genus B$ denotes the total genus of $B$.)
But since $B$ is homeomorphic to the subsurface $B'$ of $A$, we have
$\genus B=\genus B'\le\genus A$.  This contradiction completes the proof of
(2).

To prove (3), note that if $B$ is homotopic into $A$
then by (1) it is isotopic into $A$.  The time-$1$ map of the
isotopy is a homeomorphism onto a subsurface of $A$ which is clearly
large since it is isotopic to $B$. It now follows from (2) that $A$
and $B$ are isotopic.
\EndProof

A map $f$ from a finite polyhedron $K$ to a surface $S$ is called {\it
large} if for each component $K_0$ of $K$, the subgroup
$f_\#(\pi_1(K_0))$ of $\pi_1(S)$ is non-abelian.

\Proposition
\tag{EssInt}
Suppose that $A$ and $B$ are large
subsurfaces of a compact orientable surface $S$. Then up to non-ambient isotopy
there is a unique large subsurface $C$ of $S$
with the following property.

\part{$(\ast)$} Any large map from a polyhedron into
$S$ is homotopic into $C$ if and only if it is homotopic into $A$ and
homotopic into $B$.

Furthermore, if $C\subset\int S$ is a large
subsurface satisfying $(\ast)$, then there are subsurfaces
$A_0\subset\int S$ and $B_0\subset\int S$, isotopic to $A$ and $B$,
such that $\partial A_0$ and $\partial B_0$ meet transversely and
$C={\cal L}(A_0\cap B_0)$.

\EndProposition

\Proof 
We may assume without loss of generality that $A,B\subset\int S$.
We prove uniqueness first. Suppose that $C$ and $C'$ are large
 subsurfaces of $S$ satisfying ($\ast$).  Then the
inclusion $C \to S$ is homotopic into $C'$, and vice versa. Lemma
\xref{SurfaceTwo}(3) now implies that $C$ and $C'$ are isotopic in $S$.

To complete the proof of the Proposition, it suffices to show that
there is a subsurface $B_0$ of $\int S$ isotopic to $B$, such that
$\partial A$ and $\partial B_0$ meet transversely and such that ${\cal
L}(A\cap B_0)$ satisfies $(\ast)$. We define $B_0$ as follows: among
all subsurfaces of $\int S$ which are isotopic to $B$, and whose
boundaries meet $\partial A$ transversely, we choose $B_0$ so that the
number of points of $\partial A\cap\partial B_0$ is as small as
possible. We set $C={\cal L}(A\cap B_0)$.

It is clear that if a large map from a connected polyhedron $K$ into
$S$ is homotopic into $C$ then it is homotopic into $A$ and homotopic
into $B$. The proof of $(\ast)$ will be completed by showing that,
conversely, if a large map $f:K\to S$ is homotopic into $A$ and
homotopic into $B$ then it is homotopic into $C$. We may assume
without loss of generality that $K$ is non-empty and connected. We fix
a component $A_1$ of $A$ such that $f$ is homotopic into $A_1$.

Let $\tilde S$ denote the covering space of $\int S$ corresponding to
the subgroup $\im(\pi_1(A_1)\to\pi_1(\int S))$. Then $\tilde S$ has a
subsurface $\tilde A$ such that the covering projection $p:\tilde
S\to\int S$ maps $\tilde A$ homeomorphically onto $A_1$.  Set
$X=\tilde S-\int\tilde A$. Since $A$ is $\pi_1$-injective, each
component of $X$ is a half-open annulus meeting $\tilde A$ in a
component of $\partial\tilde A$.

Since $f:K\to S$ is homotopic into $A_1$ it admits a lift $\tilde f$
to the covering space $\tilde S$.  Since $f$ is also homotopic into
$B_0$, there is a homotopy from  $\tilde f$ to a map $g:K\to \tilde S$
such that $g(K)$ is contained in a component
$\tilde B_0$ of $p^{-1}(B_0)$.  

We claim that $\tilde B_0\cap\tilde A\ne\emptyset$. Suppose not; then
$g(K)\subset\tilde B_0\subset X$. But
$g_\sharp(\pi_1(K))\subset\pi_1(\tilde S)$ is non-abelian since $f$ is
a large map. Since each component of $X$ is a half-open annulus, we
have a contradiction, and the claim is proved.

We next claim that if $Z$ is any component of $\tilde B_0\cap X$, then
$Z$ deforms into $W=Z\cap\partial X\subset\partial\tilde A$. (This
means that the identity map of $Z$ is homotopic in $Z$, rel $W$, to a
map whose image is contained in $W$.) To prove the claim, first note
that $W$ is a $1$-manifold in the boundary of the $2$-manifold
$Z$. Furthermore, since $\tilde B_0\cap\tilde A\ne\emptyset$, we have
$W\ne\emptyset$. Hence, in order to prove that $Z$ deforms into $W$,
it suffices to show that if $\alpha\subset Z$ is an arc with endpoints
in $W$, then $\alpha$ is parallel in $Z$ to an arc in $W$.  Since $X$
is a half-open annulus, $\alpha$ is parallel in $X$ to an arc
$\beta\subset\partial X$. We must show that the disk $D\subset X$
bounded by $\alpha\cup\beta$ is contained in $Z$.

If $D\not\subset Z$ then $D$ contains a component $\gamma$ of the
frontier of $Z$ relative to $X$. Note that
$\gamma\subset\partial\tilde B_0$, and that $\gamma$ is a properly
embedded $1$-manifold in $X$. But $\gamma$ cannot be a simple closed
curve, since $B_0$ is $\pi_1$-injective in $S$; hence $\gamma$ is an
arc whose endpoints lie in $\beta$. Thus there is a disk in $\tilde S$
whose boundary is the union of the arc $\gamma\subset p^{-1}(\partial
B_0)$ and a sub-arc of $\beta\subset p^{_1}(\partial A)$. Among all
disks in $\tilde S$ bounded by the union of an arc in $p^{-1}(\partial
A)$ with an arc in $p^{-1}(\partial B_0)$, choose one, say $D'$, which
is minimal with respect to inclusion. According to [\cite{Epstein},
Lemma 1.6], $p|D'$ is one-to-one and hence $ p( D')$ is a disk in $S$
whose boundary consists of two arcs, one in $B_0$ and one in $A$,
which meet at their endpoints.  Since $B_0$ is large we have $\int p(
D')\cap B_0=\emptyset$. We may therefore
isotope
$B_0$ across the disk
$p(D')$ to obtain a surface whose boundary meets $\partial A$
in fewer points than $\partial B_0$. This contradicts our choice of
$B_0$, and the claim is proved.

We have shown that  $\tilde B_0\cap\tilde A\ne\emptyset$ and that
every component of $\tilde B_0\cap X$ deforms into its intersection
with $\partial X=\partial\tilde A$. It follows that $\tilde B_0$
deforms into  $\tilde B_0\cap\tilde A$. Since $\tilde f$ is
homotopic to $g$ and $g(K)\subset \tilde B_0$, it follows that $\tilde f$ is
homotopic in $\tilde S$ to a map of $K$ into  $\tilde B_0\cap\tilde
A$.  In particular, $f$ is homotopic in $S$ to a map $f'$ whose image is
contained in $B_0\cap A_1\subset B_0\cap A$. As $f$ is large it
follows that $f(K)\subset{\cal L}(B_0\cap A)=C$.
\EndProof

\Definition
\tag{DefinitionOfEssentialIntersection}
If $A$ and $B$ are large subsurfaces of a compact orientable surface $S$ then
the subsurface provided by Proposition \xref{EssInt}, which is
well-defined up to non-ambient isotopy, will be called the {\it
 large intersection} of $A$ and $B$ and will be denoted $A
\wedge_{\cal L} B$.  
Clearly $A
\wedge_{\cal L} B$ is isotopic to $B \wedge_{\cal L} A $, and $A
\wedge_{\cal L} (B \wedge_{\cal L} C)$ is isotopic to $(A \wedge_{\cal
  L} B) \wedge_{\cal L} C$.
\EndDefinition

The following result will be needed in the next two sections.

\Proposition
\tag{Containment}
Suppose that $A$ is a large subsurface of a compact orientable surface
$S$, and that $f$ and $g$ are large maps of a polyhedron $K$ into
$A$. If $f$ and $g$ are homotopic in $S$ then they are homotopic in
$A$.
\EndProposition

\Proof We may assume that $K$ and $S$ are connected. Let $A_0$
denote the component of $A$ containing $f(K)$, and let
$p\colon\tilde S\to S$ denote the largest covering space of $S$ to
which the inclusion map $i\colon A_0\to S$ lifts. If $\tilde
i:A_0\to\tilde S$ is the lift of $i$ then $\tilde A_0=\tilde
i(A_0)$ is a deformation retract of $\tilde S$, and each component
of $\tilde S-\tilde A_0$ has cyclic (and possibly trivial)
fundamental group. Let $\tilde f:K\to\tilde S$ be a lift of $f$
with $\tilde f(K)\subset\tilde A_0$, and let $\tilde g:K\to\tilde
S$ be a lift of $g$ which is homotopic to $\tilde f$ in $\tilde
S$. Since $g$ is large, $\tilde g(\pi_1(K))$ is not cyclic; thus
$\tilde g(K)$ cannot be contained in a component $\tilde
A_0'\ne\tilde A_0$ of $p^{-1}(A_0)$. Hence $\tilde
g(K)\subset\tilde A_0$, and $\tilde f$ and $\tilde g$ are
homotopic in $\tilde A_0$. The conclusion follows. \EndProof

\Lemma
\tag{Invariant}
Let $A$ be a large subsurface of a compact orientable surface
$B$, and suppose that $B$ admits an involution $\tau$ such that
$\tau(A)$ is isotopic to $A$ in $B$.  Then $A$ is isotopic to a
subsurface $A'$ of $B$ such that $\tau(A') = A'$.
\EndLemma

\Proof
Choose a negatively curved metric on $B$ which is invariant under the
involution $\tau$ and has the property that the boundary
components of $B$ are geodesics. Given a large subsurface $C$ of $B$
we divide its boundary components into two types:

\part {Type 1:} those which are not isotopic to any other boundary
component of $C$.

\part {Type 2:} those which are isotopic to some other boundary
component of $C$.

There is a small positive real number $t_0$ such that $A$ can be isotoped
in $B$ to a surface $A'$ whose type 1 boundary components
are geodesics and whose type 2 boundary components form the boundary
of a $t_0$-bicollar of a finite union of geodesics. Then $\partial A'$ is
invariant under the involution $\tau$. To see that $\tau(A') = A'$ let
$A_0'$ be
a component of $A'$. If $\tau(A_0') \not \subset A'$ then the interior of
$\tau(A_0')$ is disjoint from $A'$ and hence we may isotope $\tau(A_0')$
into the complement of $A'$.
Let $p: \tilde B \to B$ be the maximal connected cover for which the
inclusion map $i: \tau(A_0') \to B$ lifts to a map $\tilde i: \tau(A_0')
\to \tilde B$ and define $Y = \tilde i(\tau(A_0'))$. Then $X = \tilde B
- \int Y$ is a disjoint union of half-open annuli. By hypothesis
$\tilde i$ is homotopic to a map with image in $X$ so that in
particular $Y$ may be homotoped into an annulus.  It follows from Lemma
\xref{SurfaceTwo} that $Y$ is isotopic to a subsurface of an annulus.
Since $Y$ is $\pi_1$-injective and homeomorphic to $A_0$, this
contradicts our assumption that $A_0$ is large.  Thus $\tau(A_0')
\subset A'$ and we deduce that $\tau(A') = A'$.
\EndProof

\paragraph
\tag{OuterParts}
Let $A$ be a subsurface of a compact orientable surface $S$. By
an {\it outer component } of $A$ we will mean a component of $A$ which
contains a simple closed curve that is isotopic to a boundary
component of $S$.  We will say that $A$ is {\it outer} if every
component of $A$ is outer.  We define the {\it outer part } of $A$,
denoted $\dot A$, to be the union of all outer components of $A$,
i.e. the largest outer subsurface of $A$.  Thus $A$ is an outer
subsurface if and only if $A = \dot A$.

If $A$ and $B$ are subsurfaces of $S$ with $A\subset B$ note that
$\dot A \subset \dot B$.  If $A$ and $B$ are large subsurfaces of
$S$, we define $A\dotlint B$ to be the outer part of $A\lint B$.
It is a formal consequence of this definition that a large outer
subsurface of $S$ is homotopic into both $A$ and $B$ if and only
if it is homotopic into $A\dotlint B$.  It follows from this,
together with Lemma \xref{SurfaceTwo}, that if $A$, $B$ and $C$
are large subsurfaces of $S$, then  $A\dotlint B$ is homotopic to
$B\dotlint A$ and $(A\dotlint B)\dotlint C$  is isotopic to
$A\dotlint (B\dotlint C)$.

\Lemma
\tag{DottyLemma}
Suppose that $A$ and $B$ are large subsurfaces of
a compact orientable surface $S$.  Then $A \dotlint B$ is isotopic in $S$ to
$\dot A \dotlint \dot B$.
\EndLemma

\Proof
Since each component of $\dot A \dotlint \dot B$ contains a simple
closed curve that is isotopic to a boundary component of $S$, it is
clear that $A \dotlint B \supset \dot A \dotlint \dot B$.  To prove
the reverse inclusion, let $C$ be a component of $A \dotlint B$.  We
have that $C$ is isotopic into $X \lint Y$ where $X$ and $Y$ are
components of $A$ and $B$ respectively.  Since $C$ contains a simple
closed curve isotopic to a boundary component of $S$, so do $X$ and
$Y$.  Thus $X$ and $Y$ are in fact components of $\dot A$ and $\dot B$
respectively.  This shows that $C$ is isotopic into a component of
$\dot A \lint \dot B$.  But, since $C$ contains a simple closed curve
isotopic to a
boundary component of $S$, it is isotopic into a component of $\dot A
\dotlint \dot B$.  Since $C$ was an arbitrary component of
$A \dotlint B$, we have shown that $A \dotlint B$ is homotopic into
$\dot A \dotlint \dot B$.  The lemma now follows from
Lemma \xref{SurfaceTwo}.
\EndProof

\deepsection{Reduced homotopies and the characteristic pair}
\tag{GeneralBounds}
The final result of this section, Theorem \xref{LengthBound}, concerns a pair
$(M,F)$ where $M$ is a simple knot manifold and $F\subset M$ is an
essential surface in $M$ which is not semi-fiber.  The theorem
provides an upper bound, in terms of the genus and the number of
boundary components of $F$, for the length of a reduced homotopy in
$(M,F)$ having a large time-$0$ map.

\subsection{Splitting surfaces.}
\tag{SplittingSurfaces}
\Definition
Let $M$ be simple knot manifold. A {\it splitting surface in $M$} is a
transversely oriented essential surface $\tilde F\subset M$ such that
$M_{\tilde F}$ is a disjoint union of two compact submanifolds
$M_{\tilde F}^+$ and $M_{\tilde F}^-$ with the property that
$N_\epsilon(\tilde F)\subset M_{\tilde F}^\epsilon$ for
$\epsilon=\pm1$.
\EndDefinition

It is easy to see that $M_{\tilde F}^+$ and $M_{\tilde F}^-$ are
uniquely determined by $\tilde F$.  Note that any transversely
oriented, separating, connected, essential surface in $M$ is a
splitting surface.  In general, however, a splitting surface need not
be connected.

Since $\tilde F$ comes equipped with a transverse orientation, any
orientation of $M$ induces an orientation of $\tilde F$.  An
orientation of $\tilde F$ will be called {\it consistent} if it is
induced from an orientation of $M$ in this way.  If $\tilde F$ has $n$
components then it has $2^n$ possible orientations; but since $M$
is connected, only two of these orientations are consistent.
  
If $\tilde F$ is a splitting surface then for $\epsilon = \pm 1$ the
natural map from $M_{\tilde F}^\epsilon$ to $M$ is injective and we shall
identify $M_{\tilde F}^+$ and $M_{\tilde F}^-$ with submanifolds of
$M$.  In particular $\tilde F_\epsilon$ is identified with $\tilde F$
via the homeomorphism $i_\epsilon$ (see \xref{Splittings}).
Furthermore we shall identify basic homotopies in $(M,F)$ with the
corresponding homotopies in $(M_{\tilde F}, \tilde F_-\cup\tilde F_+)$
(see \xref{HomotopyDef}.)

An orientation of $M$ determines an orientation of $M_{\tilde
F}^\epsilon$ by restriction, and this orientation of $M_{\tilde
F}^\epsilon$ induces an orientation of $\tilde F_\epsilon$. The two
orientations of $\tilde F_\epsilon$ which arise in this way are
identified with the two consistent orientations of $\tilde F$.

\subsection{Supports of reduced homotopies}
\tag{EssHom}
Throughout this subsection we shall assume that $M$ is a simple knot
manifold, and that $\tilde F$ is a splitting surface in $M$.

\paragraph
\tag{JSJStuff}
It follows from the characteristic submanifold theory of Jaco-Shalen
[\cite{JacoShalen}] and Johannson [\cite{Johannson}] that for each $\epsilon\in
\{\pm1\}$ there is an $(I, \partial I)$-bundle pair $(\Sigma^\epsilon,
\Phi^\epsilon)\subset (M_{\tilde F}^\epsilon, \tilde F)$,
well-defined up to ambient isotopy in $(M_{\tilde F}^\epsilon, \tilde
 F)$,
such that

\part{(1)}  the frontier of $\Sigma^\epsilon$ in $M_{\tilde F}^\epsilon$
consists of essential annuli in $(M_{\tilde F}^\epsilon,\tilde F)$;

\part{(2)} no component $(\sigma, \phi)$ of $(\Sigma^\epsilon,
\Phi^\epsilon)$ is
 homotopic (as a pair) into $(\Sigma^\epsilon - \sigma,
\Phi^\epsilon
- \phi)$; and

\part{(3)} if $K$ is a polyhedron and $H\colon (K\times I,K\times\bdry
I)\to(M_{\tilde F}^\epsilon,\tilde F)$ is an essential basic homotopy such that
$H_0:K\to \tilde F$ is a large map, then $H$ is homotopic as a map of
pairs to a homotopy whose image lies in $\Sigma^\epsilon$.

The $I$-fibers of the $I$-bundle structure on $\Sigma^\epsilon$ can be used
to build an essential basic homotopy of $\Phi^\epsilon$ in
$\Sigma^\epsilon$. Indeed, define a {\it fundamental homotopy} of
$\Phi^\epsilon$
to be any essential basic homotopy $$H_{\Sigma^\epsilon} \colon (\Phi^\epsilon
\times I, \Phi^\epsilon \times \bdry I) \to
(\Sigma^\epsilon, \Phi^\epsilon) \subset (M_{\tilde F}^\epsilon,
 \tilde F)$$
 satisfying the following conditions.

\part{(1)}The time-$0$ map of $H_{\Sigma^\epsilon}$ is the identity map of
$\Phi^\epsilon$.

\part{(2)} For each $x \in \Phi^\epsilon$, $H_{\Sigma^\epsilon}(\{x\}
\times I)$ is an
$I$-fiber of $\Sigma^\epsilon$.

\part{(3)} For each component $\phi$ of $\Phi^\epsilon$ which is contained in a
trivial $I$-bundle component $\sigma$ of $\Sigma^\epsilon$, the map
$H_{\Sigma^\epsilon}$ restricts to a homeomorphism of $\phi \times I$ onto
$\sigma$.

\part{(4)} For each component $\phi$ of $\Phi^\epsilon$ which is contained in
a non-trivial $I$-bundle component $\sigma$ of $\Sigma^\epsilon$, the map
$H_{\Sigma^\epsilon}$ restricts to a $2$-fold covering map from $\phi\times
I$ to $\sigma$.

It is clear that a fundamental homotopy exists and is unique up to
composition with a fiber-preserving homeomorphism of $\Sigma^\epsilon$
which restricts to the identity on $\Phi^\epsilon$. Note also that a
fundamental homotopy is a two-sheeted covering map from
$\Phi^\epsilon\times I$ to $\Sigma^\epsilon$. The involution
of $\Phi^\epsilon$ induced by its $\partial I$-bundle structure is the
restriction of the deck transformation of this two-sheeted covering
map.  This involution will be denoted by $\tau_\epsilon$.
Thus $$H_{\Sigma^\epsilon}(x,1) = \tau_\epsilon(x)$$
for every $x\in\Phi^\epsilon$.

The notation $\Sigma^\epsilon$, $\Phi^\epsilon$ and $\tau_\epsilon$
will be used throughout this section.

\Lemma
\tag{Orientation}
Let $\tilde F$ be given a consistent orientation.  Then
the embedding $\tau_\epsilon:\Phi^\epsilon\to\tilde F$ reverses
orientation. (See \xref{terminology}.)
\EndLemma
\Proof
By \xref{SplittingSurfaces}, the orientation of $\Phi^\epsilon$ is induced
from an orientation of $M_{\tilde F}^\epsilon$.  Hence, if we orient
$\Sigma^\epsilon$ by restricting the orientation of $M_{\tilde
F}^\epsilon$, then the orientation of $\Phi^\epsilon$ is induced by
the orientation of $\Sigma^\epsilon$.  The involution $\tau_\epsilon$
extends to an involution $\sigma$ of $\Sigma^\epsilon$ which leaves
each $I$-fiber invariant while reversing its orientation.  Since
$\sigma$ clearly reverses the orientation of $\Sigma^\epsilon$ it
follows that $\tau_\epsilon$ reverses the orientation of
$\Phi^\epsilon$.
\EndProof

\Definition
\tag{StandardDef}
An essential basic homotopy $H \colon (K \times I, K \times\bdry I)
\to (\Sigma^\epsilon, \Phi^\epsilon)$ is said to be {\it standard} if it is
of the form $H(x, t) = H_{\Sigma^\epsilon}(f(x), t)$ where
$H_{\Sigma^\epsilon}$ is a fundamental homotopy and $f$ is some map
from $ K$ to $\Phi^\epsilon$. Note that we then have $H_0=f$.  Note
also that if $f'\colon K \to
\Phi^\epsilon$ is homotopic to $f$ in $\Phi^\epsilon$ then the
standard homotopies $H(x, t) = H_{\Sigma^\epsilon}(f(x), t)$ and $H'(x, t) =
H_{\Sigma^\epsilon}(f'(x), t)$ are homotopic as maps of pairs $(K \times I,
K \times\bdry I) \to (\Sigma^\epsilon, \Phi^\epsilon)$.
\EndDefinition

\Lemma
\tag{Standard}
If $H \colon (K \times I, K \times \bdry I) \to
(M_{\tilde F}^\epsilon, \tilde F)$ is an
essential basic homotopy such that $H_0\colon K\to \tilde F$ is a large map,
then $H$ is homotopic as a map of pairs to a standard
essential basic homotopy.
\EndLemma

\Proof Appealing to the characteristic submanifold theory
[\cite{JacoShalen},\cite{Johannson}], we may assume, without loss of
generality, that $H(K \times I) \subset \Sigma^\epsilon$ and
$H(K\times \bdry I) \subset \Phi^\epsilon$.  It is not hard to see that
$H$ lifts to a map $\tilde H \colon (K \times I, K \times \bdry I) \to
(\Phi^\epsilon \times I, \Phi^\epsilon \times \{0, 1\})$ such that
$H=H_{\Sigma^\epsilon} \circ \tilde H$.  Write $\tilde H(x, t) =
(H'(x,t), T(x, t)) \in \Phi^\epsilon \times I$ and observe that
$T(x,0) = 0$ for each $x \in K$, while the essentiality of $H$ implies that
$ T(x,1) = 1$ for each $x$.

Let $J \colon  ((K \times  I) \times I, (K \times  I) \times \{0,1\}) \to
(\Sigma^\epsilon, \Phi^\epsilon)$ be given by
$$J((x, t), s) = H_{\Sigma^\epsilon}(H'(x,(1-s)t), (1-s)T(x,t) + st).$$
Then
$$\eqalign{
J((x, t), 0) &= (H_{\Sigma^\epsilon} \circ \tilde H)(x, t) = H(x, t)\hbox{,
and}\cr
J((x, t), 1) &= H_{\Sigma^\epsilon}(H'(x,0), t)= H_{\Sigma^\epsilon}(H_0(x),t),
}$$
while for each $x \in K$,
$$\eqalign{
J((x, 0), s) &= H(x, 0) \in \Phi^\epsilon \hbox{, and}\cr
J((x, 1), s) &= H_{\Sigma^\epsilon}(H'(x, 1-s), 1) \in \Phi^\epsilon.
}$$
Thus $H=J_0$ is homotopic as a map of pairs to the standard essential
basic homotopy $J_1$.
\EndProof

\paragraph
\tag{DefinitionOfPhiOne}
We define $\Phi_1^\epsilon={\cal L}(\Phi^\epsilon) \subset \tilde F$.
Note that the free involution $\tau_\epsilon$ restricts to a free
involution of $\Phi_1^\epsilon$.  We shall denote this restriction
by $\tau_\epsilon$ as well.

\Lemma
\tag{PhiOne}
The surface $\Phi_1^\epsilon$ has the following property.

For any large map $f \colon K \to \tilde F$, there exists an essential
basic homotopy $H$ in the pair $(M,\tilde F)$ on the $\epsilon$-side with
$H_0 =
f$ if and only if $f$ is homotopic in $\tilde F$ to a map with image in
$\Phi_1^\epsilon$.

Furthermore, any large subsurface of $F$ with this
property is isotopic to $\Phi_1^\epsilon$.
\EndLemma

\Proof
If $f \colon K \to \tilde F$ is a large map homotopic in $\tilde F$ to
some $f'$ for which $f'(K) \subset \Phi_1^\epsilon$ and $H_{\Sigma^\epsilon}
\colon (\Phi^\epsilon\times [0,1], \Phi^\epsilon \times \{0,1\}) \to
(\Sigma^\epsilon,
\Phi^\epsilon)$  is a fundamental homotopy, then $H \colon (K \times [0,1],
K \times
\{0, 1 \}) \to (\Sigma^\epsilon,\Phi^\epsilon) \subset(M_{\tilde F}^\epsilon, \tilde
F)$ defined
by
$$H(x, t) = H_{\Sigma^\epsilon}(f'(x), t)$$ is an essential basic
homotopy in $(M, \tilde F)$ on the $\epsilon$-side, with $H_0=f'$.
The desired essential basic homotopy with time-$0$ map $f$ is now
readily constructed.

Conversely let $H \colon (K \times [0,1], K \times \{0, 1\}) \to
(M,\tilde F)$ be an essential basic homotopy in $(M, \tilde F)$ on the
$\epsilon$-side with $H_0=f$.  According to \xref{HomotopyDef}, any
essential basic homotopy in $(M,\tilde F)$ which starts on the
$\epsilon$-side corresponds to an essential homotopy in the pair
$(M_{\tilde F}^\epsilon, \tilde F)$.  So it follows from Lemma
\xref{Standard} that $H$ is homotopic as a map of pairs to a
standard essential basic homotopy $H'$.  It follows from \xref{JSJStuff} and
\xref{StandardDef} that $H'_t(K) \subset \Phi^\epsilon$ for
$t=0,1$, that $H'_0$ is homotopic to $f$ and that $H'_1 =\t_\epsilon \circ
H'_0$.  Since $f$, and hence $H'_0$, is large, the image of $H'_0$ is
contained in $\Phi_1^\epsilon$.

Finally suppose that $\Psi_1^\epsilon \subset \tilde F$ is another large
subsurface of $\tilde F$ which has the stated property.
Then the inclusion $\Phi_1^\epsilon\to \tilde F$ is homotopic into
$\Psi_1^\epsilon$ and
vice versa. Thus by Lemma \xref{SurfaceTwo}, $\Psi_1^\epsilon$ is isotopic
to $\Phi_1^\epsilon$ in $\tilde F$.  \EndProof

\Lemma
\tag{HSigma}
Let $H \colon (K \times I, K \times \bdry I) \to (M, \tilde F)$ be an
essential basic homotopy on the $\epsilon$-side such that $H_0$ is a large
map from $K$ to $\tilde F$.  If $f\colon K \to \Phi_1^\epsilon$ is any map
which is homotopic in $\tilde F$ to $H_0$, then $H_1$ is homotopic in
$\tilde F$ to
$\t_\epsilon\circ f$.
\EndLemma

\Proof
By Lemma \xref{Standard} the homotopy $H$ is homotopic as a map of
pairs to a standard essential basic homotopy $H'$. By
Proposition \xref{Containment} $H_0'$ is homotopic to $f$ in
$\Phi_1^\epsilon$ and so by the remark at
the end of Definition \xref{StandardDef} we may assume that $H'_0 =
f$. Thus $H_1$ is homotopic in $\tilde F$ to $H'_1 = \tau_\epsilon\circ f$.
\EndProof

Our next goal is to extend Lemma \xref{PhiOne} to reduced homotopies
of length $n$ in $(M, \tilde F)$.

\Proposition
\tag{Surfaces}
For each fixed $\epsilon\in\{\pm 1\}$, there is a sequence of large
(possibly empty) subsurfaces $(\Phi_k^\epsilon)_{k
\geq 0}$ of $\tilde F$, such that
$\Phi_0^\epsilon = \tilde F$, $\Phi_1^\epsilon$ is the surface defined
in Subsection \xref{DefinitionOfPhiOne}, and for each $k \geq 0$ we have:

\part {$(1)$} 
$\Phi_k^\epsilon \supset \Phi_{k+1}^\epsilon$; and

\part {$(2)$} 
a large map $f \colon K \to \tilde F$ is homotopic in $\tilde F$ to a
map with image in $\Phi_k^\epsilon$ if and only if there exists a
reduced homotopy $H$ of length $k$ starting on the $\epsilon$-side
with $H_0 = f$.

Furthermore condition $(2)$ determines $\Phi_k^\epsilon$, up to
isotopy, among the class of large subsurfaces of $\tilde F$.
\EndProposition

\Proof
We construct the surfaces inductively in such a way that (1) and (2)
hold.  Set $\Phi_0^\epsilon=\tilde F$ and let $\Phi_1^\epsilon$ be
defined as in \xref{DefinitionOfPhiOne}.  According to Lemma
\xref{PhiOne}, condition $(2)$ holds for $k = 1$.  Let $m\geq 2$ be
given, and suppose that for $\epsilon=\pm1$ we have defined large
subsurfaces
$$\tilde F=\Phi_0^\epsilon \supset \Phi_1^\epsilon \supset
\Phi_2^\epsilon\supset
\cdots \supset \Phi_{m-1}^\epsilon$$
such that condition $(2)$ holds for $k< m$ and for $\epsilon=\pm1$.
As we observed in \xref{DefinitionOfEssentialIntersection}, there is a
surface $A^\epsilon_m\subset\Phi_1^\epsilon$ which is isotopic to the
large intersection $\Phi_1^\epsilon\wedge_{\cal L}
\Phi_{m-1}^{-\epsilon}$. For $\epsilon=\pm1$ we define
$$\Psi_m^\epsilon =\t_\epsilon(A^\epsilon_m).$$ 
Note that by \xref{DefinitionOfPhiOne}, $\Phi_1^\epsilon$ is invariant under
the map $\t_\epsilon$ and so we have
$$\Psi_m^\epsilon \subset \Phi_1^\epsilon.$$

For $\epsilon=\pm1$, we claim that:
\part{$(*)$} a large
map $f \colon K \to \tilde F$ is homotopic in $\tilde F$ to a map with
image in $\Psi_m^\epsilon$ if and only if there exists a reduced
homotopy $H$ of length $m$ starting on the $\epsilon$-side with $H_0 =
f$.

To prove this, we first consider a large map $f \colon K \to \tilde F$
and assume that there exists a reduced homotopy $H \colon K\times I
\to (M, \tilde F)$ of length $m$ starting on the $\epsilon$-side with
$H_0 = f$.  Write $H$ as a composition of $m$ essential basic
homotopies $H^1, \ldots, H^m$ where $H^i$ starts on the
$(-1)^{i-1}\epsilon$ side for $i=1,\ldots,m$.  Define $f_1\colon K \to
\tilde F$ by $f_1 = H^1_1 = H^2_0$ and let $H'$ be the composition of
$H^2, \ldots, H^m$.  Then $f_1$ is large and $H'$ is a reduced
homotopy of length $m-1$ starting on the $-\epsilon$-side, with
$H'_0=f_1$. Our inductive hypothesis implies that $f_1$ is homotopic
to a map with image in $\Phi_{m-1}^{-\epsilon}$.  On the other hand,
$H^1$ is an essential basic homotopy on the $\epsilon$-side with
$H^1_0=f$ and $H_1^1=f_1$; hence by Lemma \xref{HSigma} we see that
$f_1$ is homotopic in $\tilde F$ to a map $f_1'$ whose image lies in
$\Phi_1^\epsilon$, and that $f$ is homotopic in $\tilde F$ to
$\tau_\epsilon\circ f_1'$.  Proposition
\xref{EssInt} now implies that $f_1'$ is homotopic in $\tilde F$ to a
map $f_1''$ whose image lies in $A^\epsilon_m$.  Proposition
\xref{Containment} implies that $f_1'$ is homotopic to $f_1''$ in
$\Phi_1^\epsilon$.  Composing this homotopy with $\tau_\epsilon$ we
see that $f$ is homotopic to a map with image in
$\t_\epsilon(A^\epsilon_m) = \Psi_m^\epsilon$.

To prove the converse observe that the fundamental homotopy
$$H = H_{\Sigma^\epsilon}\colon (\Phi^\epsilon \times [0, 1], \Phi^\epsilon
\times \{0,
1\}) \to (\Sigma^\epsilon, \Phi^\epsilon)$$
is an essential basic homotopy $H$ in
$(M, \tilde F)$ on the $\epsilon$-side such that $H_1$ is the inclusion
$i_0\colon
\Phi_1^\epsilon \to \tilde F$, and $H_1=\t_\epsilon| \Phi_1^\epsilon$.  Since
$\Psi_m^\epsilon \subset \Phi_1^\epsilon$, we may restrict $H$ to
$\Psi_m^\epsilon \times I$ to obtain an essential basic homotopy $H^1$
in $(M, \tilde F)$ on the $\epsilon$-side.  The time-$0$ map of $H^1$
is the inclusion $\Psi_m^\epsilon\to \tilde F$ and the time-$1$ map of
$H^1$ is a homeomorphism $\psi\colon \Psi_m^\epsilon \to
A^\epsilon_m$. By Proposition \xref{EssInt} and Definition
\xref{DefinitionOfEssentialIntersection}, $A_m^\epsilon$ is isotopic
to a subsurface of $\Phi_{m-1}^{-\epsilon}$. Since our induction
hypothesis implies that the inclusion $\Phi_{m-1}^{-\epsilon} \to
\tilde F$ is the time-$0$ map of a reduced homotopy of length $m-1$
starting on the $-\epsilon$ side, it follows that $\psi$ is also the
time-$0$ map of a reduced homotopy $H'$ in $(M, \tilde F)$ of length
$m-1$ starting on the $-\epsilon$-side.  The composition of the two
homotopies $H_1$ and $H'$ is a reduced homotopy of length $m$ starting
on the $\epsilon$-side whose time-$0$ map is the inclusion
$\Psi_m^\epsilon \to \tilde F$.  This completes the proof of $(*)$.

It follows from $(*)$ that there is a reduced homotopy of length $m$
starting on the $\epsilon$-side with time-$0$ map equal to the
inclusion map of $\Psi_m^\epsilon$ into $\tilde F$.  Clearly there
also exists a reduced homotopy of length $m-1$ starting on the same
side and having the same time-$0$ map.  Thus we may apply Property (2)
again to conclude that the inclusion map of $\Psi_m^\epsilon$ is
homotopic to a map with image in $\Phi_{m-1}^\epsilon$.  By Lemma
\xref{SurfaceTwo}, $\Psi_m^\epsilon$ is isotopic to a subsurface
$\Phi_m^\epsilon$ of $\Phi_{m-1}^\epsilon$. It follows from $(*)$ that
condition (2) holds for $k=m$. The induction is now complete.

To prove the last assertion of the proposition, suppose that
$\Psi_k^\epsilon
\subset \tilde F$ is another large surface in $\tilde F$
which
satisfies condition (2) of the theorem.  Taking $f$ to be the
inclusion map of either surface into $\tilde F$, we see that
$\Phi_k^\epsilon$ is homotopic into $\Psi_k^\epsilon$ and vice
versa. Thus by Lemma \xref{SurfaceTwo}, $\Psi_k^\epsilon$ is isotopic
to $\Phi_k^\epsilon$ in $\tilde F$.
\EndProof

\subsection{Time-1 maps of reduced homotopies}
\tag{GeneralBound}

Throughout this section $M$ will denote a simple knot manifold and
$\tilde F$ will denote a splitting surface in $M$.  We will define
$\Phi^\epsilon$ and $\tau_\epsilon$ as in Subsection \xref{EssHom}.
We shall also fix subsurfaces $\Phi_k^\epsilon$ of $\tilde F$ for
which the conclusions of Proposition \xref{Surfaces} hold. A
crucial ingredient in obtaining a bound of the type given by Theorem
\xref{LengthBound} is provided by Proposition \xref{StrictContainment}, which
asserts that if the inclusion $\Phi_{k+2}^\epsilon \subseteq
\Phi_k^\epsilon$ is a homotopy equivalence, then $\tilde F$ is a
semi-fiber. The estimate is strengthened by using Proposition
\xref{OddK}, which asserts that for each odd $k$, $\Phi_k^\epsilon$
admits a fixed-point free involution.

\Proposition
\tag{TimeOneMaps}
For each $\epsilon \in \{\pm1\}$ and each $k\ge0$ there exists a
homeomorphism $h_k^\epsilon\colon \Phi_k^\epsilon\to
\Phi_k^{(-1)^{k+1}\epsilon}$ which has the following property:
\part{($\ast$)}For any reduced homotopy $H$ of length $k$ with domain $K$
which starts on the $\epsilon$ side and has a large
time-$0$ map, there exists a map $f\colon K \to \Phi_k^\epsilon$
such that $H_0$ is homotopic in $\tilde F$ to $f$
and $H_1$ is homotopic in $\tilde F$ to $h_k^\epsilon\circ f$.

The homeomorphism $h_k^\epsilon$ is unique up to isotopy.  Furthermore,
if $\tilde F$ is given a consistent orientation then the embedding
$h_k^\epsilon\colon \Phi_k^\epsilon\to\tilde F$ reverses orientation
if $k$ is odd and preserves orientation if $k$ is even.  (See
\xref{terminology}.)  We may take $h_1^\epsilon$ to be the map
$\tau_\epsilon|\Phi_1^\epsilon$ described in \xref{DefinitionOfPhiOne}.
\EndProposition

\Proof 
We construct the maps $h_k^\epsilon$ by induction on $k$.
For each $\epsilon \in\{\pm1\}$ define $h_0^\epsilon$ to be the
identity map of $\tilde F$ and $h_1^\epsilon$ to be the map $\tau_\epsilon$.
It is clear that property ($\ast$) holds for $h_0^\epsilon$ and Lemmas \xref{PhiOne}
and \xref{HSigma} imply that it holds for $h_1^\epsilon$.

Suppose that for $\epsilon=\pm1$ we have constructed a homeomorphism
$h_k^\epsilon\colon \Phi_k^\epsilon\to \Phi_k^{(-1)^{k+1}\epsilon}$
with property ($\ast$).  In particular, since $\Phi_k^\epsilon$ is a
$\pi_1$-injective subsurface, $h_k^\epsilon$ is a $\pi_1$-injective
map.  By Proposition \xref{Surfaces}, for $\epsilon=\pm1$, we may choose a
reduced homotopy ${\cal H}$ of length $k+1$ starting on the $\epsilon$
side such that ${\cal H}_0$ is the inclusion $\Phi_{k+1}^\epsilon\to
\tilde F$.  Write ${\cal H}$ as the composition of a reduced homotopy ${\cal
H}'$ of length $k$ and an essential basic homotopy ${\cal H}''$.  The
induction hypothesis implies that the time-$1$ map of ${\cal H}'$ is
homotopic to the embedding $h_k^\epsilon|\Phi_{k+1}^\epsilon$.  By
reversing the time variable of ${\cal H}'$ and applying Proposition
\xref{Surfaces} we see that $h_k^\epsilon|\Phi_{k+1}^\epsilon$ is
homotopic to an embedding $h
\colon \Phi_{k+1}^\epsilon \to \Phi_k^{(-1)^{k+1}\epsilon} \subset
\Phi_1^{(-1)^{k+1}\epsilon}$.  We define
$\theta^\epsilon\colon\Phi_{k+1}^\epsilon \to \tilde F$ to be
$\tau_{(-1)^{k}\epsilon}\circ h$.  Note that $h$ is $\pi_1$-injective
since $\Phi_{k+1}^\epsilon$ is a $\pi_1$-injective subsurface and
since $h_k^\epsilon$ is a $\pi_1$-injective map.
Furthermore since $\tau_{(-1)^{k+1}\epsilon}$ is an involution
of a $\pi_1$-injective subsurface, it follows that $\theta^\epsilon$
is a $\pi_1$-injective map for $\epsilon=\pm1$.

Since $\tau_{(-1)^{k}\epsilon}$ is an involution of
$\Phi_1^{(-1)^{k}\epsilon}$ we know that $\theta^\epsilon$ is the
time-$0$ map of a standard basic essential homotopy whose time-$1$ map
is $\tau_{(-1)^{k}\epsilon}\circ \theta^\epsilon = h$.  Since $h$ is
in turn the time-$0$ map of a reduced homotopy of length $k$ starting
on the $(-1)^{k+1}\epsilon$-side, it follows that $\theta^\epsilon$ is
the time-$0$ map of a length $k+1$ homotopy starting on the
$(-1)^k\epsilon$-side. It thus follows
from Proposition \xref{Surfaces} that the embedding
$\theta^\epsilon$ is homotopic into $\Phi_{k+1}^{(-1)^k\epsilon}$,
i.e. that the surface $\theta^\epsilon(\Phi_{k+1}^\epsilon)$ is
homotopic into $\Phi_{k+1}^{(-1)^k\epsilon}$.

Replacing $\epsilon$ by $(-1)^k\epsilon$ we find that
$\theta^{(-1)^k\epsilon}(\Phi_{k+1}^{(-1)^k\epsilon})$ is homotopic
into $\Phi_{k+1}^\epsilon$.  By Lemma \xref{SurfaceTwo}(1),
$\theta^{(-1)^k\epsilon}(\Phi_{k+1}^{(-1)^k\epsilon})$ is isotopic to
a subsurface of $\Phi_{k+1}^\epsilon$, which is $\pi_1$-injective
since $\theta^{(-1)^k\epsilon}$ is a $\pi_1$-injective map.  In
particular $\Phi_{k+1}^{(-1)^k\epsilon}$ is homeomorphic to a
$\pi_1$-injective subsurface of $\Phi_{k+1}^\epsilon$, and hence to a
$\pi_1$-injective subsurface of
$\theta^\epsilon(\Phi_{k+1}^\epsilon)$.  Apply Lemma
\xref{SurfaceTwo}(2), taking $A =
\theta^\epsilon(\Phi_{k+1}^\epsilon)$ and $B =
\Phi_{k+1}^{(-1)^k\epsilon}$.
It follows that $A$ and $B$ are isotopic and hence
that $\theta^\epsilon$ is isotopic in $\tilde F$ to a homeomorphism
$h_{k+1}^\epsilon \colon
\Phi_{k+1}^\epsilon \to
\Phi_{k+1}^{(-1)^k\epsilon}$.

We now show that $h_{k+1}^\epsilon$ has property ($\ast$).  Let $H$ be any
reduced homotopy of length $k+1$ starting on the $\epsilon$ side such
that $H_0$ is a large map.  By \xref{Surfaces}, $H_0$ is homotopic in
$\tilde F$ to a map $f: K \to \Phi_{k+1}^\epsilon$.  Write $H$ as the
composition of a reduced homotopy $H'$ of length $k$ and a basic
essential homotopy $H''$.  Since $h_k^\epsilon$ has property ($\ast$), the
map $H''_0 = H'_1$ is homotopic in $\tilde F$ to $h_k^\epsilon\circ
f$, which is in turn homotopic to $h\circ f$.  Then Lemma
\xref{HSigma} implies that $H_1 = H''_1$ is homotopic to
$\tau_{(-1)^{k}\epsilon}\circ h\circ f =
\theta^\epsilon\circ f$, and hence to $h_{k+1}^\epsilon\circ f$.  This
establishes ($\ast$) and completes the inductive definition of the
$h_j^\epsilon$.

By Lemma \xref{Orientation}, if $\tilde F$ is given a consistent
orientation, the involutions $\tau_{\pm1}$ are orientation reversing
embeddings of $\phi_{\pm1}$ into $\tilde F$.  It follows from the
inductive construction that $h_k^{\pm 1}$ reverses orientation if $k$
is odd and preserves orientation if $k$ is even.

It remains to prove that a homeomorphism satisfying ($\ast$) is unique up
to isotopy. Suppose that
$h,h':\Phi_k^\epsilon\to\Phi_k^{(-1)^{k+1}\epsilon}$ both satisfy
($\ast$). We apply Proposition \xref{Surfaces}, taking $K=\Phi_k^\epsilon$
and taking $f:\Phi_k^\epsilon\to\tilde F$ to be the inclusion. This
gives a length-$k$ reduced homotopy $H$ starting on the $\epsilon$
side with time-$0$ map $f$. By property ($\ast$), $h$ and $h'$ are both
homotopic in $\tilde F$ to $H_1$, and hence to each other. By
\xref{Containment}, $h$ and $h'$ are homotopic as maps
from $h,h':\Phi_k^\epsilon$ to $\Phi_k^{(-1)^{k+1}\epsilon}$.  It
follows from Theorem 6.4 and Theorem A.4 of [\cite{Epstein}] that they are
isotopic.
\EndProof

For the rest of this section, it will be understood that for each $k$
and each $\epsilon$ we have fixed homeomorphisms
$h_k^\epsilon:\Phi_k^\epsilon\to\Phi_k^{(-1)^{k+1}\epsilon}$ for which
condition ($\ast$) of Proposition \xref{TimeOneMaps} holds.

\Lemma\tag{BetterTimeOne}
Let $\epsilon \in \{\pm1\}$ and $k\ge0$ be given.  Suppose that $H$ is
a reduced homotopy in $(M,\tilde F)$ of length $k$ with domain $K$
which starts on the $\epsilon$ side.  Supppose that $f=H_0$ is large
and that $f(K)\subset\Phi_k^\epsilon$.  Then $H_1$ is homotopic in
$\tilde F$ to $h_k^\epsilon\circ f$.
\EndProposition

\Proof
By condition ($\ast$) of \xref{TimeOneMaps} there is a map $f':K\to
\Phi_k^\epsilon$ such that $H_0 = f$ is homotopic in $\tilde F$ to $f'$ and
such that $H_1$ is homotopic in  $\tilde F$ to $h_k^\epsilon\circ f'$.
It follows from Proposition \xref{Containment} that $f$ is homotopic
to $f'$ in $\Phi_k^\epsilon$.  Thus $h_k^\epsilon\circ f$ is homotopic
in $\tilde F$ to $h_k^\epsilon\circ f'$, and hence to $H_1$.
\EndProof 

\Proposition
\tag{Inverses}
Let $\epsilon \in \{\pm1\}$ and $k\ge0$ be given.
Then the homeomorphism
$h_k^{(-1)^{k+1}\epsilon}:\Phi_k^{(-1)^{k+1}\epsilon}\to\Phi_k^\epsilon$
is isotopic to the inverse of
$h_k^\epsilon:\Phi_k^\epsilon\to\Phi_k^{(-1)^{k+1}\epsilon}$.
\EndProposition

\Proof
It follows from Proposition \xref{Surfaces} that there is a reduced
homotopy $H$ of length $k$ starting on the $\epsilon$-side such that
$H_0$ is the inclusion $\iota:\Phi_k^\epsilon\to \tilde F$.  By Lemma
\xref{BetterTimeOne} we have that $H_1$ is homotopic in $\tilde F$ to
$h_k^\epsilon$.  Applying Lemma \xref{BetterTimeOne} to the homotopy
$H'$ obtained by reversing the time variable of $H$, we see that
$\iota = H_0 = H'_1$ is homotopic in $\tilde F$ to the composition
$h_k^{(-1)^{k+1}\epsilon}\circ H'_0$.  Since $H'_0 = H_1 \sim
h_k^\epsilon$ we have that $h_k^{(-1)^{k+1}\epsilon}\circ
h_k^\epsilon$ is homotopic in $\tilde F$ to the inclusion $\iota$.  It
now follows from Proposition \xref{Containment} that the
self-homeomorphism $h_k^{(-1)^{k+1}\epsilon}\circ h_k^\epsilon$ of
$\Phi_k^\epsilon$ is homotopic to the identity in $\Phi_k^\epsilon$,
and therefore isotopic to the identity by Theorem 6.4 and Theorem A.4
of [\cite{Epstein}].
\EndProof

\Proposition\tag{PropertyThree}
Let $i$ and $j$ be non-negative integers, and set $k=i+j$. Then for
each $\epsilon\in\{\pm1\}$, the map $h_i^\epsilon|\Phi_k^\epsilon$ is
homotopic in $\tilde F$  to an embedding $g_i^\epsilon\colon \Phi_k^\epsilon\to
\Phi_{j}^{(-1)^{i}\epsilon}$ such that
$h_{j}^{(-1)^{i}\epsilon}\circ g_i^\epsilon$ is homotopic in $\tilde
F$ to $h_k^\epsilon$.
\EndProposition

\Proof Let $H$ be a reduced homotopy of length $k$ starting on the
$\epsilon$ side such that $H_0$ is the inclusion
$\Phi_k^\epsilon\to\tilde F$.  Write $H$ as the composition of a
reduced homotopy $H'$ of length $i$ starting on the $\epsilon$ side
and a reduced homotopy $H''$ of length $j$ starting on the
$(-1)^i\epsilon$ side.  Applying Lemma \xref{BetterTimeOne}, with the
roles of $k$, $H$ and $f$ played respectively by $i$, $H'$ and the
inclusion map $\Phi_k^\epsilon \to \Phi_i^\epsilon$, we find that
$H'_1 = H_0''$ is homotopic to the embedding
$h_i^\epsilon|\Phi_k^\epsilon$.

On the other hand it follows from Proposition \xref{TimeOneMaps}
that $H_0''$ is homotopic in $\tilde F$ to a map $g_i^\epsilon\colon
\Phi_k^\epsilon\to \Phi_{j}^{(-1)^{i}\epsilon}$.  Since $g_i^\epsilon$
is homotopic to the embedding $h_i^\epsilon|\Phi_k^\epsilon$ it
follows from part (1) of Lemma \xref{SurfaceTwo} that we may take
$g_i^\epsilon$ to be an embedding.  After modifying the homotopy
$H$ we may assume that $H'_1 = H''_0 = g_i^\epsilon$.

Applying Lemma \xref{BetterTimeOne} again, with $j$, $H''$ and
$g_i^\epsilon$ playing the roles of $k$, $H$ and $f$, we
conclude that $H_1 = H_1''$ is homotopic in $\tilde F$ to 
$h_j^{(-1)^{i}\epsilon}\circ g_i^\epsilon$.  

Finally, applying Lemma \xref{BetterTimeOne} directly to the homotopy
$H$, we see that $H_1$ is homotopic in $\tilde F$ to $h_k^\epsilon$.
The conclusion of the Proposition follows.
\EndProof

\Proposition\tag{PropertyTwo}
Let $i$ and $j$ be non-negative integers, and set $k=i+j$.  Then
for each $\epsilon \in \{\pm1\}$
the subsurface $h_i^\epsilon(\Phi_k^\epsilon)$ is isotopic in $\tilde F$ to
$\Phi_{i}^{(-1)^{i+1}\epsilon}\wedge_{\cal L}  \Phi_{j}^{(-1)^i\epsilon}$.
\EndProposition

\Proof
We
have $h_i^\epsilon(\Phi^\epsilon_k)\subset
h_i^\epsilon(\Phi_i^\epsilon)=\Phi_i^{(-1)^{i+1}\epsilon}$.
On the other hand it follows from Proposition \xref{PropertyThree}
that $h_i^\epsilon(\Phi^\epsilon_k)$ is isotopic in $\tilde F$ to
a subsurface of  $\Phi_{j}^{(-1)^i\epsilon}$.
Hence by \xref{DefinitionOfEssentialIntersection} the subsurface 
$h_i^\epsilon(\Phi^\epsilon_k)$ is isotopic to a subsurface of
$\Phi_{i}^{(-1)^{i+1}\epsilon}\wedge_{\cal L}  \Phi_{j}^{(-1)^i\epsilon}$.

To complete the proof of the proposition, it now suffices by Lemma
\xref{SurfaceTwo}(2) to show that the large surface
$\Phi_i^{(-1)^{i+1}\epsilon}\wedge_{\cal
L}\Phi_{j}^{(-1)^i\epsilon}$ is homeomorphic to a $\pi_1$-injective
subsurface of $\Phi_k^\epsilon$ and hence of
$h_i^\epsilon(\Phi_k^\epsilon)$.  To prove this, note that
$\Phi_i^{(-1)^{i+1}\epsilon}\wedge_{\cal L}
\Phi_{j}^{(-1)^i\epsilon}$ is isotopic to a subsurface $A$ of
$\Phi_i^{(-1)^{i+1}\epsilon}$.  By Proposition \xref{Surfaces} there is a
reduced homotopy of length $i$ whose time-$0$ map is the inclusion $A
\to \tilde F$ and whose time-$1$ map is homotopic to
$h_i^{(-1)^{i+1}\epsilon}\colon A\to \tilde F$.  By reversing the time
variable we obtain a homotopy $H$ such that $H_0$ is homotopic in
$\tilde F$ to $h_i^{(-1)^{i+1}\epsilon}\colon A\to \tilde F$, and
$H_1$ is the inclusion $A\to \tilde F$.  Since $A$ is isotopic to a
subsurface of $\Phi_{j}^{(-1)^i\epsilon}$ there is also a reduced
homotopy $H''$ of length $j$ starting on the ${(-1)^i\epsilon}$ side
whose time-$0$ map is the inclusion $A\to \tilde F$.  The composition
of $H$ and $H'$ is a reduced homotopy of length $k$ starting on the
$\epsilon$ side and having time-$0$ map homotopic to
$h_i^{(-1)^{i+1}\epsilon}|A \colon A \to \tilde F$.  In particular, by
Proposition \xref{Surfaces} and Lemma
\xref{SurfaceTwo}, $h_i^{(-1)^{i+1}\epsilon}(A)$ is isotopic to a
subsurface of $\Phi_k^\epsilon$.
\EndProof

One of the main results of this section is that $\Phi_{2k+1}^\epsilon$ admits a
free involution. The proof of this fact is based on our next lemma.

\Lemma
\tag{InvariantUpToIsotopy}
The subsurface $h_k^\epsilon(\Phi_{2k+1}^\epsilon)$ is isotopic in
$\tilde F$ to a subsurface $A$ of $\Phi_1^{(-1)^k\epsilon}$ with the
property that $\tau_{(-1)^k\epsilon}(A)$ is isotopic to $A$ in
$\Phi_1^{(-1)^k\epsilon}$.
\EndLemma

\Proof
By Proposition \xref{PropertyTwo} we have that
$h_k^\epsilon(\Phi_{2k+1}^\epsilon)$ is isotopic to
$\Phi_k^{(-1)^{k+1}\epsilon}\wedge_{\cal
L}\Phi_{k+1}^{(-1)^k\epsilon}$.  
In particular, $h_k^\epsilon(\Phi_{2k+1}^\epsilon)$ is isotopic
to a subsurface $A$ of $\Phi_{k+1}^{(-1)^k\epsilon}\subset
\Phi_1^{(-1)^k\epsilon}$.
Since $\tau_{(-1)^k\epsilon}(A) = h_1^{(-1)^k\epsilon}(A) \subset
h_1^{(-1)^k\epsilon}(\Phi_{k+1}^{(-1)^k\epsilon})$,
Proposition \xref{PropertyThree} implies that $\tau_{(-1)^k\epsilon}(A)$
is isotopic to a subsurface of $\Phi_k^{(-1)^{k+1}\epsilon}$.
On the other hand, since $A$ is isotopic to a subsurface of
$\Phi_k^{(-1)^{k+1}\epsilon}$ we know that the inclusion $A\to \tilde F$
is the time-0 map of a reduced homotopy $H''$ of length $k$ starting on
the ${(-1)^{k+1}\epsilon}$ side.  Since $A\subset
\Phi_1^{(-1)^{k+1}\epsilon}$, by reversing the time variable in the
basic homotopy provided by Lemma \xref{HSigma} we obtain a basic
homotopy $H'$ whose time-$0$ map is the inclusion
$\tau_{(-1)^k\epsilon}|A$ and whose time-$1$ map is the
inclusion $A\to \tilde F$.  The composition of $H'$ and $H''$ is a reduced
homotopy of length $k+1$ starting on the ${(-1)^k\epsilon}$ side
whose time-$0$ map is $\tau_{(-1)^k\epsilon}|A$.
It follows from Proposition \xref{Surfaces} that
$\tau_{(-1)^k\epsilon}(A)$ is homotopic into
$\Phi_{k+1}^{(-1)^k\epsilon}$.  By \xref{DefinitionOfEssentialIntersection}
we have that $\tau_{(-1)^k\epsilon}(A)$ is homotopic into
$\Phi_k^{(-1)^{k+1}\epsilon}\wedge_{\cal
L}\Phi_{k+1}^{(-1)^k\epsilon}$ which in turn is isotopic to $A$ in
$\tilde F$.
By Lemma \xref{SurfaceTwo}(2) and Proposition \xref{Containment},
$\tau_{(-1)^k\epsilon}(A)$ is isotopic to $A$ in $\Phi_1^{(-1)^k\epsilon}$.
\EndProof

\Proposition
\tag{OddK}
The surface $h^\epsilon_k(\Phi_{2k+1}^\epsilon)$ is
isotopic in $\tilde F$ to a subsurface of $\Phi_1^{(-1)^k\epsilon}$ which is
invariant under the free involution $\tau_{(-1)^k\epsilon}$.  In
particular, $\Phi_{2k+1}^\epsilon$ admits a free involution.
\EndProposition

\Proof
This is an immediate consequence of \xref{InvariantUpToIsotopy} and
\xref{Invariant}.
\EndProof

\Corollary
\tag{EvenChar}
For each odd integer $k>0$, the Euler characteristic
$\chi(\Phi_{k}^\epsilon)$ is even.
\EndCorollary
\NoProof

\Proposition
\tag{StrictContainment}
If $\Phi_k^\epsilon$ and $\Phi_{k+2}^\epsilon$ are isotopic in $\tilde
F$ for a given $k\ge0$ then either $\Phi_k^\epsilon = \emptyset$ or
$\tilde F$ is a semi-fiber.
\EndProposition

\Proof
In this proof all isotopies will be understood to take place in
$\tilde F$ unless specified otherwise.

First we claim that ($1$) $\Phi_k^\epsilon$ is isotopic to
$\Phi_m^\epsilon$ for all $m \ge k$.  Since $\Phi_{m-1}^\epsilon
\supset\Phi_{m}^\epsilon\supset\Phi_{m+1}^\epsilon$, and
$\Phi_{m}^\epsilon$ is $\pi_1$-injective, it suffices to
consider the case where $m-k$ is even.
The case $m=k+2$ holds by hypothesis.
Thus we need only show that if $m-k \ge 2$ is even and
$\Phi_{m}^\epsilon$ is isotopic to $\Phi_{m-2}^\epsilon$ then
$\Phi_{m+2}^\epsilon$ is isotopic to $\Phi_{m}^\epsilon$.
By Lemma \xref{PropertyTwo} we have that
$h_2^\epsilon(\Phi_{m+2}^\epsilon)$  and
$h_2^\epsilon(\Phi_{m}^\epsilon)$ are respectively
isotopic to $\Phi_2^{-\epsilon}
\wedge_{\cal L} \Phi_{m}^\epsilon$ and
$\Phi_2^{-\epsilon} \wedge_{\cal L} \Phi_{m-2}^\epsilon$.  These two
surfaces are isotopic by the induction hypothesis.  It follows that
$\Phi_{m+2}^\epsilon$ is isotopic to $\Phi_{m}^\epsilon$.  Claim ($1$)
follows.

Next we claim that ($1'$) $\Phi_m^{-\epsilon}$ is isotopic to
$\Phi_{k+1}^{-\epsilon}$ for all $m \ge k+1$.  By Lemma
\xref{PropertyTwo} we have for any $m > k+1$ that
$h_1^{-\epsilon}(\Phi_{m+1}^{-\epsilon})$ is isotopic to
$\Phi_1^{-\epsilon} \wedge_{\cal L} \Phi_m^\epsilon$, while
$h_1^{-\epsilon}(\Phi_{k+1}^{-\epsilon})$ is isotopic to
$\Phi_1^{-\epsilon} \wedge_{\cal L} \Phi_k^\epsilon$.
These two surfaces are
isotopic by ($1$), so ($1'$) follows.

Next we claim that ($2$) $h_j^\epsilon(\Phi_{k+j}^\epsilon)$
is isotopic to $\Phi_k^\epsilon$
for every even integer $j \ge 0$.  By Lemma
\xref{PropertyTwo} we have that
$h_j^\epsilon(\Phi_{k+j}^\epsilon)$ is isotopic to $\Phi_j^{-\epsilon}
\wedge_{\cal L} \Phi_k^\epsilon$.  In particular
$h_j^\epsilon(\Phi_{k+j}^\epsilon)$ is isotopic to a subsurface of
$\Phi_k^\epsilon$.  Now since $\Phi_{k+j}^\epsilon$ is isotopic to
$\Phi_k^\epsilon$ by ($1$), Lemma
\xref{SurfaceTwo}($2$) implies that $h_j^\epsilon(\Phi_{k+j}^\epsilon)$
is isotopic to $\Phi_k^\epsilon$.

A similar argument using ($1'$) shows that
($2'$) $h_j^{-\epsilon}(\Phi_{k+1+j}^{-\epsilon})$ is
isotopic to $\Phi_{k+1}^{-\epsilon}$ for all $j \ge 0$.

More generally we claim that ($3$) if $j$ is even and $l \ge \max(k,j)$ then
$h_j^\epsilon(\Phi_l^\epsilon)$ is isotopic to $\Phi_k^\epsilon$.  To
show this we first note that by ($1$),
$\Phi_l^\epsilon$ is isotopic in $\tilde F$ to $\Phi_{k+j}^\epsilon$
and so by Proposition \xref{Containment} $\Phi_l^\epsilon$ is isotopic to
$\Phi_{k+j}^\epsilon$ in $\Phi_j^\epsilon$.  Hence
$h_j^\epsilon(\Phi_l^\epsilon)$ is isotopic in $\tilde F$ to
$h_k^\epsilon(\Phi_{k+j}^\epsilon)$, which we have already shown is isotopic to
$\Phi_k^\epsilon$.

In the same way, using ($2'$), we see that ($3'$) if $j$ is even and $l
\ge \max(k+1,j)$, then $h_j^{-\epsilon}(\Phi_l^{-\epsilon})$ is
isotopic to $\Phi_{k+1}^{-\epsilon}$.

We claim that ($4$) for all $j \ge 0$ the surface $\Phi_k^\epsilon$ is
isotopic to a subsurface of $\Phi_j^{-\epsilon}$ and that ($4'$)
$\Phi_{k+1}^{-\epsilon}$ is isotopic to a subsurface of
$\Phi_j^{\epsilon}$.  Since the $\Phi_j^{\pm\epsilon}$ are
nested up to isotopy, we need only prove this for even $j$.  By
($2$) we have that $h_j^\epsilon(\Phi_{k+j}^\epsilon)$ is isotopic to
$\Phi_k^\epsilon$.  Lemma \xref{PropertyTwo} implies that
$h_j^\epsilon(\Phi_{k+j}^\epsilon)$ is isotopic to
$\Phi_j^{-\epsilon}\wedge_{\cal L}\Phi_k^\epsilon$.  The fact that
$\Phi_k^\epsilon$ is isotopic to $\Phi_j^{-\epsilon}\wedge_{\cal
L}\Phi_k^\epsilon$ implies that $\Phi_k^\epsilon$ is isotopic to a
subsurface of $\Phi_j^{-\epsilon}$, as asserted by ($4$).  A similar argument
using ($2'$) proves ($4'$).

We next claim that ($5$) if $j$ is odd then
$h_j^{\epsilon}(\Phi_k^{\epsilon})$ is isotopic to
$\Phi_{k+1}^{-\epsilon}$ and ($5'$) if $j$ is odd then
$h_j^{-\epsilon}(\Phi_{k+1}^{-\epsilon})$ is isotopic to
$\Phi_{k}^{\epsilon}$.  First by ($1$) we have that
$h_j^{\epsilon}(\Phi_k^{\epsilon})$ is isotopic to
$h_j^\epsilon(\Phi_{k+1+j}^\epsilon)$ which, by Proposition
\xref{PropertyTwo}, is isotopic to $\Phi_j^\epsilon \wedge_{\cal
L} \Phi_{k+1}^{-\epsilon}$.  Claim ($4'$) implies that this essential
intersection is isotopic to $\Phi_{k+1}^{-\epsilon}$.  The proof of
($5'$) is similar but uses ($1'$) and ($4$) in place of ($1$) and
($4'$).

Now fix an even integer $m > k$.  By Proposition $\OddK$ we have that
$h_m^\epsilon(\Phi_{2m+1}^\epsilon)$ is isotopic to a
$\tau_\epsilon$-invariant subsurface $B_\epsilon$ of
$\Phi_1^\epsilon$.
By ($3$) we have that $\Phi_k^\epsilon$ is isotopic to
$h_m^\epsilon(\Phi_{2m+1}^\epsilon)$ and hence to $B_\epsilon$.
Similarly, using ($3'$), we see that
$\Phi_{k+1}^{-\epsilon}$ is isotopic to a $\tau_{-\epsilon}$-invariant
subsurface $B_{-\epsilon}$ of $\Phi_1^{-\epsilon}$.  Moreover, since
$\tau_\epsilon = h_1^\epsilon$, we have that $\Phi_k^\epsilon$ is
isotopic to $\tau_\epsilon(\Phi_k^\epsilon) =
h_1^\epsilon(\Phi_{k}^\epsilon)$ which in turn is isotopic to
$\Phi_{k+1}^{-\epsilon}$ by ($5$).

Thus we have defined surfaces $B_\epsilon\subset\Phi_1^\epsilon$ and
$B_{-\epsilon}\subset\Phi_1^{-\epsilon}$ which are invariant under
$\tau_\epsilon$ and $\tau_{-\epsilon}$ respectively and are both
isotopic to $\Phi_k^\epsilon$.  By \xref{JSJStuff} it follows that
there exist $I$-pairs $(E_\epsilon, B_\epsilon)\subset
(M_{\tilde F}^\epsilon,\tilde F)$ whose associated $\partial I$-subbundles are
both isotopic to the subsurface $\Phi_k^\epsilon$ of $\tilde F$.
After modifying these $I$-pairs by isotopies we obtain a semi-fibered
submanifold $N = E_+ \cup E_-$ contained in $M$ whose semi-fiber $B$
is isotopic to the large subsurface $\Phi_k^\epsilon$ of $\tilde F$.
We may take $N$ to be contained in $\int M$.

We now are ready to show that if $\Phi_k^\epsilon$ is non-empty then
$\tilde F$ is a semi-fiber.  First we argue that each component of $\bdry N$ is a
$\pi_1$-injective torus in $M$.  Since $N$ has a large semi-fiber, it
is clear that each component of $\partial N$ is $\pi_1$-injective in
$N$.  It therefore suffices to show that $N$ is $\pi_1$-injective in
$M$. A homotopically non-trivial loop in $N$ which is contained in $B$
is homotopically non-trivial in $M$ because $B$ is $\pi_1$-injective
in $\tilde F$.  Now consider a loop $\alpha$ in $N$ which is not
homotopic to a loop in $B$.  After modifying $\alpha$ by a free
homotopy we may take it to be a composition of paths $\alpha_1,
\ldots, \alpha_{2n}$ for some $n>0$ such that each $\alpha_i$ is a
path in $E_{(-1)^i}$ which has its endpoints in $B$ and is not
fixed-endpoint homotopic in $N$ to a path in $B$.  The map of pairs
$\alpha_i\colon(I,\bdry I)\to(E_{(-1)^i},B_{(-1)^i})$ must be
homotopic to a homeomorphism onto a fiber.  Since the $I$-bundles
$E_+$ and $E_-$ are essential in $M_{\tilde F}$ the path $\alpha_i$ is
not fixed-endpoint homotopic in $M$ to a path in $\tilde F$.  It
follows that $N$ is $\pi_1$-injective in $M$.

Since $M$ is a simple knot manifold every component of $\partial N$ is
boundary parallel.  Furthermore, since $N$ has a large semi-fiber, it
cannot be homeomorphic to $S^1\times S^1 \times I$.  Thus $C =
\overline{M- N}$ is a collar on $\partial M$.  The surface $A
= \overline{\tilde F- B}$ is a $\pi_1$-injective subsurface of
$\tilde F$ without disk component which is properly embedded in $C$.
In particular $A$ is $\pi_1$-injective in $C$ and hence each component
of $A$ is an annulus.  Furthermore, by [\cite{Waldhausen}], either $A$ has a
boundary parallel component or $A$ is vertical in the sense that it is
mapped by some homeomorphism of $C$ onto $S^1\times S^1 \times I$
which maps $A$ to $X\times S^1 \times I$ where $X$ is a finite subset
of $S^1$.  Assume that $A$ has a boundary parallel component $A_0$.
Then $A_0\subset
\tilde F$ is isotopic relative to $\partial A_0$ to an annulus $A_0'$
in $\partial N$.  But the torus $\partial N$ is a union of essential
annuli which are components of the frontiers of the $I$-pairs
$(E_\epsilon, B_\epsilon)$.  Thus the inclusion map of the annulus
$A_0'$ can be regarded as a reduced homotopy between the inclusion
maps of the two boundary components of $A_0$.  In particular, $A_0'$
is not homotopic into $\tilde F$.  This is a contradiction, hence $A$
is vertical.  It follows easily that the semi-fibration of $N$ can be
extended over $C$ to obtain a semi-fibration of $M$ with semi-fiber
$\tilde F$.
\EndProof

\subsection{Bounding the length of a reduced homotopy}

\Theorem
\tag{LengthBound}
Let $F$ be a connected essential surface in a simple knot manifold $M$.
Suppose that $F$ is not a semi-fiber.  Let $g$ and $m$ denote
respectively the genus and number of boundary components of $F$.  Then
any reduced homotopy in the pair $(M,F)$ having a large time-$0$ map
has length at most $8g +3m - 8$.
\EndTheorem

For the proof of this theorem we will need the following construction.

\paragraph
\tag{AssociatedSplittingSurface}
Suppose that $F$ is any connected, nonseparating, essential surface in
$M$, and let $\tilde F$ denote the boundary of a regular neighborhood
$N$ of $F$. It is clear that $\tilde F$ has a unique transverse
orientation such that $\tilde F$ is a splitting surface and
$M_{(\tilde F)}^+=N$.  We shall call $\tilde F$, equipped with this
transverse orientation, a {\it splitting surface associated to $F$}.

If $F$ is any connected, separating, essential surface in $M$, we
refer to $F$ itself, equipped with either transverse orientation, as a
{\it splitting surface associated to $F$.}

\paragraph
\demo{Proof of Theorem \xref{LengthBound}}
Let $\tilde F$ be a splitting surface associated to $F$.  Since $F$ is
not a semi-fiber, it is clear that $\tilde F$ is not a semi-fiber.  We
define the subsurfaces $\Phi_k^\epsilon$ for $k\ge 0$ and
$\epsilon\in\{\pm1\}$ as in Section \xref{EssHom}.
Let $H$ be a reduced homotopy of length $l$ in the pair $(M,F)$
starting on the $\epsilon$ side.
Let us set $\tilde l=l$ if $F$ separates $M$, and $\tilde l=2l$ if $F$
does not separate $M$. Then the homotopy $H$ determines a reduced homotopy
$\tilde H$ of length $\tilde l$ in the pair $(M,\tilde F)$ such that
$\tilde H_0$ is large. This is obvious if $F$ separates; if $F$ does
not separate, it follows from the fact that the two components of $\tilde
F$ cobound a product.   Let $\tilde g$ and $\tilde
m$ denote respectively the total genus and number of boundary components of
$\tilde F$.  We have $\tilde g = g$ and $\tilde m = m$ in the
separating case, and $\tilde
g = 2g$ and $\tilde m = 2m$ in the nonseparating case. Hence it
suffices to show that  the length $\tilde l$ of $\tilde H$ is at
most $8\tilde g+3\tilde m - 8$ in the separating case, and is at
most $8\tilde g+3\tilde m - 16$ in the nonseparating case.

Set $n = \left[{\tilde l+1\over 2}\right]$, and consider the subsurfaces
$$\Phi_1^\epsilon\supset \Phi_3^\epsilon\supset
\ldots \supset \Phi_{2n-1}^\epsilon\supset \Phi_{2n+1}^\epsilon $$ of
$\tilde F$. Since
the reduced homotopy $\tilde H$ has length $\tilde l\ge2n-1$, it
follows from Proposition \xref{Surfaces} that
$\Phi_{2n-1}^\epsilon\ne\emptyset$.
Since $\tilde F$ is not a semi-fiber, it now follows from Proposition
\xref{StrictContainment} that
$\Phi_{2i-1}^\epsilon$ is not a regular neighborhood
of $\Phi_{2i+1}^\epsilon$ for $i=1, \ldots, n$.  Furthermore, by
Proposition \xref{EvenChar}, each of these surfaces has even Euler
characteristic.  Using these facts, we will show that $n \le 4\tilde g
+ 3\tilde m/2 - 4$ in the separating case, and that $n \le 4\tilde g +
3\tilde m/2 - 8$ in the nonseparating case; this implies the desired
conclusion, since $\tilde l\le 2n$.

If $A$ is any large subsurface of $\tilde F$ then we
will set $c(A) = \genus(A) - 3\chi(A)/2 - |A|$.  Note that $c(A)$ is always
non-negative, and is an integer if $A$ has even Euler characteristic.
We have $c(\tilde F) = 4\tilde g + 3\tilde m/2 - 4$ in the separating
case, and $c(\tilde F) = 4\tilde g + 3\tilde m/2 -8$ in the
nonseparating case. Hence it suffices to show that $n\le c(\tilde F)$.
Thus the proof reduces to the following general claim: if $A$ and $B$
are large subsurfaces of $\tilde F$, each of which
has even Euler characteristic, and if $B$ is contained in the interior
of $A$, then $c(B) < c(A)$ unless $A$ is a regular neighborhood of
$B$.

To prove the claim, it suffices to show that if $A_0$ is a component
of $A$, and if we set $B_0 = B\cap A_0$, then $c(B_0) < c(A_0)$ unless
$A_0$ is a regular neighborhood of $B_0$.  Note that we have
$\genus(B_0) \le \genus(A_0)$ and, since $A$ and $B$ are large,
$\chi(B_0) \ge \chi(A_0)$.  Thus we need only consider the two cases
where $B_0=\emptyset$ and where $B_0$ is connected.  The case
$B_0=\emptyset$ is easy because $c(A)>0$ for any non-empty large
subsurface $A$.  For the case where $B_0$ is connected we may assume
that $\genus(B_0) = \genus(A_0)$ and $\chi(B_0) =
\chi(A_0)$, and we must show that $A_0$ is a regular neighborhood of
$B_0$.  Since $A$ and $B$ are large, no
component of $\overline{A_0 - B_0}$ is a disk.  Thus the
condition $\chi(B_0) = \chi(A_0)$ implies that each component of
$\overline{A_0 - B_0}$ is an annulus.  None of these annuli
can separate $A_0$ since $B_0$ is connected.  On the other hand, since
$\genus(B_0) = \genus(A_0)$ there cannot exist a simple closed curve
in $A_0$ which has non-zero intersection number with a core curve of
an annulus component of $\overline{A_0 - B_0}$.  It follows
that each component of $\overline{A_0 - B_0}$ is a collar on a
boundary component of $A_0$ and hence that $A_0$ is a regular
neighborhood of $B_0$.  This completes the proof of the claim, and of
the theorem.
\enddemo
\EndProof

\deepsection{Boundary slopes of essential surfaces and singular surfaces}
\tag{DottedPhis}

Most of the work in this section is devoted to proving Theorem
\xref{DottedLengthBound}, which is a refinement of Theorem \xref{LengthBound} 
and gives a bound on the length of a reduced homotopy whose time-$0$
map is an essential path. Combining Theorem \xref{DottedLengthBound}
with the results in Section \xref{ReducedHomotopies} we obtain a proof
of Theorem \xref{HighGenusTheorem} and its corollaries
\xref{HighGenusDisk}, \xref{HighGenusSeifert} and
\xref{CameronCorollary}. These results give bounds on
$\Delta(\alpha,\beta)$ where $\beta$ is a boundary slope and $\alpha$
is either another boundary slope, a very small filling slope or a
Seifert-fibered filling slope.

We will need the following definition in this section.

\Definition
Let $\tilde F$ be a splitting surface in a simple knot manifold
$M$. We say that $\tilde F$ {\it admits a long rectangle} if
there exists a reduced homotopy $H:I\times I\to M$ in the pair
$(M,\tilde F)$ having length at least $|\partial \tilde F|$ in the
pair $(M,\tilde F)$ such that $H_t(\partial I) \subset \partial M$ for
all $t\in I$ and $H_0$ is an essential path.
\EndDefinition

\subsection{Reduced homotopies and outer subsurfaces}
\tag{DottedEssHom}
\paragraph
\tag{Hypotheses}
{\it Hypotheses.\/}
Throughout Subsection \xref{DottedEssHom} we will assume that $M$ is a
simple knot manifold and that $\tilde F$ is a splitting surface in $M$
which admits a long rectangle.

We will define $\Phi_{\pm1}$ and $\tau_{\pm1}$ as in Subsection
\xref{EssHom}.  For every $k\ge0$ we fix subsurfaces $\Phi_k^{\pm1}$ of
$\tilde F$ for which the conclusions of Proposition \xref{Surfaces}
hold and homeomorphisms $h_k^{\pm1}$ for which the conclusions of
Proposition \xref{TimeOneMaps} hold.

\Lemma \tag{ArcLemma}
For each $\epsilon \in \{\pm1\}$, every component of $\partial \tilde
F$ is isotopic in $\tilde F$ to a unique boundary component of
$\Phi_1^\epsilon$.  Furthermore, if $c$ and $c'$ are components of
$\partial \tilde F$ which cobound an annulus component of $M_{\tilde
F}^\epsilon\cap\partial M$, then the boundary components of
$\Phi_1^\epsilon$ which are isotopic to $c$ and $c'$ are interchanged
by $\tau_\epsilon$.
\EndLemma

\Proof
Set $m=|\partial\tilde F|$.  The existence of a long rectangle means
that there is a reduced homotopy $H$ of length $m$, such that
$H_t(\partial I) \subset \partial M$ for all $t\in I$ and $f=H_0$ is
an essential path.  As in \xref{ArcToGlasses} we extend $f$ to an
admissible pair of glasses $\hat f:\Gamma\to\tilde F$, and extend $H$
to a length-$m$ reduced homotopy $\hat H$ with $\hat H_0 = \hat f$.
By construction the homotopy $\hat H$ has the property that for each
of the ``rims'' $\bar l_i$, we have $\hat H_t(\bar l_i)\subset
\partial M$ for all $t\in [0,1]$.

Write $\hat H$ as a composition of essential basic homotopies $\hat
H^1, \dots, \hat H^m$.  Because the $\hat H^i$ are essential
homotopies, $\hat H^i$ maps $\bar l_1\times I$ to an annulus component
$A_i$ of $M_{\tilde F}\cap\partial M$ and maps the two components of
$\bar l_1\times \partial I$ to distinct components of $\partial A_i$.
Since $\hat H$ is reduced the annuli $A_i$ and
$A_{i+1}$ are on opposite sides of their common boundary component
$\hat H_1^i(\bar l_1) = \hat H_0^{i+1}(\bar l_1)$ for each $i=1,
\ldots, m-1$.  Since $m$ is the number of components of $\tilde F\cap
\partial M$ it follows that the annuli $A_i$ are distinct and that
every annulus component of $M_{\tilde F}\cap\partial M$ appears as one
of the $A_i$. If we set $A_0=A_m$, then for each component $c$ of
$\partial \tilde F$ there is some $i$ with $0 \le i < m$ such that
$c$ is the common boundary curve of the two annuli $A_i$ and
$A_{i+1}$, one of which is contained in $M_{\tilde F}^+$ and the other
in $M_{\tilde F}^-$.

To prove the first assertion of the lemma let $c$ be a component of
$\partial \tilde F$ and let $\epsilon\in\{\pm1\}$ be given.  The curve
$c$ is a boundary component of some annulus $A_j\subset M_{\tilde
F}^\epsilon\cap\partial M$.  Lemma \xref{PhiOne} implies
that $\hat H^j_0$ is homotopic to a map from $\Gamma$ to
$\Phi_1^\epsilon$.  In particular $c$ is homotopic in $\tilde F$ to a
(singular) curve in $\Phi_1^\epsilon$.  Since $c$ is a boundary
component of $\tilde F$ it follows that $c$ is homotopic to a boundary
component of $\Phi_1^\epsilon$.

To prove the second assertion, suppose that $c$ and $c'$ are boundary
curves of $\tilde F$ which cobound an annulus component $A$ of
$M_{\tilde F}^\epsilon\cap\partial M$.  Then we have $A = A_j$ for
some $j$.  Let $\gamma$ be the boundary component of $\Phi_1^\epsilon$
which is homotopic to $c$.  Lemma \xref{PhiOne} implies that $\hat
H^j_0$ is homotopic to a map $g:\Gamma\to\Phi_1^\epsilon$ such that
$g(\bar l_1) = \gamma$.  Applying Proposition \xref{HSigma} to the
homotopy $H^j$, with $g$ in place of $f$ and $\Gamma$ in place of $K$,
we conclude that $c' = \hat H^j_1(\bar l_1)$ is homotopic in $\tilde
F$ to $\tau_\epsilon(\gamma)$.  The assertion follows immediately.
\EndProof

\Lemma 
\tag{BoundaryToBoundary}
If $\epsilon \in \{\pm1\}$ and if $\gamma$ is a simple closed curve in
$\Phi_k^\epsilon$ which is isotopic in $\tilde F$ to a component of
$\partial \tilde F$, then $h_k^\epsilon(\gamma)$ is also isotopic to
some component of $\partial \tilde F$.
\EndLemma

\Proof
For $k=0$ the assertion is trivial. To prove the lemma for $k=1$,
we may assume that $\gamma$ is a component of $\partial
\Phi_k^\epsilon$.  Let $c$ denote the component of $\partial \tilde F$
which is isotopic in $\tilde F$ to $\gamma$, and let $A$ denote the
annulus component of $M_{\tilde F}^\epsilon\cap\partial M$ having $c$
as a boundary curve.  If $c'$ denotes the other boundary curve of $A$,
then it follows from Lemma \xref{ArcLemma} that $c'$ is isotopic in
$\tilde F$ to a boundary component $\gamma'$ of $\Phi_1^\epsilon$ and
that $h_1^\epsilon(\gamma) = \tau_\epsilon(\gamma) = \gamma'$.  This
proves the assertion in this case.

Now assume that $k>1$ and that the assertion holds with $k$ replaced
by $k-1$, both for $\epsilon=1$ and $\epsilon=-1$.  Suppose
we are given $\epsilon\in\{\pm1\}$ and a simple closed curve $\gamma\subset
\Phi_k^\epsilon\subset\Phi_{k-1}^\epsilon$ which is isotopic in
$\tilde F$ to a component of $\partial \tilde F$.  We will apply
Proposition \xref {PropertyThree} with $i = k-1$ and $j=1$.  According
to \xref {PropertyThree}, the map $h_{k-1}^\epsilon|\Phi_k^\epsilon$
is homotopic to a map $g_{k-1}^\epsilon:\Phi_{k-1}^\epsilon
\to \Phi_1^{-\epsilon}$ such that $h_k^\epsilon$ is homotopic to
$h_1^{(-1)^{k-1}\epsilon}\circ g_{k-1}^\epsilon$.  Now
$g_{k-1}^\epsilon(\gamma)$ is homotopic to $h_{k-1}^\epsilon(\gamma)$
which, by the induction hypothesis, is homotopic to some component of
$\partial \tilde F$.  By the case $k=1$ of the proposition we know
that $h_1^{(-1)^{k-1}\epsilon}(h_{k-1}^\epsilon(\gamma))$ is homotopic
to a component of $\partial \tilde F$.  Since $h_k^\epsilon$ is
homotopic to $h_1^{(-1)^{k-1}\epsilon}\circ g_{k-1}^\epsilon$ this
shows that $h_k^\epsilon(\gamma)$ is homotopic to a component of
$\partial \tilde F$.
\EndProof
 
Recall from \xref{OuterParts} that if $A$ is a subsurface of a compact
orientable surface $S$ then the outer part of $A$ is denoted $\dot A$. 

\Lemma
\tag{OuterToOuter}
For either $\epsilon\in\{\pm1\}$ and for any large subsurface $A$ of
$\Phi_k^\epsilon$, the outer part of $h_k^\epsilon(A)$ is
$h_k^\epsilon(\dot A)$.
\EndLemma

\Proof
Set $B=h_k^\epsilon(A)$. It follows from Lemma \xref{BoundaryToBoundary} that
$h_k^\epsilon(\dot A)\subset\dot B$.  Let
$g:\Phi_k^{(-1)^{k+1}\epsilon}\to\Phi_k^\epsilon$ denote the inverse
of $h_k^\epsilon$, so that $g(B)=A$.  According to Proposition
\xref{Inverses} the map $g$ is isotopic to $h_k^{(-1)^{k+1}\epsilon}$
as a map from $\Phi_k^{(-1)^{k+1}\epsilon}$ to $\Phi_k^\epsilon$.  It
therefore follows from Lemma \xref{ArcLemma} that $g(\dot
B)\subset\dot A$, i.e. that $\dot B\subset h_k^\epsilon(\dot
A)$. \EndProof

\paragraph
\tag{Dottedh}
We now consider the outer parts of the surfaces $\Phi_k^\epsilon$ which,
according to our conventions, are denoted $\dot\Phi_k^\epsilon$.
Note that since
$$\tilde F=\Phi_0^\epsilon \supset \Phi_1^\epsilon \supset
\Phi_2^\epsilon\supset
\cdots $$
for $\epsilon\in\{\pm1\}$, it follows from \xref{OuterParts} that
$$\tilde F=\dot\Phi_0^\epsilon \supset \dot\Phi_1^\epsilon \supset
\dot\Phi_2^\epsilon\supset
\cdots .$$
It follows from Lemma \xref{OuterToOuter} that $h_k^\epsilon$
restricts to a homeomorphism from $\dot\Phi_k^\epsilon$ to
$\dot\Phi_k^{(-1)^{k+1}\epsilon}$.  This homeomorphism will be denoted
by $\dot h_k^\epsilon$.  In particular the involution $\tau_\epsilon =
h_1^\epsilon$ of $\Phi_1^\epsilon$ restricts to an involution
$\dot\tau_\epsilon = \dot h_1^\epsilon$ of $\dot \Phi_1^\epsilon$.  It
also follows from Lemma \xref{OuterToOuter} that if $A$ is any large
subsurface of $\Phi_k^\epsilon$ then $\dot h_k^\epsilon(\dot A)$ is
the outer part of $h_k^\epsilon(A)$.  Note also that, according to
Lemma \xref{TimeOneMaps}, if $\tilde F$ is given a consistent
orientation $\dot h_k^\epsilon:\dot\Phi_k^\epsilon\to \tilde F$ is
orientation reversing if $k$ is odd and orientation preserving if $k$
is even.

The following six results,
\xref{DottedPropertyThree} -- \xref{DottedStrictContainment}, are
analogues of \xref{PropertyThree}--\xref{StrictContainment}.
(Of course these results, unlike their counterparts in Section 4,
depend on the hypothesis stated at the beginning of the section that
$\tilde F$ admits a long rectangle.)

\Proposition\tag{DottedPropertyThree}
Let $i$ and $j$ be non-negative integers, and set $k=i+j$. Then for
each $\epsilon\in\{\pm1\}$, the map $\dot h_i^\epsilon|\dot\Phi_k^\epsilon$ is
homotopic in $\tilde F$  to an embedding $\dot g_i^\epsilon\colon \dot\Phi_k^\epsilon\to
\dot\Phi_{j}^{(-1)^{i}\epsilon}$ such that
$\dot h_{j}^{(-1)^{i}\epsilon}\circ \dot g_i^\epsilon$ is homotopic in $\tilde
F$ to $\dot h_k^\epsilon$.
\EndProposition

\Proof
By Proposition \xref{PropertyThree} we have that 
$h_i^\epsilon|\Phi_k^\epsilon$ is
homotopic in $\tilde F$  to an embedding $g_i^\epsilon\colon \Phi_k^\epsilon\to
\Phi_{j}^{(-1)^{i}\epsilon}$ such that
$h_{j}^{(-1)^{i}\epsilon}\circ g_i^\epsilon$ is homotopic in $\tilde
F$ to $h_k^\epsilon$.  Set $\dot g_i^\epsilon =
g_i^\epsilon|\dot\Phi_k^\epsilon$, and set $A=\dot
g_i^\epsilon(\dot\Phi_k^\epsilon)$.  To complete the proof it suffices
to show that $A\subset\dot\Phi_{j}^{(-1)^{i}\epsilon}$.

Since $h_{j}^{(-1)^{i}\epsilon}\circ g_i^\epsilon$ is homotopic to
$h_k^\epsilon$, the subsurface
$h_{j}^{(-1)^{i}\epsilon}(A)$
is homotopic into $h_k^\epsilon(\dot\Phi_k^\epsilon)$ which, by
\xref{Dottedh}, is equal to $\dot\Phi_k^{(-1)^{k+1}\epsilon}$.
In particular $h_{j}^{(-1)^{i}\epsilon}(A)$ is homotopic into
$\dot\Phi_{j}^{(-1)^{k+1}\epsilon}$.  Since
$h_{j}^{(-1)^{i}\epsilon}(A)$ is a large subsurface of
$h_j^{(-1)^i\epsilon}(\Phi_j^{(-1)^i\epsilon})=\Phi_j^{(-1)^{k+1}\epsilon}$,
and since 
$\dot\Phi_{j}^{(-1)^{k+1}\epsilon}$ is a union of components of
$\Phi_{j}^{(-1)^{k+1}\epsilon}$, it follows that
$h_{j}^{(-1)^{i}\epsilon}(A)$ is contained in
$\dot\Phi_{j}^{(-1)^{k+1}\epsilon}$, which by \xref{Dottedh} is
equal to $h_{j}^{(-1)^{i}\epsilon}(\dot\Phi_j^{(-1)^i\epsilon})$.
We therefore have $A\subset\dot\Phi_{j}^{(-1)^{i}\epsilon}$, as required.
\EndProof

\Proposition\tag{DottedPropertyTwo}
Let $i$ and $j$ be non-negative integers, and set $k=i+j$.  Then for
each $\epsilon \in \{\pm1\}$ the subsurface $\dot h_i^\epsilon(\dot
\Phi_k^\epsilon)$ is isotopic in $\tilde F$ to $\dot
\Phi_{i}^{(-1)^{i+1}\epsilon}\dot\wedge_{\cal L}
\dot\Phi_{j}^{(-1)^i\epsilon}$.
\EndProposition

\Proof
We may assume by Proposition \xref{PropertyTwo} that
$h_i^\epsilon(\Phi^\epsilon_k) =
\Phi_{i}^{(-1)^{i+1}\epsilon}\lint\Phi_{j}^{(-1)^i\epsilon}$.  It
follows from \xref{Dottedh} that
$\dot h_i^\epsilon(\dot\Phi^\epsilon_k)$
is the outer part of
$h_i^\epsilon(\Phi^\epsilon_k)$
and is therefore equal to the surface 
$\Phi_{i}^{(-1)^{i+1}\epsilon}\dotlint\Phi_{j}^{(-1)^i\epsilon}$.
By Lemma \xref{DottyLemma} this surface is isotopic in $\tilde F$ to
$\dot\Phi_{i}^{(-1)^{i+1}\epsilon}\dotlint\dot\Phi_{j}^{(-1)^i\epsilon}$.
\EndProof

\Lemma
\tag{DottedInvariant}
For any non-negative integer $k$ and for each $\epsilon\in\{\pm1\}$
the subsurface $\dot h_k^\epsilon(\dot\Phi_{2k+1}^\epsilon)$ is isotopic in
$\tilde F$ to a subsurface $A$ of $\dot\Phi_1^{(-1)^k\epsilon}$ with the
property that $\dot\tau_{(-1)^k\epsilon}(A)$ is a subsurface of
$\dot\Phi_1^{(-1)^k\epsilon}$ which is
isotopic to $A$ in $\dot\Phi_1^{(-1)^k\epsilon}$.
\EndLemma

\Proof
By Proposition \xref{InvariantUpToIsotopy} we know that
$h_k^\epsilon(\Phi_{2k+1}^\epsilon)$ is isotopic in $\tilde F$ to a
subsurface $A_0$ of $\Phi_1^{(-1)^k\epsilon}$ with the property that
$\tau_{(-1)^k\epsilon}(A_0)$ is isotopic to $A_0$ in
$\Phi_1^{(-1)^k\epsilon}$.  It follows from \xref{Dottedh} that $\dot
h_k^\epsilon(\dot\Phi_{2k+1}^\epsilon)$ is isotopic in $\tilde F$ to
$\dot A_0$.  It also follows from \xref{Dottedh} that
$\dot\tau_{(-1)^k\epsilon}(\dot A_0) = h_1^{(-1)^k\epsilon}(\dot A_0)=
\dot h_1^{(-1)^k\epsilon}(\dot A_0)$ is equal to the outer part of
$h_1^{(-1)^k\epsilon}(A_0)=\dot\tau_{(-1)^k\epsilon}(A_0)$.  Since
$\dot\tau_{(-1)^k\epsilon}(A_0)$ is isotopic to $A_0$ it follows that
$\dot\tau_{(-1)^k\epsilon}(\dot A_0)$ is isotopic in $\tilde F$ to $\dot
A_0$.  Since the subsurfaces $\dot\tau_{(-1)^k\epsilon}(\dot A_0)$ and
$\dot A_0$ of $\dot\Phi_1^\epsilon$ are isotopic in $\tilde F$,
it follows from Lemma \xref{SurfaceTwo} that they are isotopic in
$\dot\Phi_1^{(-1)^k\epsilon}$.  Thus we may take $A = \dot A_0$.
\EndProof 

\Proposition
\tag{DottedOddK}
For any non-negative integer $k$ and for each $\epsilon\in\{\pm1\}$
the surface
$\dot h^\epsilon_k(\dot\Phi_{2k+1}^\epsilon)$ is isotopic in $\tilde F$ to
a subsurface of $\dot\Phi_1^{(-1)^k\epsilon}$ which is invariant under
the free involution $\dot\tau_{(-1)^k\epsilon}$.  In particular,
$\dot\Phi_{2k+1}^\epsilon$ admits a free involution.
\EndProposition

\Proof
This is an immediate consequence of \xref{DottedInvariant} and
\xref{Invariant}.
\EndProof

\Corollary
\tag{DottedEvenChar}
For each $\epsilon\in\{\pm1\}$ and each odd integer $k>0$, the Euler
characteristic $\chi(\dot\Phi_{k}^\epsilon)$ is even.
\EndCorollary
\NoProof

\Proposition
\tag{DottedStrictContainment}
Let $k$ be a non-negative integer and let $\epsilon\in\{\pm1\}$ be given.
If $\dot\Phi_k^\epsilon$ and
$\dot\Phi_{k+2}^\epsilon$ are isotopic in $\tilde F$ then either
$\dot\Phi_k^\epsilon = \emptyset$ or $\tilde F$ is a semi-fiber.
\EndProposition
\Proof
This is formally identical with the proof of Proposition
\xref{StrictContainment}.  All occurrences of $\Phi$, $h$, $\tau$ and
$\lint$ are replaced by $\dot\Phi$, $\dot h$, $\dot\tau$ and
$\dotlint$ respectively.  References to Proposition \xref{PropertyTwo}
and Lemma \xref{OddK} are replaced by references to Proposition
\xref{DottedPropertyTwo} and Lemma \xref{DottedOddK}
\EndProof

\subsection{The distance bounds}
The next result is a strengthened version of Theorem \xref{LengthBound}
that applies to a reduced homotopy whose time-$0$ map is an essential
path.

\Theorem
\tag{DottedLengthBound}
Let $F$ be a connected essential surface in a simple knot manifold $M$.
Suppose that $F$ is not a semi-fiber.  Set $g = \genus(F)$ and $m =
|\bdry F|$.  Let $H$ be any reduced homotopy in the pair $(M,F)$ such
that $H_0$ is an essential path in $F$ and $H_t(\partial I)\subset
\partial M$ for each $t\in I$. Then the length of $H$ is at
most $4g +3m - 4$.
\EndTheorem

\Proof
The proof will be similar to that of Theorem \xref{LengthBound}.  Let $\tilde
F$ be a splitting surface associated to $F$.  Since $F$ is not a
semi-fiber it follows that $\tilde F$ is not a semi-fiber.

Let $l$ denote the length of $H$.  We may assume $l>m$.
Set $\tilde l = l$ if $F$ is separating and $\tilde l = 2l$ if $F$ is
nonseparating.  As in the proof of Theorem \xref{LengthBound}, the homotopy
$H$ determines a homotopy $\tilde H$ of length $\tilde l$ in the pair
$(M,\tilde F)$ such that $\tilde H_0$ is an essential path.
Let $\tilde g$ and $\tilde m$ denote respectively the total genus and number
of boundary components of $\tilde F$.  We have
$\tilde m = m$ and $\tilde g = g$ in the separating case and
$\tilde m = 2m$, $\tilde g = 2g$ in the nonseparating case.
It therefore suffices to prove that $\tilde l
\le 4\tilde g + 3\tilde m - 4$ if $F$ is separating and $\tilde l
\le 4\tilde g + 3\tilde m - 8$ if $F$ is nonseparating.

Since $l>m$ we also have $\tilde l > \tilde m$.  In particular this
means that $\tilde F$ admits a long rectangle and hence
that the hypotheses stated in \xref{Hypotheses} hold.  We define the
subsurfaces $\Phi_k^{\pm1}$ and $\dot\Phi_k^{\pm1}$ for $k\ge 0$ as in
\xref{Hypotheses} and \xref{Dottedh}.

Let us say that a subsurface $A$ of $\tilde F$ is {\it allowable} if
$A$ is a large subsurface with even Euler characteristic, and if $A$
is an outer subsurface (see \xref{OuterParts}).  If $A$ is any large
subsurface of $\tilde F$, let $\nu(A)$ denote the number of components
of $\partial \tilde F$ which are homotopic into $A$.  Note that
$\nu(A)>0$ for any non-empty allowable subsurface $A$.  We set
$$c(A) = \genus(A) - {\chi(A)\over 2} - |A| + \nu(A) =
2\genus(A) + {|\partial A|\over 2} - 2|A| + \nu(A).$$
Note that $c(A)$ is non-negative and integer-valued if $A$ is
allowable.  Moreover, $c(A)>0$ if $A$ is nonempty and allowable.

Set $n = \left[{\tilde l + 1\over 2}\right]$ and 
define $\epsilon\in\{\pm1\}$ by the condition that the homotopy
$\tilde H$ starts on the $\epsilon$-side. 
Consider the subsurfaces
$$\dot\Phi_1^\epsilon\supset \dot\Phi_3^\epsilon\supset
\ldots \supset \dot\Phi^\epsilon_{2n-1}\supset \dot\Phi^\epsilon_{2n+1} .$$
Since the hypotheses stated in \xref{Hypotheses} hold, Corollary
\xref{DottedEvenChar} implies that each of these surfaces has even
Euler characteristic, and, in view of the definition of the
$\dot\Phi_k^\epsilon$, it follows that each of these surfaces is
allowable.  On the other hand, since the reduced homotopy $\tilde H$
has length $\tilde l \ge 2n-1$, it follows from \xref{ArcToGlasses}
and Proposition \xref{Surfaces} that there is an admissible pair of
glasses $\gamma:\Gamma\to\tilde F$ which is homotopic in $\tilde F$ to
a map from $\Gamma$ to $\Phi_{2n-1}^\epsilon$.  In particular there is
a map $\alpha:S^1 \to \partial\tilde F$ which is homotopic in $\tilde
F$ to a map from $S^1$ to a component $A$ of $\Phi_{2n-1}^\epsilon$.
It follows that $A$ must be an outer component of
$\Phi_{2n-1}^\epsilon$, and hence that $\dot\Phi_{2n-1}^\epsilon
\not=\emptyset$.  Since $\tilde F$ is not a semi-fiber and since
the hypotheses stated in \xref{Hypotheses} hold, we may apply
Proposition \xref{DottedStrictContainment} to conclude that
$\dot\Phi_{2i-1}^\epsilon$ is not a regular neighborhood of
$\dot\Phi_{2i+1}^\epsilon$ for $i=1, \ldots, n$. We will show that in
this situation $n \le 2\tilde g + 3\tilde m/2 - 2$ if $F$ is
separating and $n \le 2\tilde g + 3\tilde m/2 - 4$ if $F$ is
non-separating.  Since $\tilde l \le 2n$ this will imply the theorem.

We have $c(\tilde F) = 2\tilde g + 3\tilde m/2 - 2$ if $F$ is
separating and $c(\tilde F) = 2\tilde g + 3\tilde m/2 - 4$ if $F$ is
nonseparating.  Hence it will suffice to show that if $A$ and $B$ are
two allowable subsurfaces with $B\subset \int A$, then $c(B) < c(A)$
unless $A$ is a regular neighborhood of $B$.

As in the proof of Theorem \xref{LengthBound} it suffices to show, for a
component $A_0$ of $A$ and the subsurface $B_0 = B\cap A_0$,
that $c(B_0) < c(A_0)$ if
$A_0$ is not a regular neighborhood of $B_0$.  We have
that $\genus(B_0) \le \genus(A_0)$, that $\chi(B_0) \ge \chi(A_0)$, and
$\nu(B_0)
\le \nu(A_0)$.  Thus we need only consider the two cases where $|B_0|
\le |A_0| = 1$, i.e. where $B_0$ is empty and where $B_0$ is
connected.  The case $B_0=\emptyset$ is easy since $A_0
\not=\emptyset$ implies $c(A_0) > 0$.  For the case where $B_0$ is
connected we observe that if $c(B_0) = c(A_0)$ then we have
$\genus(B_0) = \genus(A_0)$ and $\chi(B_0) = \chi(A_0)$; it then
follows as in the proof of Theorem \xref{LengthBound} that $A_0$ is a regular
neighborhood of $B_0$.
\EndProof

We are now ready to state and prove one of the main results of this
paper.  Recall that the function $N(s,n,v)$ was defined in \xref{Nsnvf}.

\Theorem
\tag{HighGenusTheorem}
Let $M$ be a simple knot manifold and $F\subset M$ an essential
bounded surface with boundary slope $\beta$ which is not a semi-fiber.
Let $(S,X,h)$ be a singular surface which is well-positioned with
respect to $F$ and has boundary slope $\alpha$. Set $s=\genus S$,
$n=|\partial S-X|$, $v=|X|$, $g =
\genus F$ and $m=|\partial F|$. Then
$$\Delta(\alpha,\beta)\le\left({4g-3\over m} + 3\right)N(s,n,v).$$
\EndTheorem

\Proof
According to Proposition \xref{GraphProposition}, 
there exists an essential homotopy $H:I\times I\to M$ having
length $$l \ge {m\D(\a,\b)\over N(s,n,v)}-1. \leqno{(1)}$$ such that
$H_0$ is an essential path in $F$ and
$H_t(\partial I)\subset\partial M$ for all $t\in I$.
By Proposition \xref{DottedLengthBound} we have that
$$ l \le 4g + 3m - 4. \leqno{(2)}$$
The conclusion follows from the inequalities (1) and (2).
\EndProof 

\Corollary
\tag{HighGenusDisk}
Let $M$ be a simple knot manifold and $F\subset M$ an essential
bounded surface with boundary slope $\beta$ which is not a
semi-fiber. Set $g=\genus F$ and $m=|\partial F|$.  Let $\alpha$ be a
slope in $\partial M$. If $M(\alpha)$ is very small, or more generally
if $F\subset M\subset M(\alpha)$ is not $\pi_1$-injective in
$M(\alpha)$, then
$$\D(\a,\b)\le{20g-15\over m}+15.$$
\EndCorollary

\Proof
We invoke Corollary \xref{Disk} to obtain a singular surface
$(S,X,h)$, well-positioned with respect to $F$, such that $\genus S=0$
and $|X|=1$. The conclusion now follows from Theorem
\xref{HighGenusTheorem} because for any $v\ge1$ we have
$N(0,1,v)\le 5$.
\EndProof

\Corollary
\tag{HighGenusSeifert}
Let $M$ be a simple knot manifold and $F\subset M$ an essential
bounded surface with boundary slope $\beta$ which is not a
semi-fiber. Let $\alpha$ be a slope in $\partial M$. Set $g=\genus F$
and $m=|\partial F|$.  If $M(\alpha)$ is a Seifert fibered space or if
there exists a $\pi_1$-injective map from $S^1\times S^1$ to $M$ then
$$\D(\alpha,\beta)\le{24g-18\over m}+18.$$
\EndCorollary

\Proof
We invoke Corollary \xref{Seifert} to obtain a singular surface
$(S,X,h)$, well-positioned with respect to $F$, such that either
$\genus S=0$ and $|X|=1$, or $\genus S=1$ and $|X|=0$. The conclusion
now follows from Theorem \xref{HighGenusTheorem} because for any $v\ge1$ we
have $N(0,1,v)\le 5$ and $N(1,0,v)=6$.
\EndProof

\Corollary
\tag{CameronCorollary}
Let $M$ be a simple knot manifold.  Suppose that, for $i = 1,2$,
that $F_i\subset M$ is an essential bounded surface of genus
$g_i$ with boundary slope $\beta_i$.  Let $m_i = |\partial F_i|$.
If $F_1$ is not a semi-fiber then we have
$$\Delta(\beta_1,\beta_2)\le
\left({4g_1 - 3\over m_1}+3\right)\left(\left[{12g_2 - 12\over
m_2}\right]+ 6\right).$$
\EndCorollary

\Proof
We apply Proposition \xref{Cameron} with $F=F_1$ and $S=F_2$.  This
gives a singular surface $(F_2,\partial F_2, h)$ which is
well-positioned with respect to $F_1$.  Theorem
\xref{HighGenusTheorem} then implies that
$$\Delta(\beta_1,\beta_2)\le\left({4g_1 - 3\over m_1}+3\right)N(g_2,0,m_2).$$ 
Note that since $M$ is a simple knot manifold the surface $F_2$ cannot
be a disk or an annulus.  It then follows from \xref{Nsnvf} that
$$N(s,n,v) = \left[{12g_2 - 12\over m_2}\right]+6 .$$
\EndProof

Ian Agol has informed us that a slightly stronger estimate follows
from the techniques in his paper [\cite{Agol}].  By combining his Theorem 5.1
with the proof of his Theorem 8.1 he can show under the
hypotheses of Corollary \xref{CameronCorollary} that
$$\Delta(\beta_1,\beta_2)\le {36\over3.35}\left({2g_1 - 2\over
m_1}+1\right)\left({2g_2 - 2\over m_2}+1\right).$$
In particular the coefficient of $g_1g_2/m_1m_1$ is less than
43 for this estimate, while in the estimate provided by Corollary
\xref{CameronCorollary} the corresponding coefficient is 48.
Agol's methods depend on the rigorous computational results of
Cao and Meyerhoff [\cite{CaoMeyerhoff}].

\Corollary
\tag{HighGenusSphere}
Let $M$ be a simple knot manifold and $F\subset M$ an essential
bounded surface with boundary slope $\beta$ which is not a
semi-fiber. Set $g=\genus F$ and $m=|\partial F|$.  If $\alpha$
is the boundary slope of an essential planar surface in $M$
then
$$\D(\a,\b) \le {20g-15\over m}+15.$$
\EndCorollary

\Proof
If $M(\alpha)$ is reducible then there is an essential planar surface
$F_2$ with boundary slope $\alpha$.  The result follows from
\xref{CameronCorollary} by taking $F = F_1$, $\beta = \beta_1$,
and $\alpha=\beta_2$.  
\EndProof

Still another corollary to Theorem \xref{HighGenusTheorem} can be
obtained by using Corollary \xref{Genusg}.  The reader is invited to
formulate the statement.

\deepsection{Tight surfaces}
\tag{TightSurfaces}
The goal of this section is to prove Theorem \xref{PlanarTheorem},
which provides a major improvement on conclusion of Theorem
\xref{DottedLengthBound} in the special case where $F$ is planar,
i.e. $g=0$.  This leads to corresponding improvements to the
corollaries of the previous section in cases where the planarity
assumption hold.

The techniques used in this section make use of some variants of the
surfaces $\dot\Phi_k^\epsilon$ which are denoted
$\breve\Phi_k^\epsilon$.  We start with some preliminaries which
are needed for the definition of these surfaces.

\subsection{Perfect Surfaces}
In this subsection $S$ will denote a compact orientable surface of
negative Euler characteristic.  If $A$ is a subsurface of $S$ we will
denote the frontier of $A$ by $\frontier A$.

\Definition
\tag{Perfect}
A subsurface $A$ of $S$ will be said to be {\it perfect} if
\part{(i)} $A$ is $\pi_1$-injective;
\part{(ii)} $A$ contains $\partial S$; and
\part{(iii)} every component of $A$ contains a component of
$\partial S$.
\EndDefinition

Thus if $A$ is a perfect subsurface of $S$ then each component of $A$
is either a large subsurface of $S$ which contains a component of
$\partial S$ or a regular neighborhood of a component of $\partial
S$.

Suppose that $A$ is a large outer subsurface of a compact orientable
surface $S$.  Let $A'$ denote a subsurface in the non-ambient isotopy
class of $A$ such that for every component $C$ of $\partial A'$,
either $C\subset\partial S$, or $C\subset\int S$ and $C$ is not
parallel to any component of $\partial S$.  The surface $A'$ is unique
up to ambient isotopy.  We define a {\it perfection} of $A$, denoted
${\cal P}(A)$, to be a surface of the form $A'\cup N$ where $N$ is a
regular neighborhood of the union of all components of $\partial S$
which are not contained in $\partial A'$, and $N\cap A'=\emptyset$.
Note that ${\cal P}(A)$ is a perfect surface and that the ambient
isotopy class of ${\cal P}(A)$ is uniquely determined by the
non-ambient isotopy class of $A$.  Moreover, if $c$ is a component of
$\partial S$ and if $C$ is the component containing $c$ then
no other component of ${\cal P}(A)$ contains a curve isotopic to $c$.

\paragraph
\tag{InversePerfection}
There is an obvious inverse operation to perfection: if $B$ is a
perfect subsurface of $S$ then ${\cal L}(B)$ is a large outer subsurface of
$S$.  It is clear that if $A$ is any large outer subsurface of $S$ then
${\cal L}({\cal P}(A))$ is (non-ambiently) isotopic to $A$, and that
if $B$ is any perfect subsurface of $S$ then $B$ is ambiently
isotopic to ${\cal P}({\cal L}(B))$.  Thus we have a natural bijective
correspondence between non-ambient isotopy classes of large outer
subsurfaces of $S$ and ambient isotopy classes of perfect subsurfaces
of $S$.

The following lemma will be needed in Subsection \xref{Smiles}.

\Lemma
\tag{PerfectInvariance}
Let $P$ be a perfect subsurface of $S$ and suppose that $\tau$ is a
free involution of $P$.  Suppose that $A\subset P$ is a large outer
subsurface of $S$ which is invariant under $\tau$.  Then $A$ has a
perfection which is contained in $P$ and invariant under $\tau$.
\EndLemma

\Proof
By replacing $A$ with a slightly smaller $\tau$-invariant surface we
may assume that $A$ is contained in $\int S$.  Let $A'$ be the union
of $A$ with all of the annuli in $S$ which have one boundary component
in $\partial A$ and one boundary component in $\partial S$.  Since $P$
is perfect, $A$ is contained in $P$ and invariant under $\tau$.  Thus
$\partial S - \partial A'$ is also invariant under $\tau$.  We
define the required perfection of $A$ to be the union of $A'$ with a
suitably small $\tau$-invariant regular neighborhood of $\partial S
- \partial A'$.
\EndProof

\paragraph
If $A$ and $B$ are two perfect subsurfaces of $S$ then we define the
{\it perfect intersection} of $A$ and $B$, denoted $A\pint B$, to be
${\cal P}({\cal L}(A)\dot\lint{\cal L}(B))$.  Observe that a perfect
subsurface of $S$ is homotopic into $A\pint B$ if and only if it is
homotopic into both $A$ and $B$.  It follows from this together with
Lemma \xref{SurfaceTwo}, that $A\pint B$ is isotopic to $B\pint A$ and
that $(A\pint B)\pint C$ is isotopic to $A\pint (B\pint C)$.

Now let $A$ and $B$ be large outer subsurfaces of $S$.  It follows from the
definition of perfect intersection and the bijective correspondence
described in \xref{Perfect} that
$${\cal P}(A)\pint{\cal P}(B) = {\cal P}(A\dot\lint B)  .$$

\Proposition
\tag{GeometricPerfection}
Let $P$ and $Q$ be perfect surfaces of $S$.  Then there exist
(perfect) surfaces $P_1$ and $Q_1$ which are ambiently isotopic to $P$
and $Q$ respectively, such that the frontiers of $P_1$ and $Q_1$ meet
transversely and $P\pint Q$ is ambiently isotopic to the union of all
components of $P_1\cap Q_1$ which meet $\partial S$.
\EndProposition

\Proof
Set $A ={\cal L}(P)$ and $B = {\cal L}(Q)$.  According to Proposition
\xref{EssInt}, $A$ and $B$ are non-ambiently isotopic to surfaces
$A_0\subset \int S$ and $B_0\subset \int S$ such that $\partial A_0$
and $\partial B_0$ meet transversely and ${\cal L}(A_0\cap B_0)$ is
isotopic to $A\lint B$.  Let $A_1$ be the union of $A_0$ with all of
the annular components of $\overline{S- A_0}$ which meet
$\partial S$.  It is clear that $A_1$ is ambiently isotopic to $A$.
In the same way, using $B_0$, we define a subsurface $B_1$ which is
ambiently isotopic to $B$.  We claim that ${\cal L}(A_1\cap B_1)$ is
isotopic to $A\lint B$.  According to
\xref{EssInt} and \xref{DefinitionOfEssentialIntersection} it suffices
to show that property ($*$) of Proposition \xref{EssInt} holds with
$C$ replaced by ${\cal L}(A_1\cap B_1)$.  The ``only if'' part of
($*$) is clear because $A_1$ is isotopic to $A$ and $B_1$ is isotopic
to $B$.  Since the ``if'' part of Proposition \xref{EssInt} holds with
$C$ replaced by ${\cal L}(A_0\cap B_0)$, the claim will follow once we
show that ${\cal L}(A_0\cap B_0)\subset {\cal L}(A_1\cap B_1)$.  Since
$A_0\cap B_0\subset A_1\cap B_1$, it suffices to show that $A_1\cap
B_1$ is $\pi_1$-injective.  But since $\chi(S)<0$, any homotopically
trivial simple closed curve $\gamma \subset A_1 \cap B_1$ bounds a
unique disk $D$ in $S$; since $A_1$ and $B_1$ are $\pi_1$-injective
we must have $D\subset A_1$ and $D\subset B_1$ and hence $D\subset
A_1\cap B_1$.  This proves the claim.  

Thus we may take $A\lint B$ to be equal to ${\cal L}(A_1\cap B_1)$.
Let $W$ denote the union of the large components of $A_1\cap B_1$ that
meet $\partial S$.  We next assert that $W$ is $A\dot\lint B$.  To
prove this it suffices to show that every large outer component $X$ of
$A_1\cap B_1$ contains a component of $\partial S$.  By the definition
of an outer component (see \xref{OuterParts}), $X$ contains a closed
curve $\gamma\subset \int S$ which is the frontier of an annulus
$\alpha\subset S$.  Since $\gamma \subset A_1$ and $\gamma\subset
B_1$, it follows from the construction of $A_1$ and $B_1$ that
$\alpha$ is contained in both $A_1$ and $B_1$, and hence in $X$.  In
particular, $X$ contains a component of $\partial S$ as required.

Now let $N$ be regular neighborhood of $\partial S$.  We may assume
$N$ to be chosen so that each component of $N$ is disjoint from
the frontiers of $A_1$ and $B_1$.  We set $P_1 = A_1\cup N$ and
$Q_1=B_1\cup N$.  It is clear that $P_1$ and $Q_1$ are isotopic to $P$
and $Q$ respectively.  It is also clear that the union $Z$ of the
components of $P_1\cap Q_1$ which meet $\partial S$ is equal to $W\cup
N$.  Since $W$ is $A\dot\lint B$, we have that $Z= {\cal P}(W) = {\cal
P}(A\dot\lint B)$.  But by definition we have that $P\pint Q = {\cal
P}(A\dot\lint B)$.  This completes the proof.
\EndProof

\paragraph
\tag{UniqueComponent}
Suppose that $P$ and $Q$ are perfect
subsurfaces of $S$.  Let $C$ be a component of the perfect
intersection $P\pint Q$.  Since $C$ is isotopic to a subsurface of
$P$, there must be a component $P_0$ of $P$ such that $C$ is
isotopic to a subsurface of $P_0$.  We will say in this situation
that $C$ is {\it isotopically contained in} $P_0$. Let $c$ be a
component of $\partial S$ which is contained in $C$.  Then, since
there is a unique component of $P$ which contains a curve isotopic
to $c$, it follows that $P_0$ is the unique component of
$P$ which contains a surface isotopic to $C$.  Thus each component of
$P\pint Q$ is isotopically contained in a unique component of $P$.

\Definition
Let $A\subset S$ be a perfect surface.  A component $A_0$ if $A$ said
to be {\it tight} if $A$ is planar and the frontier of $A$ in $S$ is a
simple closed curve.  We define the {\it size} of a tight component
$A_0$ of $A$ to be the number of components of $\partial S$ which are
contained in $A_0$.  We will denote the size of $A_0$ by $s(A_0)$.
\EndDefinition

\Proposition
\tag{TightIntersect}
Let $P$ and $Q$ be perfect subsurfaces of $S$ and let $P_0$ be a tight
component of $P$.  Assume that every tight component of $Q$ has size
at least $s(P_0)$.  Then every component of $P\pint Q$ which is
isotopically contained in $P_0$ is tight and has size at most $s(P_0)$.
Furthermore if $P_0$ contains only one component of $P\pint Q$
then this component is isotopic to $P_0$.
\EndProposition

\Proof
By Proposition \xref{GeometricPerfection} we may assume that $P$ and
$Q$ have been chosen within their isotopy classes so that the
frontiers of $P$ and $Q$ meet transversely and $Z = P\pint Q$ is the
union of all components of $P\cap Q$ which meet $\partial S$.

Consider first the case in which $P_0$ contains at least one component
of $\frontier Q$.  Note that since $Q$ is perfect, the components of
$\frontier Q$ are homotopically non-trivial simple closed curves in
$S$.  Since $P_0$ is planar and has connected frontier, every
homotopically non-trivial simple closed curve $\gamma$ in $\int P_0$
is the frontier in $S$ of a unique subsurface $W_\gamma$ of $P_0$; the
non-triviality of $\gamma$ implies that $W_\gamma$ is not a disk and
hence that $W_\gamma\cap\partial S \not=\emptyset$.  Among all
components of $\frontier Q$ contained in $P_0$ we choose one,
$\gamma_0$, such that $W_{\gamma_0}$ is minimal with respect to
inclusion.  The minimality implies that $W_{\gamma_0}$ is either a
component of $Q$ or of $\overline{S- Q}$.  But since $Q$ is
perfect, we have $\partial S\subset Q$ and hence $W_{\gamma_0}\cap
Q\supset W_{\gamma_0}\cap\partial S \not=\emptyset$.  Hence
$W_{\gamma_0}$ must be a component of $Q$.  Since $W_{\gamma_0}$ is
contained in the planar surface $P_0$ and has connected frontier, it
is in fact a tight component of $Q$ with $s(W_{\gamma_0})\le s(P_0)$.
On the other hand the hypothesis of the proposition implies that
$s(W_{\gamma_0})\ge s(P_0)$.  It now follows that $P_0$ is a regular
neighborhood of $W_{\gamma_0}$.  Clearly $W_{\gamma_0}$ is a component
of $P\cap Q$ and hence of $Z$.  Since the annulus $\overline{P_0
- W_{\gamma_0}}$ is disjoint from $\partial S$ it follows that
$W_{\gamma_0}$ is the only component of $P_0\cap Q$ which meets
$\partial S$.  Hence $W_{\gamma_0}$ is the only component of $Z$ which
is contained in $P_0$.  Since $W_{\gamma_0}$ is isotopic to its
regular neighborhood $P_0$, both conclusions of the proposition are
established in this case.

There remains the case in which each component of $P_0\cap\frontier Q$
is a properly embedded arc in $P_0$ having both endpoints on
$\frontier P_0$.  It follows that any component of $P_0\cap Q$ is a
planar surface whose frontier in $S$ is a simple closed curve.  In
particular, if $C$ is a component of $P_0\cap Q$ which contains a
component of $\partial S$ then $C$ is a tight component of $Z$.  It is
also clear that $s(C)
\le s(P_0)$.  To prove the last assertion of the proposition in this
case, assume that $C$ is the only tight component of $Z$ which is
contained in $P_0$.  Then the frontier in $P_0$ of $C$ is a collection
of properly embedded arcs which are parallel to subarcs of $\frontier
P_0$.  It then follows that $C$ is isotopic to $P_0$.
\EndProof

\Proposition
\tag{HasTightComponent}
Suppose that $S$ is planar and let $P \not= S$ be a perfect subsurface
of $S$.  Then $P$ has a tight component.
\EndProposition

\Proof
Since $P \not= S$ we have $\frontier P \not= \emptyset$.  Since $P$ is
perfect each component of $\frontier P$ is a homotopically non-trivial
curve in $S$ and hence is the frontier of two planar subsurfaces of
$S$, neither of which is a disk.  Among all subsurfaces $A$ of $S$
such that $\frontier A$ consists of a single component of $\frontier
S$ we choose one, say $A_0$ which is minimal with respect to
inclusion.  Since $A_0$ is not a disk we have $A_0\cap\partial S
\not=\emptyset$.  The minimality implies that either $A_0$ is a
component of $P$ or of $\overline{S-P}$.  But since $P$ is
perfect, we have $\partial S\subset P$ and hence $A_0\cap P\supset
A_0\cap\partial S \not=\emptyset$.  Thus $A_0$ is a component of $P$
and by definition is tight.
\EndProof

We record here a simple lemma that will be needed in the next
subsection.

\Lemma
\tag{AnnulusLemma}
Let $A$ be an oriented annulus and $\alpha$ a component of $\partial
A$. Let $f$ and $g$ be two embeddings of $A$ into an orientable
surface $F$.  Suppose that $f(\alpha) = g(\alpha) = c$ where $c$ is a
component of $\partial F$, and that $f$ and $g$ carry the orientation
of $A$ to the same orientation of $F$. Then $f$ and $g$ are homotopic.
\EndLemma
\NoProof

\subsection{Reduced homotopies and perfect subsurfaces}
\tag{Smiles}

\paragraph
\tag{SmileHypotheses}
Throughout Subsection \xref{Smiles} we will assume that $M$ is a
simple knot manifold and that $\tilde F$ is a splitting surface in $M$
which admits a long rectangle.  Since this is the same assumption that
was made in \xref{Hypotheses}, the results of Subsection
\xref{DottedEssHom} may be applied in this subsection.

We will fix $\dot\Phi_k^{\pm1}$, $\dot h_k^{\pm1}$ and
$\dot\tau_{\pm1}$ as in Subsection
\xref{Dottedh}.  Recall that the surfaces $\dot\Phi_k^\epsilon$
are only defined up to non-ambient isotopy.  Here we suppose each
surface $\dot\Phi_k^\epsilon$ to have been normalized within its
non-ambient isotopy class so that if $C$ is a boundary component of
$\dot\Phi_k^\epsilon$ then either $C\subset \partial
\tilde F$, or $C\subset \int \tilde F$ and $C$ is not parallel to
any component of $\partial \tilde F$.  (It is clear that the surfaces
$\dot\Phi_k^\epsilon$ can be chosen to have this property in addition
to having the nestedness property stated in \xref{Dottedh}.)

By definition $\dot\Phi_k^\epsilon$ is a large outer subsurface for
each $\epsilon\in\{\pm 1\}$ and $k\ge 0$.  Using the notation of
\xref{Perfect} we set $$\breve\Phi_k^\epsilon = {\cal
P}(\dot\Phi_k^{\epsilon}).$$

Because of the way that the surfaces $\dot\Phi_k^\epsilon$ have
been normalized, ${\cal P}(\dot\Phi_k^\epsilon)$ is the disjoint union
of $\dot\Phi_k^\epsilon$ with ${\cal A}_k^\epsilon$, where ${\cal
A}_k^\epsilon$ is a regular neighborhood of the union of all
components of $\partial \tilde F$ which are not contained in
$\dot\Phi_k^{\pm1}$.

We may assume that the regular
neighborhoods ${\cal A}_k^\epsilon$ have been chosen so that
$\breve\Phi_k^\epsilon \supset {\cal A}_{k+1}^\epsilon$.
This means that for each $\epsilon\in\{\pm1\}$ we have
$$\tilde F=\breve\Phi_0^\epsilon \supset \breve\Phi_1^\epsilon \supset
\breve\Phi_2^\epsilon\supset
\cdots .$$
Note also that by \xref{InversePerfection} we have  ${\cal
L}(\breve\Phi_k^\epsilon)=\dot\Phi_k^\epsilon$ for each
$\epsilon\in\{\pm1\}$ and each $k\ge 0$.

We denote by $m$ the number of boundary components of $\tilde F$.
Since $\tilde F$ is a splitting surface the integer $m$ is even.  We
index the components of $\partial\tilde F$ as $c_t$ where $t$ ranges
over $\Z/m\Z$.  If $q$ is an integer we will denote the image of
$q$ in $\Z/m\Z$ by $\bar q$.

We assume that the indexing of components of $\partial \tilde F$ has
been done in such a way that for each $t\in\Z/m\Z$ the curves $c_t$
and $c_{t+1}$ cobound an annulus $R_t\subset \partial M$ whose
interior is disjoint from $\partial \tilde F$.  We may assume further
that the indexing is done in such a way that $R_{\bar q}\subset
M_{\tilde F}^+$ for every even integer $q$ and $R_{\bar q}\subset
M_{\tilde F}^-$ for every odd integer $q$.

\paragraph
\tag{Permutations}
For every integer $k\ge 0$ and $\epsilon\in\{\pm1\}$ we define a
permutation $\sigma_k^\epsilon$ of $\Z/m\Z$ by
$$\sigma_k^\epsilon(\bar q) = \overline {q+(-1)^q\epsilon k} .$$ (Since $m$ is
even, the coefficient $(-1)^q$ is determined by the congruence class
$\bar q$.)

We observe that if $i$ and $j$ are non-negative integers with $i+j
=k$ then for each $\epsilon\in\{\pm 1\}$ we have
$$\sigma_k^\epsilon =
\sigma_j^{(-1)^i\epsilon}\circ\sigma_i^\epsilon.
\leqno{(\subsubnumber)}$$
\tag{Composition}
We also observe that for each $k\ge 0$ and each $\epsilon\in\{\pm 1\}$
we have
$$(\sigma_k^\epsilon)^{-1} =
\sigma_k^{(-1)^k\epsilon}.\leqno{(\subsubnumber)}$$ From $(1)$ and
$(2)$ it follows that for each $k\ge 0$ and each $\epsilon\in\{\pm
1\}$ we have $$\sigma_{2k+1}^\epsilon =
(\sigma_k^\epsilon)^{-1}\circ\sigma_1^{(-1)^k\epsilon}\circ\sigma_k^\epsilon
. \leqno{(\subsubnumber)}$$\tag{Conjugation} 
It is clear from \xref{Conjugation} that $\sigma_k^\epsilon$ is a free involution
if $k$ is odd.

\Lemma
\tag{BoundaryAction}
Let $k$ be a positive integer and fix $\epsilon\in\{\pm1\}$.  If a
component $c_t$ of $\partial \tilde F$ is contained in
$\dot\Phi_k^\epsilon$ for some $t\in \Z/m\Z$, then $\dot
h_k^\epsilon(c_t) = c_{\sigma_k^\epsilon(t)}$.
\EndLemma

\Proof
The lemma is trivial in the case $k=0$.

Consider the case $k=1$.  By Lemma \xref{ArcLemma}, together with our
normalization of $\dot\Phi_1^\epsilon$, we know that every component
of $\partial \tilde F$ is contained in $\dot\Phi_1^\epsilon$ and that
$h_1^\epsilon = \tau_\epsilon$ interchanges two components $c$ and
$c'$ of $\partial \tilde F$ if and only if $c$ and $c'$ cobound an
annulus component of $M_{\tilde F}^\epsilon\cap M$.  The definition of
$\sigma_1^\epsilon$ thus implies that for every $t\in \Z/m\Z$ we have
$h_1^\epsilon(c_t) = \tau_\epsilon(c_t) = c_{\sigma_1^\epsilon(t)}$.

Now, arguing inductively, we assume that $k>1$, that $c_t$ is a
component of $\partial \tilde F$ contained in $\dot\Phi_{k}^\epsilon
\subset \dot\Phi_{k-1}^\epsilon$ and that $\dot h_{k-1}^\epsilon(c_t)
= c_{\sigma_{k-1}^\epsilon(t)}$.  Applying Proposition
\xref{DottedPropertyThree} with $i=k-1$ and $j=1$ we see $\dot
h_{k-1}^\epsilon|\dot\Phi_k^\epsilon$ is homotopic in $\tilde F$ to an
embedding $\dot g_{k-1}^\epsilon:\dot
\Phi_{k}^\epsilon \to \dot \Phi_1^{(-1)^{k-1}\epsilon}$
such that $\dot h_k^\epsilon$ is homotopic to $\dot
h_1^{(-1)^{k-1}\epsilon}\circ \dot g_{k-1}^\epsilon$.  Since
$\dot h_{k-1}^\epsilon(c_t) =
c_{\sigma_{k-1}^\epsilon(t)}$, 
the curve $\dot g_{k-1}^\epsilon(c_t)\subset \dot
\Phi_1^{(-1)^{k-1}\epsilon}$ is homotopic in $\tilde F$ to the boundary
component $c_{\sigma_{k-1}^\epsilon(t)}$ of $\tilde F$.  It follows
from the normalization of $\Phi_1^{-\epsilon}$ in
\xref{SmileHypotheses} that $c_{\sigma_{k-1}^\epsilon(t)}$ is a
boundary component of $\Phi_1^{(-1)^{k-1}\epsilon}$, and that $\dot
g_{k-1}^\epsilon(c_t)$ is homotopic in $\Phi_1^{(-1)^{k-1}\epsilon}$ to
$c_{\sigma_{k-1}^\epsilon(t)}$.  Hence $\dot h_k^\epsilon(c_t)$ is
homotopic in $\tilde F$ to
$\dot h_1^{(-1)^{k-1}\epsilon}( c_{\sigma_{k-1}^\epsilon(t)})$ which by the
case $k=1$ of the lemma is equal to
$$c_{\sigma_1^{(-1)^{k-1}\epsilon}\circ \sigma_{k-1}^\epsilon(t)} =
c_{\sigma_k^\epsilon(t)}.$$ 

Since the boundary component $\dot
h_k^\epsilon(c_t)$ of $\dot h_k^\epsilon(\dot\Phi_k^{\epsilon}) =
\dot\Phi_k^{(-1)^{k+1}\epsilon}$ is homotopic to the component
$c_{\sigma_k^\epsilon(t)}$ of $\partial \tilde F$, it follows from the
normalization of $\dot\Phi_k^{(-1)^{k+1}\epsilon}$ that $\dot
h_k^\epsilon(c_t)=c_{\sigma_k^\epsilon(t)}$.
\EndProof

\Lemma
\tag{HSmileExists}
For each
$\epsilon\in\{\pm 1\}$ and each $k\ge 0$, there exists a homeomorphism $\breve
h_k^\epsilon:\breve\Phi_k^\epsilon \to
\breve\Phi_k^{(-1)^{k+1}\epsilon}$,
such that
\part{$(1)$}  the restriction of
$\breve h_k^\epsilon$
to ${\cal L}(\breve\Phi_k^\epsilon)=\dot\Phi_k^\epsilon$ is $\dot
h_k^\epsilon$;
\part{$(2)$} if $\tilde F$ is given a consistent orientation
$\breve h_k^\epsilon:\breve\Phi_k^\epsilon\to\tilde F$ reverses orientation
if $k$ is odd and preserves orientation if $k$ is even; and
\part{$(3)$} for each $t\in \Z/m\Z$ we have
$\breve h_k^\epsilon(c_t)=c_{\sigma_k^\epsilon(t)}$.

The homeomorphism $\breve h_k^\epsilon$ is determined up to
isotopy by the properties $(1)$-$(3)$.  Furthermore we may choose
$\breve h_1^\epsilon$ within its isotopy class so that it is a free
involution.
\EndLemma

\Proof
Fix an orientation of $\breve\Phi_k^\epsilon$ which is induced from a
consistent orientation of $\tilde F$.
By Lemma \xref{BoundaryAction}
we know that if $c_t \subset \dot\Phi_k^\epsilon$ then
$\dot h_k^\epsilon(c_t) = c_{\sigma_k^\epsilon(t)}$.  In particular, the
correspondence $c_t \rightarrow c_{\sigma_k^\epsilon(t)}$
restricts to a bijection between the components of 
$\partial \tilde F\cap
\dot\Phi_k^\epsilon$ and those of $\partial \tilde F\cap
\dot\Phi_k^{(-1)^{k+1}\epsilon}$.  It therefore also restricts to a bijection
between the components of $\partial \tilde F - \partial
\dot\Phi_k^\epsilon$ and those of $\partial \tilde F - \partial
\dot\Phi_k{(-1)^{k+1}\epsilon}$.  Now $\breve\Phi_k^\epsilon$ is the
union of $\dot\Phi_k^\epsilon$ with the regular neighborhood
${\cal A}_k^\epsilon$ of
$\partial \tilde F - \partial
\dot\Phi_k^\epsilon$, and 
$\breve\Phi_k^{(-1)^{k+1}\epsilon}$ is the union of
$\dot\Phi_k^{(-1)^{k+1}\epsilon}$ with the regular neighborhood ${\cal
A}_{k+1}^{(-1)^{k+1}\epsilon}$ of $\partial \tilde F -
\partial \dot\Phi_k^{(-1)^{k+1}\epsilon}$.
Let $f_k^\epsilon$ be a homeomorphism from ${\cal A}_k^\epsilon$ to
${\cal A}_{k+1}^{(-1)^{k+1}\epsilon}$ which maps each component $c_t$
of $\partial \tilde F\cap{\cal A}_k^\epsilon$ to the component
$c_{\sigma_k^\epsilon(t)}$ of $\partial \tilde F\cap{\cal
A}_{k+1}^{(-1)^{k+1}\epsilon}$.  Since each component of ${\cal
A}_k^\epsilon$ is an annulus we may choose $f_k^\epsilon$ to be an
orientation-reversing embedding of ${\cal A}_k^\epsilon$ into $\tilde
F$ if $k$ is odd, and an orientation-preserving embedding if $k$ is
even.  These conditions determine $f_k$ up to isotopy.  We define
$\breve h_k^\epsilon$ to be the homeomorphism whose restriction to
$\dot\Phi_k^\epsilon$ is $\dot h_k^\epsilon$ and whose restriction to
${\cal A}_k^\epsilon$ is $f_k$.  Conditions ($1$) and ($3$) hold by
construction.  To see that condition (2) holds it suffices to observe
that by \xref{Dottedh} the embedding $\dot
h_k^\epsilon:\dot\Phi_k^\epsilon\to\tilde F$ reverses orientation if
$k$ is odd and preserves orientaion if $k$ is even, and that
$f_k^\epsilon$ has the same property by construction.

Since we have observed that $f_k^\epsilon$ is determined up to isotopy
by its stated properties, it follows that $\breve h_k^\epsilon$ is
determined up to isotopy by conditions (1)-(3).  Finally, since
$\sigma_1^\epsilon$ is a free involution it is clear that we may
choose $f_1^\epsilon$ within its isotopy class so that it is a free
involution.  Since $\dot h_1^\epsilon$ is a free involution by
\xref{Dottedh} it follows that $\breve h_1^\epsilon$ is a free
involution.
\EndProof

\paragraph
For the rest of Subsection \xref{Smiles} we will fix homeomorphisms
$\breve h_k^\epsilon$ satisfying the conclusions of Lemma 
\xref{HSmileExists}.
The free involution $\breve h_1^\epsilon$ will sometimes be denoted
$\breve\tau_\epsilon$.  Note that $\breve\tau_\epsilon$ is an
extension of the free involution $\dot\tau_\epsilon$ defined in
\xref{Dottedh}.  

\Lemma
\tag{PerfectImage}
Let $k$ and $i$ be integers with $k\ge i\ge 0$.
Then $\breve h_i^\epsilon(\breve\Phi_k^\epsilon)$ is ambiently isotopic
to ${\cal P}(\dot h_i^\epsilon(\dot\Phi_k^\epsilon))$.
\EndLemma

\Proof
Since $\breve\Phi_k^\epsilon$ and $\breve\Phi_i^\epsilon$ are perfect
surfaces, and since $\breve
h_i^\epsilon:\breve\Phi_i^\epsilon\to\tilde F$ is an embedding which
maps $\partial \tilde F$ onto $\partial \tilde F$ (see Lemma
\xref{BoundaryAction}), it follows that $\breve
h_i^\epsilon(\breve\Phi_k^\epsilon)$ is a perfect subsurface of
$\tilde F$.
On the other hand, by
\xref{InversePerfection} and \xref{SmileHypotheses} we have
$\dot\Phi_k^\epsilon = {\cal L}(\breve\Phi_k^\epsilon)$.
Since $\breve h^\epsilon_i$ is a
homeomorphism, $\dot h^\epsilon_i(\dot\Phi_k^\epsilon)=
\breve h^\epsilon_i(\dot\Phi_k^\epsilon)=
{\cal L}(\breve h^\epsilon_i(\breve\Phi_k^\epsilon))$.
Since $\breve h^\epsilon_i(\breve\Phi_k^\epsilon)$
is perfect, it follows from \xref{InversePerfection} that
$\breve h_i^\epsilon(\breve\Phi_k^\epsilon)$ is ambiently isotopic
to ${\cal P}(\dot h_i^\epsilon(\dot\Phi_k^\epsilon))$.
\EndProof

The next four results,
\xref{SmileyPropertyThree} -- \xref{SmileyStrictContainment},
are analogues for the surfaces $\breve\Phi_k^\epsilon$ of Propositions
\xref{DottedPropertyThree}, \xref{DottedPropertyTwo},
\xref{DottedOddK} and \xref{DottedStrictContainment}.

\Proposition
\tag{SmileyPropertyThree}
Let $i$ and $j$ be non-negative integers, and set $k=i+j$. Then for
each $\epsilon\in\{\pm1\}$, the map
$\breve h_i^\epsilon|\breve\Phi_k^\epsilon$ is homotopic in $\tilde
F$, rel $\partial \tilde F$, 
to an embedding $\breve g_i^\epsilon\colon \breve\Phi_k^\epsilon\to
\breve\Phi_{j}^{(-1)^{i}\epsilon}$ such that
$\breve h_{j}^{(-1)^{i}\epsilon}\circ \breve g_i^\epsilon$ is
homotopic in $\tilde F$ to $\breve
h_k^\epsilon$.
\EndProposition

\Proof
Let $\dot g_i^\epsilon\colon \dot\Phi_k^\epsilon\to
\dot\Phi_{j}^{(-1)^{i}\epsilon}$ be given by 
Proposition \xref{DottedPropertyThree}.  Let $c_t$ be a boundary
component of $\tilde F$ that is contained in $\dot\Phi_k^\epsilon$.
We have $\dot h_i^\epsilon(c_t)=c_{\sigma_i^\epsilon(t)}$.  Since
$\dot g_i^\epsilon$ is homotopic to $\dot h_i^\epsilon$, the map $\dot
g_i^\epsilon|c_t: c_t \to \dot\Phi_{j}^{(-1)^{i}\epsilon}$ is
homotopic in $\tilde F$ to $\dot h_i^\epsilon|c_t$.  Because of the
way that the surface $\dot\Phi_{j}^{(-1)^{i}\epsilon}$ has been
normalized, this means that $c_{\sigma_i^\epsilon(t)}$ is a boundary
curve of $\dot\Phi_{j}^{(-1)^{i}\epsilon}$ and that $\dot
g_i^\epsilon|c_t$ is isotopic to $h_i^\epsilon|c_t$ in
$\dot\Phi_{j}^{(-1)^{i}\epsilon}$.  Therefore, after modifying $\dot
g_i^\epsilon$ by a non-ambient isotopy, we may assume that, for each
$c_t$ contained in $\partial \dot\Phi_k^\epsilon$, we have $\dot
g_i^\epsilon(c_t) = c_{\sigma_i^\epsilon(t)}$ and $\dot
g_i^\epsilon|c_t = \dot h_i^\epsilon|c_t$.  Since the maps $\dot
g_i^\epsilon$ and $\dot h_i^\epsilon$ are homotopic and agree on
$\partial\tilde F\cap\dot\Phi_k^\epsilon$, they are homotopic rel
$\partial \tilde F$.

The homeomorphism $\breve h_i^\epsilon$ maps ${\cal A}_k^\epsilon$ to
a regular neighborhood of a collection of boundary curves of $\tilde
F$.  Since $\dot g_i^\epsilon$ agrees with $\breve h_i^\epsilon$ on
$\dot\Phi_k^\epsilon\cap \partial \tilde F$, each boundary
component of $\breve h_i^\epsilon({\cal A}_k^\epsilon)\cap \partial
\tilde F$ is disjoint from the image of $\dot g_i^\epsilon$.  Thus
$\breve h_i^\epsilon|{\cal A}_k^\epsilon$ is isotopic rel $\partial
\tilde F$ to an embedding
$f:{\cal A}_k^\epsilon\to \dot\Phi_{j}^{(-1)^{i}\epsilon}$ whose image
is disjoint from the image of $\dot g_i^\epsilon$.

We define $\breve g_i^\epsilon$ so that $\breve g_i^\epsilon|{\cal
A}_k^\epsilon = f$ and $\breve g_i^\epsilon|\dot\Phi_k^\epsilon = \dot
g_i^\epsilon$.  Then
$\breve g_i^\epsilon$ is homotopic rel $\partial\tilde F$ to  
$\breve h_i^\epsilon|\breve\Phi_k^\epsilon$.  Moreover we have
$\breve g_i^\epsilon(c_t) =
c_{\sigma_i^\epsilon(t)}$ for every $t\in \Z/m\Z$.  Since
$\sigma_k^\epsilon = \sigma_j^{(-1)^i\epsilon}\circ\sigma_i^\epsilon$
by \xref{Composition}, we have $\breve
h_{j}^{(-1)^{i}\epsilon}\circ \breve g_i^\epsilon(c_t) =
\breve h_k^\epsilon(c_t)$ for every $t\in
\Z/m\Z$.

Let $\tilde F$ be given a consistent orientation.  Since the
embeddings $\breve g_i^\epsilon$ and $\breve
h_i^\epsilon|\breve\Phi_k^\epsilon$ of $\breve\Phi_k^\epsilon$ into
$\tilde F$ are homotopic rel $\partial\tilde F$, they both reverse
orientation if $i$ is odd and preserve orientation if $i$ is even.  In
particular, since $i+j=k$, the embeddings $\breve
h_{j}^{(-1)^{i}\epsilon}\circ \breve g_i^\epsilon|{\cal A}_k$ and
$\breve h_k^\epsilon|{\cal A}_k$ are either both orientation
preserving or both orientation reversing.  Applying Lemma
\xref{AnnulusLemma} to the restrictions of $\breve
h_{j}^{(-1)^{i}\epsilon}\circ \breve g_i^\epsilon$ and $h_k^\epsilon$
to each component of ${\cal A}_k^\epsilon$ we conclude that $\breve
h_{j}^{(-1)^{i}\epsilon}\circ \breve g_i^\epsilon|{\cal A}_k^\epsilon$
is homotopic to $\breve h_k^\epsilon|{\cal A}_k^\epsilon$.  On the
other hand, according to Proposition
\xref{DottedPropertyThree}, $\dot h_{j}^{(-1)^{i}\epsilon}\circ \dot
g_i^\epsilon$ is homotopic to $\dot h_k^\epsilon$. Hence $\breve
h_{j}^{(-1)^{i}\epsilon}\circ \breve g_i^\epsilon$ is homotopic to
$h_k^\epsilon$.
\EndProof

\Proposition
\tag{SmileyPropertyTwo}
Let $i$ and $j$ be non-negative integers, and set $k=i+j$.  Then for
each $\epsilon \in \{\pm1\}$ the subsurface $\breve h_i^\epsilon(\breve
\Phi_k^\epsilon)$ is ambiently isotopic in $\tilde F$ to the perfect
intersection
$\breve \Phi_{i}^{(-1)^{i+1}\epsilon}\pint \breve\Phi_{j}^{(-1)^i\epsilon}$.
\EndProposition

\Proof
By definition we have 
$$\breve \Phi_{i}^{(-1)^{i+1}\epsilon}\pint
\breve\Phi_{j}^{(-1)^i\epsilon} =
{\cal P}\left({\cal L}(\breve\Phi_{i}^{(-1)^{i+1}\epsilon})\dotlint
         {\cal L}(\breve\Phi_{j}^{(-1)^i\epsilon})\right) =
{\cal P}(\dot\Phi_{i}^{(-1)^{i+1}\epsilon}\dotlint\dot\Phi_{j}^{(-1)^i\epsilon}) .$$
Combining this with Proposition \xref{DottedPropertyTwo} we conclude that
$\breve\Phi_{i}^{(-1)^{i+1}\epsilon}\pint\breve\Phi_{j}^{(-1)^i\epsilon}$
is equal to
${\cal P}(\dot h_i^\epsilon(\dot\Phi_k^\epsilon))$ which, according to
Lemma \xref{PerfectImage}, is ambiently isotopic to $\breve
h_i^\epsilon(\breve\Phi_k^\epsilon)$.
\EndProof

\Proposition
\tag{SmileyOddK}
For any non-negative integer $k$ and for each $\epsilon\in\{\pm1\}$
the surface $\breve h^\epsilon_k(\breve\Phi_{2k+1}^\epsilon)$ is
ambiently isotopic in $\tilde F$ to a subsurface of
$\breve\Phi_1^{(-1)^k\epsilon}$ which is invariant under the free
involution $\breve\tau_{(-1)^k\epsilon}$.  In particular,
$\breve\Phi_{2k+1}^\epsilon$ admits a free involution
which maps $c_t$ to $c_{\sigma_{2k+1}^\epsilon(t)}$ for each $t\in\Z/m\Z$,
and is orientation-reversing as an embedding of 
$\breve\Phi_{2k+1}^\epsilon$ into $\tilde F$, if $\tilde F$ is
given a consistent orientation.
\EndProposition

\Proof
According to Proposition \xref{DottedOddK} there is a subsurface
$A$ of $\dot\Phi_1^{(-1)^k\epsilon}\subset
\breve\Phi_1^{(-1)^k\epsilon}$ which is invariant under 
$\dot\tau_{(-1)^k\epsilon}=
\breve\tau_{(-1)^k\epsilon}|\dot\Phi_1^{(-1)^k\epsilon}$ and isotopic to 
$\dot h^\epsilon_k(\dot\Phi_{2k+1}^\epsilon)$.  We apply Lemma
\xref{PerfectInvariance}, taking $S=\tilde F$, $P =
\breve\Phi_1^{(-1)^k\epsilon}$, and $\tau =
\breve\tau_{(-1)^k\epsilon}$.  We conclude that $A$ has a perfection
$Q \subset \breve\Phi_1^{(-1)^k\epsilon}$ which is invariant under
$\breve\tau_{(-1)^k\epsilon}$.  Since $Q$ is a perfection of $\dot
h^\epsilon_k(\dot\Phi_{2k+1}^\epsilon)$, it follows from Lemma
\xref{PerfectImage} that $Q$ is ambiently isotopic to $\breve
h^\epsilon_k(\breve\Phi_{2k+1}^\epsilon)$.  This completes the proof of
the first assertion.

Let $h:\breve\Phi_{2k+1}^\epsilon \to \tilde F$ be an embedding which
is ambiently isotopic to $\breve
h^\epsilon_k|\breve\Phi_{2k+1}^\epsilon$ and maps
$\breve\Phi_{2k+1}^\epsilon$ onto $Q$.  We define a free involution
$\tau$ of $\breve\Phi_{2k+1}^\epsilon$ by $\tau =
h^{-1}\circ\breve\tau_{(-1)^k\epsilon}\circ
h$. Thus $\tau(c_t) = c_s$ where $s = (\sigma_k^\epsilon)^{-1}\circ
\sigma_1^{(-1)^k\epsilon}\circ\sigma_k^\epsilon(t)$.
It follows from \xref{Conjugation} that $s = \sigma_{2k+1}^\epsilon(t)$.

Let $\tilde F$ be given a consistent orientation.  
Lemma \xref{HSmileExists} implies that $h$ preserves orientation if $k$ is
odd, that $h$ reverses orientation if $k$ is even, and that
$\breve \tau_{(-1)^k\epsilon}$ reverses orientation.  It follows
that $\tau$ is orientation-reversing as an embedding of 
$\breve\Phi_{2k+1}^\epsilon$ into $\tilde F$.
\EndProof

\Proposition
\tag{SmileyStrictContainment}
Let $k$ be a non-negative integer and let $\epsilon\in\{\pm1\}$ be
given.  If $\breve\Phi_k^\epsilon$ and $\breve\Phi_{k+2}^\epsilon$ are
isotopic in $\tilde F$ then either $\breve\Phi_k^\epsilon$ is a
regular neighborhood of $\partial \tilde F$ or $\tilde F$ is a
semi-fiber.
\EndProposition

\Proof
By definition we have that $\breve\Phi_k^\epsilon = {\cal
P}(\dot\Phi_k^\epsilon)$ and $\breve\Phi_{k+2}^\epsilon = {\cal
P}(\dot\Phi_{k+2}^\epsilon)$.  Thus it follows from
\xref{InversePerfection} that $\breve\Phi_k^\epsilon$ is isotopic to
$\breve\Phi_{k+2}^\epsilon$ if and only if $\dot\Phi_k^\epsilon$ is
isotopic to $\dot\Phi_{k+2}^\epsilon$, and that
$\breve\Phi_k^\epsilon$ is a regular neighborhood of $\partial \tilde
F$ if and only if $\dot\Phi_k^\epsilon$ is empty.  The result
therefore follows from Proposition
\xref{DottedStrictContainment}. 
\EndProof

\subsection{Very tight surfaces.}

In this subsection we assume that $M$ is a simple knot manifold, and
that $\tilde F$ is a splitting surface for $M$ which admits a long
rectangle.

For each integer $k\ge 0$ and each $\epsilon\in\{\pm 1\}$ we will define
the perfect subsurfaces $\breve \Phi_k^\epsilon = {\cal A}_k^\epsilon
\cup \dot\Phi_k^\epsilon$, with
$\tilde F=\breve\Phi_0^\epsilon \supset \breve\Phi_1^\epsilon \supset
\breve\Phi_2^\epsilon\supset
\cdots$, as in Subsection \xref{Smiles}.

We set $m=|\partial\tilde F|$ and we assume that the components of
$\partial \tilde F$ have been indexed by elements of $\Z/m\Z$ as
described in \xref{Smiles}.  We also define the permutations
$\sigma_k^\epsilon$ as in \xref{Smiles}, and for each integer $k\ge 0$
and $\epsilon\in\{\pm 1\}$ we fix a homeomorphism $h_k^\epsilon$
satisfying the conclusions of Lemma \xref{HSmileExists}.

For $\epsilon\in\{\pm1\}$ we denote by ${\cal T}^\epsilon$ the set of
tight components of $\breve \Phi_1^\epsilon$ We define $s_0$ to be the
infimum of $s(C)$ as $C$ ranges over ${\cal T}^+\cup{\cal T}^-$.
Thus, if ${\cal T}^+ = {\cal T}^- = \emptyset$ then $s_0=+\infty$ and
otherwise $s_0$ is a strictly positive integer.  We will say that a
perfect subsurface of $\tilde F$ is {\it very tight} if it is tight
and has size at most $s_0$.  (In particular, if ${\cal T}^+ = {\cal
T}^- = \emptyset$ then any tight surface is very tight.)

If $A$ is a perfect subsurface of $\tilde F$
then we define $\VT(A)$ to be the union of the very tight components
of $A$.  

\Lemma
\tag{NumVTEven}
If $\epsilon\in\{\pm1\}$ and $k> 0$ is odd, then
$|\VT(\breve\Phi_k^\epsilon)|$ is even.
\EndLemma
\Proof
We give $\breve\Phi_k^\epsilon$ the orientation inherited from a
consistent orientation of $\tilde F$.  Since $k$ is odd, Proposition
\xref{SmileyOddK} implies that the surface $\breve\Phi_k^\epsilon$
admits an orientation-reversing free involution $\tau$ which permutes
the components of $\partial \tilde F$.  It follows that if $T$ is a
very tight component of $\breve\Phi_k^\epsilon$ then $\tau(T)$ is also
a very tight component of $\breve\Phi_k^\epsilon$, and that $\tau$
maps the (connected) frontier of $T$ to the frontier of $\tau(T)$.
Since a free orientation-reversing involution of an oriented surface
cannot leave any boundary component invariant, we conclude that no
tight component of $\breve\Phi_k^\epsilon$ can be invariant under
$\tau$.  Thus the number of tight components must be even.
\EndProof

\Lemma
\tag{TightBegetsTight}
Let $\epsilon\in\{\pm1\}$ and let $k$ be a non-negative integer.
Let $T$ be a very tight component of
$\breve \Phi_k^\epsilon$.  Then every component of $\breve
\Phi_{k+1}^\epsilon$ which is contained in $T$ is very
tight.  Moreover if $T$ contains exactly one component
$X$ of $\breve\Phi_{k+1}^\epsilon$ then $T'$ is ambiently 
isotopic to $T$.
\EndLemma

\Proof
Since $\partial \tilde F\subset \breve \Phi_k^\epsilon$ is invariant
under the homeomorphism $\breve h_k^\epsilon$, the component
$\breve h_k^\epsilon(T)$ of
$\breve\Phi_k^{(-1)^{k+1}\epsilon}$ is very tight.
According to Proposition \xref{SmileyPropertyTwo}, 
$\breve h_k^\epsilon(\breve\Phi_{k+1}^\epsilon)$
is ambiently isotopic to the perfect intersection
$\breve\Phi_k^{(-1)^{k+1}\epsilon}\pint\breve\Phi_1^{-\epsilon}$.

Thus if $X$ is a component of $\breve \Phi_{k+1}^\epsilon$ which is
contained in $T$ then $\breve h_k^\epsilon(X)$ is ambiently
isotopic to a component $Y$ of
$\breve\Phi_k^{(-1)^{k+1}\epsilon}\pint\breve\Phi_1^{-\epsilon}$.
Note that $Y$ is isotopically contained in $\breve h_k^\epsilon(T)$.
According to the definition of $s_0$,
every tight component of $\breve\Phi_1^{-\epsilon}$ has size
at least $s_0$, and $s(T) \le s_0 $ by the definition of a very tight subsurface.
Applying Proposition \xref{TightIntersect} with 
$P=\breve\Phi_k^{(-1)^{k+1}\epsilon}$,   $P_0=\breve h_k^\epsilon(T)$, and $Q =
\breve\Phi_1^{-\epsilon}$,
we conclude that if a component of
$\breve\Phi_k^{(-1)^{k+1}\epsilon}\pint\breve\Phi_1^{-\epsilon}$
is isotopically contained in $\breve h_k^\epsilon(T)$ then it is tight of size
at most $s(T)$, and therefore is very tight.  This shows that $Y$ is very tight.

Now suppose that $T$ contains exactly one component $X$ of
$\breve\Phi_{k+1}^\epsilon$.  
Again $\breve h_k^\epsilon(X)$ is ambiently
isotopic to a component $Y$ of
$\breve\Phi_k^{(-1)^{k+1}\epsilon}\pint\breve\Phi_1^{-\epsilon}$.
Since $X = \breve\Phi_{k+1}^\epsilon\cap
T$, and since $\partial \tilde F\subset\breve\Phi_{k+1}^\epsilon$,
we have $T \cap \partial\tilde F =
X \cap \partial\tilde F$.  Since $\partial\tilde F$ is
invariant under the homeomorphism $\breve h_k^\epsilon$, we have
 $\breve h_k^\epsilon(T) \cap \partial\tilde F = h_k^\epsilon(X) \cap
\partial\tilde F = Y \cap \partial\tilde F$.  But any component of
$\breve\Phi_k^{(-1)^{k+1}\epsilon}\pint\breve\Phi_1^{-\epsilon}$ which
is isotopically contained in $\breve h_k^\epsilon(T)$ must contain some
component of $\breve h_k^\epsilon(T) \cap \partial\tilde
F\subset Y$.
Thus $Y$ is the only component of
$\breve\Phi_k^{(-1)^{k+1}\epsilon}\pint\breve\Phi_1^{-\epsilon}$ which
is isotopically contained in $\breve h_k^\epsilon(T)$.  Proposition
\xref{TightIntersect} now implies that $\breve h_k^\epsilon(T)$
is ambiently isotopic to $Y$, and hence to $\breve h_k^\epsilon(X)$.
Applying the inverse of the homeomorphism $\breve h_k^\epsilon$, we
conclude that $X$ is ambiently isotopic to $T$.
\EndProof

\Lemma
\tag{VTIncreasesWeakly}
Let $\epsilon\in\{\pm1\}$ be given and let $l\ge k \ge 0$ be integers.
Then each very tight component of $\breve\Phi_k^\epsilon$ contains at least
one very tight component of $\breve\Phi_l^\epsilon$.  In particular
$|\VT(\breve\Phi_l^\epsilon)| \ge |\VT(\breve\Phi_k^\epsilon)|$.
\EndLemma

\Proof
It suffices to consider the case $l=k+1$.  Since
$\breve\Phi_k^\epsilon \supset\breve\Phi_{k+1}^\epsilon
\supset\partial\tilde F$,
and since each component
of $\breve\Phi_k^\epsilon$ meets $\partial \tilde F$,
each component of $\breve\Phi_k^\epsilon$ must contain at least
one component of $\breve\Phi_{k+1}^\epsilon$.
The assertion therefore follows from  Lemma \xref{TightBegetsTight}.
\EndProof

\Lemma
\tag{TightPartsIsotopic}
Let $\epsilon\in\{\pm1\}$ be given and let $l\ge k$ be non-negative
integers.  If $|\VT(\breve\Phi_l^\epsilon)| =
|\VT(\breve\Phi_k^\epsilon)|$, then $\VT(\breve\Phi_l^\epsilon)$ is
isotopic to $\VT(\breve\Phi_k^\epsilon)$.
\EndLemma

\Proof
If $l > k$ and $|\VT(\breve\Phi_l^\epsilon)| =
|\VT(\breve\Phi_k^\epsilon)|$, then by Lemma \xref{VTIncreasesWeakly}
we have
$$|\VT(\breve\Phi_k^\epsilon)| = |\VT(\breve\Phi_{k+1}^\epsilon)| =
\cdots = |\VT(\breve\Phi_{l}^\epsilon)| .$$ 

It thus suffices to show that if $|\VT(\breve\Phi_{k+1}^\epsilon)| =
|\VT(\breve\Phi_k^\epsilon)|$ then $\VT(\breve\Phi_{k+1}^\epsilon)$ is
isotopic to $\VT(\breve\Phi_k^\epsilon)$.

By Lemma \xref{VTIncreasesWeakly}, if
$|\VT(\breve\Phi_{k+1}^\epsilon)| = |\VT(\breve\Phi_k^\epsilon)|$ then
each very tight component of $\breve\Phi_k^\epsilon$ contains exactly
one very tight component of $\breve\Phi_{k+1}^\epsilon$, and each very
tight component of $\breve\Phi_{k+1}^\epsilon$ is contained in a very
tight component of $\breve\Phi_k^\epsilon$.  Thus Lemma
\xref{TightBegetsTight} implies that $\VT(\breve\Phi_{k+1}^\epsilon)$ is
isotopic to $\VT(\breve\Phi_k^\epsilon)$.
\EndProof

\Lemma
\tag{TightSkipsTwo}
Let $\epsilon\in\{\pm1\}$ and $t\in\Z/m\Z$ be given, and 
let $k\ge0$ be an integer.
If $c_{\sigma_2^\epsilon(t)}$ is contained in
$\VT(\breve\Phi_{k}^\epsilon)$ then $c_t$ is contained in
$\VT(\breve\Phi_{k+2}^\epsilon)$.
\EndLemma

\Proof
Let $C$ be the component of $\breve\Phi_{k+2}^\epsilon$ which contains
$c_t$.  Then $\breve h_2^\epsilon(C)$ is the component of $\breve
h_2^\epsilon(\Phi_{k+2}^\epsilon)$ which contains
$c_{\sigma_2^\epsilon(t)}$.  Since the embedding $\breve
h_2^\epsilon$ maps $\partial\tilde F$ to $\partial\tilde F$,
the subsurface $C$ is very tight if and only if $\breve
h_2^\epsilon(C)$ is very tight.
According to Proposition \xref{SmileyPropertyTwo}, $\breve
h_2^\epsilon(\Phi_{k+2}^\epsilon)$ is isotopic to the perfect
intersection $\breve\Phi_2^{-\epsilon}\pint\breve\Phi_k^\epsilon$.  By
\xref{UniqueComponent} $\breve h_2^\epsilon(C)$ is isotopically
contained in a unique component $C'$ of
$\breve\Phi_{k+2}^\epsilon$, which must be the component of
$\breve\Phi_k^\epsilon$ that contains $c_{\sigma_2^\epsilon(t)}$.
Thus $C'$ is very tight.  By Proposition
\xref{TightIntersect}, every component of
$\breve\Phi_2^{-\epsilon}\pint\breve\Phi_k^\epsilon$ which is
contained in $C'$ is very tight.  This shows that $\breve
h_2^\epsilon(C)$ is very tight as required.
\EndProof

\Lemma
\tag{VTDecreasesStrictly}
Let $\epsilon\in\{\pm1\}$ be given and let $k > 0$ be an odd integer.
If $|\VT(\breve\Phi_{k+2}^\epsilon)| = |\VT(\breve\Phi_{k}^\epsilon)|
> 0$ then $\breve\Phi_{k+2}^\epsilon$ is ambiently isotopic to
$\breve\Phi_k^\epsilon$.
\EndLemma

\Proof
Assume that 
$|\VT(\breve\Phi_{k+2}^\epsilon)| = |\VT(\breve\Phi_{k}^\epsilon)|
> 0$.

By Lemma \xref{TightPartsIsotopic} we have that
$\VT(\breve\Phi_{k+2}^\epsilon)$ is isotopic to
$\VT(\breve\Phi_{k}^\epsilon)$.  By Lemma \xref{TightSkipsTwo} this
implies that if $c_{\sigma_2^\epsilon(t)}$ is contained in
$\VT(\breve\Phi_{k}^\epsilon)$ then so is $c_t$.  Since
$\VT(\breve\Phi_{k}^\epsilon) \not=\emptyset$ it follows that either
$c_{\bar q}\subset\VT(\breve\Phi_{k}^\epsilon)$ for every even integer
$q$, or else $c_{\bar q}\subset\VT(\breve\Phi_{k}^\epsilon)$ for every odd
integer $q$.  But $\VT(\breve\Phi_{k}^\epsilon)$ is invariant under
the free involution $\tau_k^\epsilon$, which maps each boundary curve
$c_{\bar q}$ to a boundary curve $c_{\bar r}$ where $q$ and $r$ have
opposite parity.  
Thus every boundary component of $\tilde F$ is
contained in $\VT(\breve\Phi_{k}^\epsilon)$.  Since every component of
$\breve\Phi_{k}^\epsilon$ contains a component of $\partial\tilde F$,
it follows that $\VT(\breve\Phi_{k}^\epsilon) = \breve\Phi_{k}^\epsilon$.
Since $\breve\Phi_{k+2}^\epsilon \subset \breve\Phi_{k}^\epsilon$,
and since Lemma \xref{TightBegetsTight} implies that
every component of $\breve\Phi_{k+2}^\epsilon$ which is contained in a
tight component of $\breve\Phi_{k}^\epsilon$ is tight, we conclude
that $\VT(\breve\Phi_{k+2}^\epsilon) = \breve\Phi_{k+2}^\epsilon$.
Hence $\breve\Phi_{k}^\epsilon$ is ambiently isotopic to
$\breve\Phi_{k+2}^\epsilon$.
\EndProof

\Lemma
\tag{NumVTIncreasesStrictly}
Suppose that $\tilde F$ is not a semi-fiber.
Let $\epsilon\in\{\pm1\}$ be given and let $k > 0$ be an odd integer.
If $|\VT(\breve\Phi_{k+2}^\epsilon)| = |\VT(\breve\Phi_{k}^\epsilon)| > 0$
then $\breve\Phi_k^\epsilon$ is a regular neighborhood of $\partial \tilde
F$.
\EndLemma

\Proof
This is an immediate consequence of Proposition
\xref{SmileyStrictContainment}
and Lemma \xref{VTDecreasesStrictly}.
\EndProof

\Proposition
\tag{BigSmile}
Suppose that $\tilde F$ is not a semi-fiber. Let $p>0$ be an odd
integer and suppose that either $\breve\Phi_p^+$ or $\breve\Phi_p^-$
has a tight component.  Then either 
$\breve\Phi_{p+m-2}^+$ or $\breve\Phi_{p+m-2}^-$ is a regular
neighborhood of $\partial \tilde F$.
\EndProposition

\Proof
We can assume without loss of generality that $p$ is the smallest
odd integer such that either $\breve\Phi_p^+$ or $\breve\Phi_p^-$ has a
tight component.  We claim that either $\breve\Phi_p^+$ or
$\breve\Phi_p^-$ has a very tight component.  If $p=1$ then, by the
definition of $s_0$, either $\breve\Phi_1^+$ or $\breve\Phi_1^-$ has a
component of size $s_0$, which is very tight by the definition of a
very tight component.  If $p>1$ then $s_0=\infty$ and any tight
component is very tight, so the claim is true in this case as well.

Now fix $\epsilon\in\{\pm 1\}$ such that
$\breve\Phi_p^\epsilon$ has a very tight component.  
It follows from Lemma
\xref{NumVTEven} that $|\VT(\breve\Phi_p^\epsilon)| \ge 2$.
We will show that
$\breve\Phi_{p+m-2}^\epsilon$ is a regular neighborhood of $\partial
\tilde F$.  
If there is an even integer $k$ with $0 < k < m-2$ such that
$\breve\Phi_{p+k}^\epsilon$ is a regular neighborhood of $\partial
\tilde F$ then the conclusion holds because
$\partial\tilde F \subset\breve\Phi_{p+m-2}^\epsilon
\subset\breve\Phi_{p+k}^\epsilon$.

Now suppose that there is no even integer $k$ with $0 < k < m-2$ such
that $\breve\Phi_{p+k}^\epsilon$ is a regular neighborhood of
$\partial \tilde F$.  Since $\tilde F$ is not a semi-fiber, for all
even $k$ with $0 < k < m-2$ we have that
$|\VT(\breve\Phi_{p+k+2}^\epsilon)| >
|\VT(\breve\Phi_{p+k}^\epsilon)|$ by Lemma
\xref{VTIncreasesWeakly} and Lemma
\xref{NumVTIncreasesStrictly},
and hence, by Lemma \xref{NumVTEven}, that
$|\VT(\breve\Phi_{p+k+2}^\epsilon)| \ge
|\VT(\breve\Phi_{p+k}^\epsilon)| + 2$.
Since $|\VT(\breve\Phi_p^\epsilon)| \ge 2$ it follows that
$|\VT(\breve\Phi_{p+m-2}^\epsilon)| \ge m$.

In particular $\breve\Phi_{p+m-2}$ has at least $m$ tight components.
Each tight component has size at least 1, and the sum of the sizes is
at most $m = |\partial\tilde F|$.  Thus 
$\breve\Phi_{p+m-2}$ has $m$ tight components of size
exactly 1, which are therefore regular neighborhoods of components of
$\partial\tilde F$.  Since $\breve\Phi_{p+m-2}^\epsilon$ is perfect it must be
a regular neighborhood of $\partial\tilde F$.
\EndProof

\Corollary
\tag{BigSmileCorollary}
Suppose that $\tilde F$ is not a semi-fiber. Let $p>0$ be an odd
integer and suppose that either $\breve\Phi_p^+$ or $\breve\Phi_p^-$
has a tight component.  Then both 
$\breve\Phi_{p+m-1}^+$ and $\breve\Phi_{p+m-1}^-$ are regular
neighborhoods of $\partial \tilde F$.
\EndProposition

\Proof
By Proposition \xref{BigSmile} we know that, for some
$\epsilon\in\{\pm1\}$, $\breve\Phi_{p+m-2}^\epsilon$ is a regular
neighborhood of $\partial \tilde F$.  Since the
$\breve\Phi_k^\epsilon$ are nested perfect surfaces,
$\breve\Phi_{p+m-1}^\epsilon$ is also a regular neighborhood of
$\partial \tilde F$.  By Proposition \xref{SmileyPropertyTwo} the
subsurface $\breve\Phi_{p+m-1}^{-\epsilon}$ is mapped homeomorphically
by $\breve h_1^{-\epsilon}$ to
$\breve\Phi_1^{-\epsilon}\pint\breve\Phi_{p+m-2}^\epsilon$, which is
isotopically contained in $\breve\Phi_{p+m-2}^\epsilon$.
Since $\partial\tilde F$ is invariant under $\breve h_1^{-\epsilon}$
it follows that 
$\breve\Phi_{p+m-1}^{-\epsilon}$ is a regular neighborhood of
$\partial\tilde F$.
\EndProof

\subsection{Planar essential surfaces and their boundary slopes}

\Theorem
\tag{TightLengthBound}
Let $F$ be an essential planar surface in a simple knot manifold $M$.
Suppose that $F$ is not a semi-fiber.  Set $m =
|\bdry F|$ and let $H$ be any reduced homotopy in the pair $(M,F)$ such
that $H_0$ is an essential path in $F$ and $H_t(\partial I)\subset
\partial M$ for each $t\in I$. Then the length of $H$ is at
most $m-1$.
\EndTheorem

\Proof
We shall assume that the length of $H$ is at least $m$ and derive a
contradiction.

Let $\tilde F$ be the splitting surface associated to $F$
(see \xref{AssociatedSplittingSurface}).  Set $\tilde m = |\partial
\tilde F|$, so $\tilde m = m$ if $F$ is separating and $\tilde m = 2m$
if $F$ is non-separating.  The homotopy $H$ determines a reduced
homotopy $\tilde H$ in the pair $(M,\tilde F)$ of length at least
$\tilde m$ such that $\tilde H_0$ is an essential path in $\tilde F$
and $\tilde H_t(\partial I)\subset \partial M$ for each $t\in I$.  In
particular, $\tilde F$ admits a long rectangle.  Thus the assumptions
of Subsections \xref{DottedEssHom} and \xref{Smiles} hold in our
situation and we may use the notation and apply the results from those
subsections.

Since the surface $F$ is not a semi-fiber, it follows that $\tilde F$
is also not a semi-fiber and hence, by Proposition
\xref{StrictContainment},
that there exists $\epsilon \in \{\pm 1\}$ such $\tilde F$ is not a
regular neighborhood of $\Phi_1^\epsilon$.  This implies that
$\breve\Phi_1^\epsilon$ is a proper subsurface of $\tilde F$.  Thus
Proposition \xref{HasTightComponent} implies that
$\breve\Phi_1^\epsilon$ has a tight component.  We conclude from
Proposition \xref{BigSmile} that $\breve \Phi_{\tilde m}^+$ and
$\breve \Phi_{\tilde m}^-$ are regular neighborhoods of $\partial
\tilde F$.  Hence $\dot \Phi_{\tilde m}^+$ and
$\dot \Phi_{\tilde m}^-$ are empty.

On the other hand, since the reduced homotopy $\tilde H$
has length at least $\tilde m$, it follows from \xref{ArcToGlasses}
and Proposition \xref{Surfaces} that there is an admissible pair of
glasses $\gamma:\Gamma\to\tilde F$ which is homotopic in $\tilde F$ to
a map from $\Gamma$ to $\Phi_{\tilde m}^\epsilon$.  In particular there is
a map $\alpha:S^1 \to \partial\tilde F$ which is homotopic in $\tilde
F$ to a map from $S^1$ to a component $A$ of $\Phi_{\tilde m}^\epsilon$.
It follows that $A$ must be an outer component of
$\Phi_{\tilde m}^\epsilon$, and hence that $\dot\Phi_{\tilde m}^\epsilon
\not=\emptyset$.  This contradiction completes the proof.
\EndProof

\Theorem
\tag{PlanarTheorem}
Let $M$ be a simple knot manifold and let $F\subset M$ be
an essential planar surface with boundary slope $\beta$ which is
not a semi-fiber.  Let $(S,X,h)$ be a singular
surface which is well-positioned with respect to $F$ and has boundary
slope $\alpha$. Set $s=\genus S$, $n=|\partial S-X|$, $v=|X|$. Then
$$\Delta(\alpha,\beta)\le N(s,n,v).$$
\EndTheorem

\Proof
Set $m=|\partial F|$. 
According to Proposition \xref{GraphProposition}, 
there exists an essential homotopy $H:I\times I\to M$ having
length $$l \ge {m\D(\a,\b)\over N(s,n,v)}-1. \leqno(1)$$ such that
$H_0$ is an essential path in $F$ and $H_t(\partial I)\subset
\partial M$ for all $t\in I$.
By Proposition \xref{TightLengthBound} we have that $$ l \le m - 1
. \leqno{(2)}.$$ The conclusion follows from the inequalities (1) and
(2).
\EndProof 

\Corollary
\tag{PlanarDisk}
Let $M$ be a simple knot manifold and $F\subset M$ an essential
planar surface with boundary slope $\beta$ which is not a semi-fiber.
Let $\alpha$ be a slope in $\partial M$. If
$M(\alpha)$ is very small, or more generally if $F\subset M\subset
M(\alpha)$ is not $\pi_1$-injective in $M(\alpha)$, then
$$\D(\a,\b)\le 5.$$
\EndCorollary

\Proof
We invoke Corollary \xref{Disk} to obtain a singular surface
$(S,X,h)$, well-positioned with respect to $F$, such that $\genus S=0$
and $|X|=1$. The conclusion now follows from Theorem
\xref{PlanarTheorem} because for any $v\ge1$ we have
$N(0,1,v)\le 5$.
\EndProof

\Corollary
\tag{VerySmallReducible}
Let $M$ be a simple knot manifold.  Suppose that $M(\beta)$ is a
reducible manifold which is not homeomorphic to $S^1\times S^2$
or $P^3\#P^3$ and that $M(\alpha)$ is very small. 
Then $$\D(\a,\b)\le 5.$$
\EndCorollary

\Proof
Since $M(\beta)$ is reducible, $M$ contains an essential planar
surface $F$ with boundary slope $\beta$.  If $F$ is a semi-fiber then
$M(\beta)$ is homeomorphic to either $S^1\times S^2$ or $P^3\#P^3$.
Thus the corollary follows from Corollary \xref{PlanarDisk}.
\EndProof

\Corollary
\tag{PlanarSeifert}
Let $M$ be a simple knot manifold and $F\subset M$ an essential planar
surface with boundary slope $\beta$ which is not a semi-fiber. Let
$\alpha$ be a slope in $\partial M$. If $M(\a)$ is a Seifert fibered
space or if there exists a $\pi_1$-injective map from $S^1\times S^1$
to $M$ then
$$\D(\a,\b)\le 6.$$
\EndCorollary

\Proof
We invoke Corollary \xref{Seifert} to obtain a singular surface
$(S,X,h)$, well-positioned with respect to $F$, such that either
$\genus S=0$ and $|X|=1$, or $\genus S=1$ and $|X|=0$. The conclusion
now follows from Theorem \xref{PlanarTheorem} because for any $v\ge1$ we
have $N(0,1,v)\le 5$ and $N(1,0,v)=6$.
\EndProof

\Corollary
\tag{ReducibleSeifert}
Let $M$ be a simple knot manifold.  Suppose that $M(\beta)$ is a
reducible manifold which is not homeomorphic to $S^1\times S^2$
or $P^3\#P^3$ and that $M(\alpha)$ is Seifert fibered. 
Then $$\D(\a,\b)\le 6.$$
\EndCorollary

\Proof
As in the proof of Corollary \xref{VerySmallReducible} we apply
Corollary \xref{PlanarSeifert} to the planar surface obtained by
intersecting a reducing sphere for $M(\beta)$ with $M$.
\EndProof

The following corollary to Theorem \xref{PlanarTheorem} is a special
case of a result of Gordon and Litherland [\cite{GLi}, Proposition
6.1], which has the same upper bound, but with a strict inequality and
without the assumption that the essential planar surface is not a
semi-fiber.
 
\Corollary
\tag{GordonLitherland}
Let $M$ be a simple knot manifold
and $F\subset M$ an essential planar surface with boundary slope
$\beta$ which is not a semi-fiber. Suppose that $S\subset M$ is an
essential bounded surface of genus $g$ with boundary slope
$\alpha$.  Set $m = |\partial S|$. Then we have
$$\Delta(\alpha,\beta)\le \left[{12g - 12\over m}\right]+6.$$
\EndCorollary

\Proof
We apply Proposition \xref{Cameron} to obtain
a singular surface $(S',\partial S', h)$ which is
well-positioned with respect to $F$. 
Theorem \xref{PlanarTheorem} then implies that
 $$\Delta(\alpha,\beta)\le N(g,0,m) .$$
Note that since $M$ is a simple knot manifold the surface $S'$ cannot
be a disk or an annulus.  It then follows from \xref{Nsnvf} that
$$N(s,n,v) = \left[{12g - 12\over m}\right]+6 .$$
\EndProof

%
%

\section{Seifert fibered surgeries}
\tag{SeifertSurgeries}
According to Corollary \xref{ReducibleSeifert}, if $M$ is a simple
knot manifold and if $\alpha$ and $\beta$ are slopes such that
$M(\alpha)$ is Seifert fibered while $M(\beta)$ is reducible
but is not $S^1\times S^2$ or $P^3\#P^3$,
then $\Delta(\alpha,\beta) \le 6$. In fact we know of no examples
where $\Delta(\alpha, \beta) > 3$.  In this section, building on
Corollary \xref{ReducibleSeifert}, we prove a result, Proposition
\xref{ThreeToSix}, which gives restrictions on the possible Seifert
fibrations of $M(\alpha)$ in the cases where $\Delta(\alpha,\beta) >
3$, which shows that this situation is not generic.  The proof uses
the character variety of $M$ and some observations from algebraic
number theory. We will see that in the generic situation, $3$ is an
upper bound for the distance between $\alpha$ and $\beta$. A similar
result, Proposition \xref{SixToTen} applies to the case where
$M(\alpha)$ is a Seifert fibered space that contains an incompressible
torus, $M(\alpha)$ is Seifert fibered and $\Delta(\alpha,\beta) > 5$.
Here, in place of Corollary \xref{ReducibleSeifert} we use a theorem
of Agol [\cite{Agol}] and Lackenby [\cite{Lackenby}] which implies that
$\Delta(\alpha,\beta) \le 10$ in this situation.

For any integer $n \geq 1$, we set $\zeta_n = e^{{2\pi i \over n}}$.

\Lemma
\tag{Divides}
Let $m, n \geq 1$ be integers. If $\zeta_n + \bar \zeta_n \in
\Q(\zeta_m)$, then one of the following three conditions holds.

\part{(i)} $n \in \{1,2,3,4,6\}$.

\part{(ii)} $n$ divides $m$.

\part{(iii)} ${n \over 2}$ is an odd integer dividing $m$.

\EndLemma

\Proof Without loss of generality we take
$n \not \in \{1,2,3,4,6\}$, so $\zeta_n + \bar \zeta_n \not \in \Q$.
Let $d = {\rm gcd}(n,m)$.  By hypothesis $\zeta_n + \bar \zeta_n
\in \Q(\zeta_m) \cap \Q(\zeta_n) = \Q(\zeta_d)$ ([\cite{FrohlichTaylor},
VI.2.8]) so that $\zeta_n + \bar \zeta_n \in \Q(\zeta_d)_\R$.
Thus $\Q(\zeta_n)_\R = \Q(\zeta_n + \bar \zeta_n)
\subset \Q(\zeta_d)_\R$. But clearly $\Q(\zeta_d)_\R
\subset \Q(\zeta_n)_\R$ and therefore
$\Q(\zeta_d)_\R = \Q(\zeta_n)_\R$.  Since $\zeta_n
+ \bar \zeta_n \not \in \Q$, we have $d, n > 2$. Moreover, since
$[\Q(\zeta_k + \bar \zeta_k)\colon \Q] = {\phi(k) \over 2}$
if $k > 2$ (cf. [\cite{FrohlichTaylor}, Theorem 44]), we have $\phi(n) =
\phi(d)$. Finally since $d|n$, either $d = n$ or $d$ is odd and $2d =
n$. In other words either $n$ divides $m$ or $n$ is even, ${n \over
2}$ is odd and $n$ divides $2m$.
\EndProof

In what follows we let
$$\Delta(a,b,c) = \langle x,y \; | \; x^a, y^b, (xy)^c \rangle$$
denote the $(a,b,c)$ triangle group.

\Lemma
\tag{Order}
Let $\rho\colon \Delta(a,b,c) \to PSL_2(\C)$ be a homomorphism
and suppose that the image of $\rho$ contains an element of order $n
< \infty$. Then either $n \in \{1,2,3\}$ or $n$ divides the least
common multiple of $a,b,c$.
\EndLemma

\Proof
There are matrices $A, B, C \in SL_2(\C)$ whose orders divide
$2a,2b,2c$ respectively so that $\rho(x) = \pm A, \rho(y) = \pm B,
\rho(xy) = \pm C$.  Then ${\rm trace}(A) = \zeta_{2a}^j + \bar
\zeta_{2a}^j, {\rm trace}(B) = \zeta_{2b}^k +
\bar \zeta_{2b}^k$ and ${\rm trace}(C) = \zeta_{2c}^l + \bar \zeta_{2c}^l$
for some integers $j,k,l$. Since the trace of any word in $A, B$ is an
integral polynomial in the traces of $A, B$ and $C$ [\cite{CullerShalen}, proof of
proposition 1.4.1], such a trace lies in the field
$\Q(\zeta_{2a}, \zeta_{2b}, \zeta_{2c}) = \Q(\zeta_h)$ where $h =
2{\rm lcm}(a,b,c)$ [\cite{FrohlichTaylor}, VI.2.8].

Let $W\in SL_2(\C)$ be a matrix of order $2n$ whose image $[W]$ in
$PSL(2,\C)$ is an element of order $n$ in the image of $\rho$.
Then ${\rm trace}(W) = \zeta_{2n}^m + \bar \zeta_{2n}^m$ for some
$m$ relatively prime to $2n$. Fix a word $w$ so that the
element $w(x,y) \in \Delta(a,b,c)$ satisfies $[W] = \rho(w(x,y)) =
[w(A,B)]$. Then by the previous paragraph we have $$\zeta_{2n}^m +
\bar \zeta_{2n}^m = {\rm trace}(W) \in \{\pm {\rm trace}(w(A,
B))\} \subset \Q(\zeta_h).$$ As $m$ is relatively prime to $2n$ we
have $$\zeta_{2n} + \bar \zeta_{2n} \in \Q(\zeta_{2n}^m)_\R =
\Q(\zeta_{2n}^m + \bar\zeta_{2n}^m) \subset \Q(\zeta_h).$$ Lemma
\xref{Divides} now yields the desired conclusion. \EndProof

\Lemma
\tag{LCMLemma}
Let $M$ be a simple knot manifold.
Fix slopes $\alpha$ and $\beta$ on $\partial M$. Suppose that
$M(\beta)$ is a connected sum of two lens spaces whose
fundamental groups have orders $p, q \geq 2$, and that $M(\alpha)$ is a
Seifert fibered space whose base orbifold has the form $S^2(a,b,c)$
where $a, b, c \geq 2$. If $\Delta(\alpha, \beta) > 3$, then
$\Delta(\alpha, \beta)$ divides $\lcm(a,b,c)$.
\EndLemma

\Proof
Fix a point on $\partial M$ so that we have homomorphisms
$H_1(\partial M) \cong \pi_1(\partial M) \to \pi_1(M)$. In this way
each slope $r$ on $\partial M$ determines an element $\gamma(r)$ of
$\pi_1(M)$ well-defined up to taking an inverse.

There is a curve $X_0$ contained in the $PSL_2(\C)$-character
variety of $\pi_1(M)$ containing the character of an irreducible
representation and consisting of characters $\chi_\rho$ of
representations $\rho\colon \pi_1(M) \to PSL_2(\C)$ which factor
through $\pi_1(M(\beta)) \cong \Z/p * \Z/q$ [\cite{BoyerZhang}, Example 3.2].
For each slope $r$ let $f_r\colon X_0 \to \C$ be the regular
function $f_r(\chi_\rho) = {\rm trace}(\rho(\gamma(r)))^2 -
4$. Evidently $f_{\beta}$ is identically zero. We claim that for each $r
\ne \beta$ and ideal point $x$ of $X_0$, $f_r$ has a pole at $x$. If
this were not the case, there would be a closed essential surface $S
\subset M$ which remains essential in either $M(\beta)$ or $M(\alpha)$
[\cite{BoyerZhang}, Proposition 4.10].  But $S$ compresses in both
$M(\beta)$ and $M(\alpha)$. This is obvious for $M(\beta)$, while if
$S$ is essential in $M(\alpha)$, then $S$ is a fiber in some
realization of $M(\alpha)$ as a surface bundle over the circle (see
eg. [\cite{Jaco}, VI.34]) and so it is non-separating in $M$. But then
$b_1(M) \geq 2$, contrary to the fact that $b_1(M(\beta)) = 0$. Thus
$S$ compresses in $M(\alpha)$ and therefore $f_r$ has a pole at $x$.

Let $r$ be a slope so that $\Delta(r, \beta) = 1$.
From the previous paragraph there is a character $\chi_\rho \in
X_0$ at which $f_r$ takes the value $(e^{\pi i/\Delta(\alpha,
\beta)} + e^{-\pi i/\Delta(\alpha, \beta)})^2$.  The
representation $\rho$ may be taken to factor through
$\pi_1(M(\beta))$ and to have a non-diagonalisable image (cf. the
method of proof of [\cite{CGLS}, Lemma 1.5.10]). Since $\Delta(\alpha,
\beta) > 1$, $\rho(\gamma(r))$ has order $\Delta(\alpha, \beta)$,
while by construction $\rho(\gamma(\beta)) = \pm I$. It follows
that $\rho(\gamma(\alpha)) = \pm I$ and therefore $\rho$ factors
through a representation $\rho_1\colon \pi_1(M(\alpha)) \to
PSL_2(\C)$. In fact, $\rho$ further factors through
$\pi_1^{orb}(S^2(a,b,c)) \cong \Delta(a,b,c)$. To see this, first
observe that since $b_1(M) = 1$, $\chi_\rho$ is a non-trivial
character [\cite{Boyer}, Proposition 2.8] and therefore
[\cite{AbdelghaniBoyer}, Lemma 3.1] implies that $\rho$ factors as
claimed. In conclusion we have produced a homomorphism
$\Delta(a,b,c) \to PSL_2(\C)$ which contains an element of order
$\Delta(\alpha, \beta) > 3$ in its image. Apply Lemma \xref{Order}
to see that $\Delta(\alpha, \beta)$ divides $\lcm(a,b,c)$.
\EndProof

\Proposition
\tag{ThreeToSix}
Let $M$ be a simple knot manifold and fix slopes $\alpha$ and $\beta$
on $\partial M$. Suppose that $M(\beta)$ is reducible, but not
homeomorphic to $S^1 \times S^2$ or $P^3 \# P^3$, and that $M(\alpha)$
is a Seifert fibered manifold. If $\Delta(\alpha,\beta) > 3$, then

\part{(i)} $M(\beta)$ is a connected sum of two lens spaces.

\part{(ii)} $M(\alpha)$ admits a Seifert structure whose base orbifold is the
$2$-sphere with exactly three exceptional fibers whose orders $a,b,c$
are either a Platonic or hyperbolic triple.

\part{(iii)}$\Delta(\alpha, \beta)$ is equal to 4, 5 or 6 and divides
$\lcm(a,b,c)$.
\EndProposition

\Proof
We first show that $\Delta(\alpha, \beta)\leq 1$ when $b_1(M) \geq 2$.
In this case, according to [\cite{Gabai}], the slope $\beta$ is the unique
degenerating slope for a closed non-separating essential surface $S_*$
(which is Thurston norm minimizing in the homology class it
represents) of genus larger than $1$ in $M$, i.e. $S_*$ will remain
incompressible in $M(\delta)$ for any slope $\delta$ except for
$\delta=\beta$.  Hence the irreducible manifold $M(\alpha)$ cannot be
very small.  If it is Seifert fibered we have
$\Delta(\alpha,\beta)\leq 1$ by [\cite{BoyerGordonZhang}, Proposition 5.1].

We may therefore assume that $b_1(M) = 1$.  If $M(\alpha)$
contains an embedded incompressible torus then $\Delta(\alpha,
\beta) \leq 3$ by [Oh], [Wu]. If $M(\alpha)$ is geometrically
atoroidal then it admits a Seifert structure with three or fewer
exceptional fibers and whose base orbifold has the $2$-sphere for
underlying space. When there are no more than two exceptional fibers
it is known that $\Delta(\alpha, \beta) \leq 1$ [\cite{BoyerZhang},
Theorem 1.2(1)], while if the base orbifold of $M(\alpha)$ has the
form $S^2(a,b,c)$ where $a,b,c \geq 2$ it is known that
$\Delta(\alpha, \beta) \leq 3$ if $(a,b,c)$ is a Euclidean triple
[\cite{Boyer}, Theorem C]. Thus (ii) holds.  By [\cite{CGLS}, Theorem
2.0.3] $M(\beta)$ is a connected sum of two lens spaces, so (i)
holds. Finally (iii) is a consequence of Theorem
\xref{ReducibleSeifert} and Lemma \xref{LCMLemma}.
\EndProof

\Corollary
\tag{ThreeToSixCorollary}
Let $M$ be a simple knot manifold and fix slopes $\alpha$ and $\beta$
on $\partial M$.
If $M(\beta)$ is reducible, though not $S^1 \times S^2$ or $P^3 \# P^3$,
and $M(\alpha)$ is a Seifert fibered space, then $\Delta(\alpha, \beta) \leq
5$ unless perhaps $M(\beta) \cong P^3 \# L(p,q)$ and $M(\alpha)$ is a small
Seifert manifold with base orbifold $S^2(a,b,c)$ where $(a,b,c)$ is a
hyperbolic triple and $6$ divides $\lcm(a,b,c)$.
\EndCorollary

\Proof The corollary follows from Theorem \xref{ReducibleSeifert}, the previous
proposition and the fact that $\Delta(\alpha, \beta)
\leq 5$ if $(a,b,c)$ is a Platonic triple [\cite{BoyerZhang}, Theorem 1.2(2)].
\EndProof

There is another situation when we can use the same method to sharpen
the known distance bounds.

\Lemma
Let $M$ be a simple knot manifold and fix slopes $\alpha$ and $\beta$
on $\partial M$. Suppose that $M(\beta)$ is a Seifert fibered manifold
which contains an embedded incompressible torus and that $M(\alpha)$
admits a Seifert fibration whose base orbifold has the form
$S^2(a,b,c)$ where $a, b, c \geq 2$. If $\Delta(\alpha, \beta) > 5$
then $\Delta(\alpha, \beta)$ divides $\lcm(a,b,c)$.
\EndLemma

\Proof
By [\cite{BoyerGordonZhang}, Theorems 1.1 and 1.7] we may assume that
$b_1(M) = 1$ and that $M(\beta)$ has base orbifold a Klein bottle,
$S^2(2,2,2,2)$, or $P^2(p,q)$ for some integers $p, q \geq 2$. In each
case there is a curve $X_0 \subset X(\pi_1(M))$ containing the
character of an irreducible representation and consisting of
characters of representations $\rho\colon \pi_1(M) \to PSL_2(\C)$
which factor through $\pi_1(M(\beta))$ [\cite{BoyerZhang}, Lemma 8.7].  If
for some slope $r \ne \beta$ on $\partial M$ and ideal point $x$ of
$X_0$, $f_r$ is finite at $x$, then there is a closed, essential
surface $S \subset M$ which is incompressible in at least one of
$M(\beta)$ and $M(\alpha)$ [\cite{BoyerZhang}, Proposition 4.10]. It was
shown in [\cite{BoyerZhang}, Claim, page 786] that $S$ compresses in
$M(\beta)$, so it must be essential in $M(\alpha)$. But this would
imply that $b_1(M) \geq 2$, contrary to our assumptions (cf. the proof
of Lemma \xref{LCMLemma}).  Thus for each slope $r \ne \beta$ on
$\partial M$, the function $f_r$ has a pole at each ideal point
$x$ of $X_0$.  We now proceed as in the last paragraph of the proof of
Lemma \xref{LCMLemma} to see that $\Delta(\alpha, \beta)$ divides
$\lcm(a,b,c)$.
\EndProof

\Proposition
\tag{SixToTen}
Let $M$ be a simple knot manifold and fix slopes $\alpha$ and $\beta$
on $\partial M$. Suppose that $M(\beta)$ is a Seifert fibered manifold
that contains an embedded incompressible torus, but is not
homeomorphic to the union of two twisted $I$-bundles over Klein
bottles, and that $M(\alpha)$ is a Seifert fibered manifold.  If
$\Delta(\alpha, \beta) > 5$, then

\part{(i)} $M(\beta)$ admits a Seifert fibration over $P^2$ with
exactly two exceptional fibers whose orders are $p$ and $q$ for some
integers $p > q \geq 2$.

\part{(ii)} $M(\alpha)$ admits a Seifert fibration over the
$2$-sphere with exactly three exceptional fibers whose orders $a,b,c$
form either a hyperbolic triple or the Euclidean triple $2,3,6$.

\part{(iii)}$\Delta(\alpha, \beta) \in \{6,7,8,9,10\}$ and $\Delta(\alpha, \beta)$
divides $\lcm(a,b,c)$.
\EndProposition

\Proof
By [\cite{BoyerGordonZhang}, Theorems 1.1 and 1.7] we may assume that
$b_1(M) = 1$ and that $M(\beta)$ has base orbifold of the form $K$,
the Klein bottle, $S^2(2,2,2,2)$, or $P^2(p,q)$ for some integers $p,
q \geq 2$. Since $M(\beta)$ is not the union of two twisted
$I$-bundles over the Klein bottle its base orbifold must be of the
form $P^2(p,q)$ for some integers $p > q \geq 2$. Thus (i) holds.

Next observe that if $M(\alpha)$ contains an incompressible torus,
then it is shown in [\cite{Gordon}] that $\Delta(\alpha, \beta) \leq
5$.  Thus $M(\alpha)$ admits a Seifert structure whose base orbifold
${\cal B}$ is a $2$-sphere with three or fewer cone points. Theorem
1.5 of [\cite{BoyerZhang}] shows that ${\cal B} = S^2(a,b,c)$ where
$(a,b,c)$ is a Euclidean or hyperbolic triple. The former possibility
is ruled out as in the proof of Theorem C of [\cite{Boyer}] unless
$a,b,c$ is the Euclidean triple $2,3,6$.  Finally by [\cite{Agol}] or
[\cite{Lackenby}], $\Delta(\alpha, \beta) \in \{6,7,8,9,10\}$ and the
previous lemma shows that it divides $\lcm(a,b,c)$.
\EndProof

\vfill\eject
\references

\key{Agol}
Agol, Ian,  ``Bounds on exceptional Dehn filling,'' {\it
Geom. Topol.} {\bf 4} (2000), 431--449 (electronic). 

\key{AbdelghaniBoyer}
L. Ben Abdelghani and S. Boyer, ``A calculation of the Culler-Shalen
seminorms associated to  small Seifert Dehn fillings,''
{\it Proc. Lond. Math. Soc.} {\bf 83} (2001), 235-256.

\key{Boyer}
S. Boyer, ``On the local structure of $SL_2(\C)$-character
varieties at reducible  characters,'' {\it Top. Appl.} {\bf 121}
(2002), 383-413.

\key{BoyerGordonZhang}
S. Boyer, C. McA. Gordon, and X. Zhang,
``Dehn fillings of large hyperbolic 3-manifolds'',
preprint.

\key{BoyerZhang}
S. Boyer, and X. Zhang, ``On Culler-Shalen seminorms and Dehn filling,''
{\it Ann. of Math.} (2) {\bf 148} (1998), 737--801.

\key{CGLS}
M. Culler, C. McA. Gordon, J. Luecke, and P.~B. Shalen,
``Dehn surgery on knots,'' {\it Ann. of Math.} (2)
{\bf 125} (1987), 237--300. 

\key{CaoMeyerhoff}
C. Cao and R. Meyerhoff, ``The orientable cusped hyperbolic
$3$-manifolds of minimum volume,'' {\it Invent. Math.} {\bf 146}
(2001), 451--478.

\key{CooperLong}
D. Cooper and D. Long, ``Virtually Haken Dehn-filling.'' {\it
J. Differential Geom.} {\bf 52} (1999), 173--187.

\key{CullerShalen}
M. Culler and P.~B. Shalen, ``Varieties of group representations and
splittings of $3$-manifolds,'' {\it Ann. of Math.} (2) {\bf 117}
(1983), 109--146.

\key{Epstein}
D. B. A. Epstein, ``Curves on $2$-manifolds and isotopies,'' {\it Acta
Math.} {\bf 115} (1966) 83--107.

\key{FrohlichTaylor}
A. Frolich and M. Taylor, {\it Algebraic Number Theory}\/, Cambridge
studies in advanced mathematics 27, Cambridge University Press, 1991.

\key{Gabai}
D. Gabai, "Foliations and the topology of
$3$-manifolds II,"  J. Diff. Geom. {\bf 26} (1987) 461-478.

\key{GordonSurvey}
C.~MacA. Gordon, ``Dehn filling: a survey.'' {\it Knot theory (Warsaw,
1995)\/}, Banach Center Publ., 42, Polish Acad. Sci., Warsaw, 1998,
129--144.

\key{Gordon}
C. MacA. Gordon, ``Toroidal Dehn surgeries on knots in lens spaces.''
{\it Math. Proc. Cambridge Philos. Soc.} {\bf 125} (1999), 433--440.

\key{GLi}
C. McA. Gordon and Litherland, ``Incompressible planar surfaces in
3-manifolds'' {\it Topology Appl.} {\bf 18} (1984), 121--144.

\key{GordonLuecke}
C. McA. Gordon and J. Luecke, ``Reducible manifolds and Dehn
surgery,'' {\it Topology} {\bf 35} (1996), 385--409.

\key{Jaco}
W. Jaco, {\it Lectures on three-manifold topology}, CBMS Regional
Conference Series in Mathematics, no. 43, American Mathematical
Society, Providence, R.I., 1980.

\key{Johannson}
K. Johannson, {\it Homotopy equivalences of $3$-manifolds with
boundaries}, Lecture Notes in Mathematics, no. 761.  Springer, Berlin,
1979.

\key{JacoShalen}
W. Jaco and P. B. Shalen, ``Seifert fibered spaces in $3$-manifolds,''
{\it Mem. Amer. Math. Soc.} {\bf 21} (1979), no. 220.

\key{Lackenby}
M. Lackenby, ``Hyperbolic Dehn surgery,'' {\it Invent. Math.} {\bf 140}
(2000), 243--282.
 
\key{Li}
T. Li, ``Immersed essential surfaces in hyperbolic 3-manifolds,'' {\it
Comm. Anal. Geom.} {\bf 10} (2002), 275--290.

\key{Oh}
S. Oh, ``Reducible and toroidal $3$-manifolds obtained by Dehn
fillings,'' {\it Topology Appl.} {\bf 75} (1997), 93--104.

\key{To}
I. Torisu, Boundary slopes for knots. Osaka J. Math.  33 (1996),
no. 1, 47--55.

\key{Waldhausen}
F. Waldhausen, ``On irreducible $3$-manifolds which are sufficiently
large,'' {\it Ann. of Math.} (2) {\bf 87} (1968), 56--88.

\key{Wu}
Y. Wu, ``Dehn fillings producing reducible manifolds and toroidal
manifolds,'' {\it Topology} {\bf 37} (1998), 95--108.

\bigskip
\bigskip
\address
D\'epartement de math\'ematiques
Universit\'e du Qu\'ebec \`a Montr\'eal
P. O. Box 8888, Postal Station Centre-ville
Montr\'eal, Qc, H3C 3P8, Canada
e-mail: boyer@math.uqam.ca
\endaddress
        
\address
Department of Mathematics, Statistics and Computer Science 
University of Illinois at Chicago 
851 South Morgan Street 
Chicago, IL 60607-7045 USA 
e-mail: culler@math.uic.edu
\endaddress

\address
Department of Mathematics, Statistics and Computer Science 
University of Illinois at Chicago 
851 South Morgan Street 
Chicago, IL 60607-7045 USA 
e-mail: shalen@math.uic.edu
\endaddress

\address
Department of Mathematics 
SUNY at Buffalo 
Buffalo, NY, 14260-2900, USA  
e-mail: xinzhang@math.buffalo.edu
\endaddress

\bye